\documentclass[10pt,journal,compsoc]{IEEEtran}

\usepackage{amsmath,amssymb,amsfonts,amsthm,cases}
\usepackage{booktabs, tabularx, array}
\usepackage{graphicx,subcaption}
\usepackage{color, transparent}
\usepackage[algoruled,linesnumbered]{algorithm2e}
\usepackage{extarrows}
\usepackage[noadjust]{cite}
\usepackage{url}

\newcommand*{\defeq}{\stackrel{\text{def}}{=}}
\DeclareMathOperator*{\argmin}{\arg\min}

\newtheorem{theorem}{\textbf{Theorem}}
\newtheorem{lemma}{\textbf{Lemma}}

\newtheorem{proposition}{\textbf{Proposition}}

\newtheorem{assumption}{\textbf{Assumption}}

\newcommand{\matf}[1]{\boldsymbol{#1}}
\newcommand{\vecf}[1]{\boldsymbol{#1}}
\newcommand*{\tran}{\mathsf{T}}

\begin{document}

\title{The Proxy Step-size Technique for Regularized Optimization on the Sphere Manifold}

\author{Fang~Bai, Adrien~Bartoli
	\IEEEcompsocitemizethanks{\IEEEcompsocthanksitem Both the authors are with the ENCOV, IGT, Institut Pascal, Universit\'e Clermont Auvergne, CHU Clermont-Ferrand, France. \hfil\break
		E-mail: fang.bai@yahoo.com; adrien.bartoli@gmail.com \hfil\break
		Corresponding author: Fang Bai
		\protect\\
		The work was supported by the ANR project TOPACS and the CLARA project AIALO.
	}%
	\thanks{Code available at:  {https://bitbucket.org/FangBai/proxystepsize-pgs}}
}


\IEEEtitleabstractindextext{%
\begin{abstract}
We give an effective solution to the regularized optimization problem $g (\boldsymbol{x}) + h (\boldsymbol{x})$, where $\boldsymbol{x}$ is constrained on the unit sphere $\Vert \boldsymbol{x} \Vert_2 = 1$. Here $g (\cdot)$ is a smooth cost with Lipschitz continuous gradient within the unit ball $\{\boldsymbol{x} : \Vert \boldsymbol{x} \Vert_2 \le 1 \}$ whereas $h (\cdot)$ is typically non-smooth but convex and absolutely homogeneous, \textit{e.g.,}~norm regularizers and their combinations. Our solution is based on the Riemannian proximal gradient, using an idea we call \textit{proxy step-size} -- a scalar variable which we prove is monotone with respect to the actual step-size within an interval. The proxy step-size exists ubiquitously for convex and absolutely homogeneous $h(\cdot)$, and decides the actual step-size and the tangent update in closed-form, thus the complete proximal gradient iteration. Based on these insights, we design a Riemannian proximal gradient method using the proxy step-size. We prove that our method converges to a critical point, guided by a line-search technique based on the $g(\cdot)$ cost only. The proposed method can be implemented in a couple of lines of code. We show its usefulness by applying nuclear norm, $\ell_1$ norm, and nuclear-spectral norm regularization to three classical computer vision problems. The improvements are consistent and backed by numerical experiments.
\end{abstract}
}

\maketitle

\IEEEdisplaynontitleabstractindextext

%
\IEEEpeerreviewmaketitle

\section{Introduction}

\IEEEPARstart{W}e start with optimization problems on the unit sphere, \textit{i.e.,} the sphere manifold:
\begin{equation}
\label{eq: the general form: optimization on sphere without regularization}
\argmin_{\matf{x}\, \in\, \mathbb{R}^n}\ 
g (\vecf{x})
\quad \mathrm{s.t.}\quad \Vert \vecf{x} \Vert_2 = 1.
\end{equation}
The sphere constraint $\Vert \vecf{x} \Vert_2 = 1$ has been widely used to model scale invariant mathematical structures, \textit{e.g.,}~the fundamental matrix~\cite{faugeras1995geometry} and the dual absolute quadric~\cite{triggs1997autocalibration} in geometric vision (see~\cite{Hartley2004} for more such applications).
In statistics, the likelihood is defined up to scale~\cite{etz2018introduction}, thus often results in problems in the form (\ref{eq: the general form: optimization on sphere without regularization}), \textit{e.g.,}~the spectral correspondence association~\cite{leordeanu2005spectral} and the wavelet density estimation~\cite{peter2008maximum}.
Some other well-known applications of problem (\ref{eq: the general form: optimization on sphere without regularization}) include the $p$-harmonic energy minimization~\cite{lai2014folding} and discretized Bose-Einstein condensates~\cite{hu2020brief}.

Other than being constrained on the sphere, some applications require $\vecf{x}$ to possess additional structures, \textit{e.g.,}~low-rank (by reorganizing the elements of $\vecf{x}$ in matrix form) or sparsity (\textit{i.e.,}~number of nonzeros), to favor its physical/geometric meaning. These additional properties can be often enforced by a dedicated regularization term $h (\vecf{x})$, leading to the following regularized optimization:
\begin{equation}
\label{eq: the general form: optimization problem studied}
\argmin_{\matf{x}\, \in\, \mathbb{R}^n}\ 
g (\vecf{x}) + h (\vecf{x})
\quad \mathrm{s.t.}\quad \Vert \vecf{x} \Vert_2 = 1.
\end{equation}
In general, the regularizer $h (\cdot)$ is non-smooth but convex and absolutely homogeneous.
Typically $h (\cdot)$ are norm functions or their combinations, in particular, $\ell_1$ norm for sparsity and nuclear norm for low-rank~\cite{fazel2002matrix}.

Problem (\ref{eq: the general form: optimization problem studied}) is difficult, because of the entangling of the non-smooth cost $h(\cdot)$ and the non-convex manifold constraint $\Vert \vecf{x} \Vert_2 = 1$.
This inhibits the direct applicability of well-studied classical methods,
\textit{i.e.,}~the Euclidean optimization techniques for non-smooth composite costs \cite{nesterov2003introductory, beck2017first, beck2009fast_FISTA, fazel2013hankel, nesterov2013gradient},
and the Riemannian optimization techniques for smooth costs \cite{absil2009optimization, boumal2020introduction, hu2020brief}.
In the literature, researchers have explored several ideas to solve non-smooth optimization problems with non-convex constraints,
\textit{e.g.,}
Riemannian subgradient methods \cite{grohs2016varepsilon, grohs2016nonsmooth, hosseini2017riemannian, hosseini2018line},
proximal point methods \cite{ferreira2002proximal, de2016new,  bento2017iteration},
operator-splitting methods
\cite{lai2014splitting, kovnatsky2016madmm, chen2016augmented, zhu2017nonconvex, wang2019global},
and more recently Riemannian proximal gradient methods
\cite{chen2020proximal, huang2021riemannian}
(see Section \ref{section: related work} for a short review of these methods). Among them, the Riemannian proximal gradient methods show significant advantages over other methods, in terms of both convergence guarantees and convergence speed \cite{chen2020proximal}.

The Riemannian proximal gradient is quite a recent research topic, mainly due to Chen \textit{et al.}'s~\cite{chen2020proximal} and Huang \textit{et al.}'s~\cite{huang2021riemannian} work on the Stiefel manifold.
These methods \cite{chen2020proximal, huang2021riemannian} are exact methods with convergence proofs to a critical point, while the others either lack such proofs or only have proofs for special forms of $h(\cdot)$~\cite{tan2019learning}.
At its core, \cite{chen2020proximal, huang2021riemannian} solve a non-smooth equation derived from the Karush–Kuhn–Tucker (KKT) system by an iterative semi-smooth Newton method (SSNM) \cite{xiao2018regularized}.
This approach requires the generalized Clarke differential, which is difficult to obtain and only applicable to special forms, \textit{e.g.,}~the $\ell_1$ norm regularization used in \cite{chen2020proximal, huang2021riemannian, xiao2018regularized}.
Besides, the use of the generalized Clarke differential contradicts the spirit of proximal gradient, as the proximal is intended to avoid the computation of sub-gradients.
It should be noted that the unit sphere is a special case of the Stiefel manifold, thus the result in \cite{chen2020proximal, huang2021riemannian} applies to problem (\ref{eq: the general form: optimization problem studied}).
However, due to the usage of the generalized Clarke differential, it is hard to implement $h(\cdot)$ even as the nuclear norm regularization, not to mention more advanced $h(\cdot)$.

In this work, we advance the proximal gradient method on the sphere manifold (PGS) by discovering a concept we call \textit{proxy step-size}, for convex and absolute homogeneous $h(\cdot)$.
Importantly, we prove that the proxy step-size is monotone with respect to the actual step-size in the working region and define one line-search iteration in closed-form. Thus the generalized Clarke differential is never required.
Based on this novel insight, we control the optimization flow using the proxy step-size, and establish the convergence proof to a critical point (\textit{i.e.,}~a solution that satisfies the first-order optimality condition).
Our final outputs are three PGS algorithms (of which two are accelerated algorithms with the Nesterov momentum technique) that retain the elegance of classical Euclidean proximal gradient methods. Our method is easy to implement and much faster than the SSNM based methods \cite{chen2020proximal, huang2021riemannian}. Our main contributions are highlighted as follows:

\begin{itemize}
\item\textbf{Section~\ref{subsection: proxy step-size}.}
We reveal the existence of the proxy step-size by exploiting the convexity and absolute homogeneity of the non-smooth cost $h(\cdot)$,
and show how it decides both the actual step-size and the tangent update in closed-form.

\item\textbf{Section~\ref{subsection. the mapping between proxy step-size and actual step-size}.}
We establish the monotonicity between the proxy step-size and the actual step-size, allowing one to control the actual step-size using the proxy step-size monotonically.

\item\textbf{Section~\ref{subsection. line-search and convergence proofs}.}
We establish the convergence proof to a critical point, using a line-search from the $g(\cdot)$ cost only, by mildly assuming $g(\cdot)$ has Lipschitz continuous gradient within the unit ball $\{\vecf{x} : \Vert \vecf{x} \Vert \le 1 \}$.

\item\textbf{Section~\ref{section. Acceleration using Nesterov Momentum Technique}.}
We present accelerated versions of our PGS algorithm using the Nesterov momentum technique~\cite{nesterov1983method,nesterov2013gradient}. We empirically show that the accelerated algorithms converge much faster.

\item\textbf{Section~\ref{section. applications. theoretical description}.}
We demonstrate our algorithms with three applications, by applying nuclear norm, $\ell_1$ norm and nuclear-spectral norm regularization to three well-know computer vision problems.
\end{itemize}

We start with related work in Section~\ref{section: related work},
and necessary background in Section \ref{section. backgrounds. Proximal Gradient on the Sphere Manifold}.
Then we formally present our proxy step-size technique in Section~\ref{section: proposed solution by using proxy step-size},
and the accelerated version in Section~\ref{section. Acceleration using Nesterov Momentum Technique}.
The three example applications and experimental results are given in
Section~\ref{section. applications. theoretical description}
and
Section~\ref{section. experimental results}, respectively.
Section~\ref{section. conclusion} concludes the paper.

\section{Related Work}
\label{section: related work}

Both optimizing a smooth cost on the manifold \cite{absil2009optimization, boumal2020introduction, hu2020brief} and proximal gradient for non-smooth optimization in the Euclidean space \cite{nesterov2003introductory, beck2017first, beck2017first, beck2009fast_FISTA, nesterov2013gradient} have been well studied in the literature.
However, there exist only a few methods for solving non-smooth cost functions on the manifold.
We review existing techniques in this regard.

\subsection{Riemannian Subgradient Methods}

Subgradient methods require one to evaluate the descent direction of the total cost directly.
The descent direction, in the non-smooth setting,
is characterized by the notion of generalized Clarke gradient, which is usually difficult to calculate even numerically in practice,
\textit{e.g.,} see \cite{dirr2007nonsmooth, borckmans2014riemannian} for certain types of functions.
Instead, researchers seek for approximations of the subgradient.
A key concept in this regard is the $\epsilon$-subgradient \cite{goldstein1977optimization}.
Grohs \textit{et al.} proposed two
Riemannian $\epsilon$-subgradient methods
based on the
line-search \cite{grohs2016varepsilon} and trust-region techniques \cite{grohs2016nonsmooth},
with convergence guarantee to a critical point.
Hosseini \textit{et al.} \cite{hosseini2018line} generalized the idea to the $\epsilon$-subgradient-oriented descent sequence which combines the idea of the BFGS algorithm.
Hosseini \textit{et al.} \cite{hosseini2017riemannian} gave a non-smooth Riemannian gradient sampling method with convergence analysis.
Despite the hassle to handle the subgradient, the subgradient method is shown to have a slow convergence rate $O(1/\sqrt{k})$ \cite{zhang2016first, bento2017iteration}.

\subsection{Proximal Point Methods}

Ferreira \textit{et al.} \cite{ferreira2002proximal} proposed
the proximal point algorithm
on the Riemannian manifold,
and Bento \textit{et al.} \cite{bento2017iteration} 
established the $O(1/k)$ convergence rate of the algorithm on the Hadamard manifold for convex cost functions.
However since every smooth function that is geodesically convex on a compact Riemannian manifold is a constant \cite{bishop1969manifolds},
the analysis in \cite{ferreira2002proximal, bento2017iteration} does not apply to compact Riemannian manifolds (\textit{e.g.,} the Stiefel manifold and the unit sphere).
Bento \textit{et al.} \cite{de2016new} gave
a convergence analysis on the general Riemannian setting by assuming the cost function satisfies the Kurdyka–Lojasiewicz inequality.
In terms of computation,
the Riemannian proximal point algorithm
\cite{ferreira2002proximal} requires one to solve a subproblem to which an efficient solution does not exist within the current research.
As a result, this line of research is largely restricted to theoretical interests at the moment.

\subsection{Operator-splitting Methods}

The hardness of problem (\ref{eq: the general form: optimization problem studied}) is caused by the composition of a non-smooth cost function and the non-convex manifold constraint.
In the convex setting, the cost function and the constraint can be handled separately by the alternating direction methods of multipliers (ADMM) \cite{boyd2010distributed}.
We recapitulate the major advancements of this technique in the non-smooth and non-convex setting.
Lai \textit{et al.} \cite{lai2014splitting} explored the splitting of orthogonality constraints (SOC) method, which handles the orthogonality constraint and the cost function separately.
Kovnatsky \textit{et al.} \cite{kovnatsky2016madmm} proposed
the manifold ADMM (MADMM) method,
which further exploits the composite structure of the smooth and non-smooth cost functions.
However, both these methods, SOC and MADMM, lack convergence proofs in their original paper.
A deeper exploitation in terms of the convergence study has been conducted by Wang \textit{et al.}~\cite{wang2019global}.
For special forms of $h(\cdot)$, the convergence of ADMM to a stationary point can be established, \textit{e.g.,}~see the stabilized ADMM (SADMM) for $\ell_1$ norm regularization on the Stiefel manifold \cite{tan2019learning}.
More recently, Chen \textit{et al.} \cite{chen2016augmented} proposed PAMAL,
the proximal alternating minimized
augmented Lagrangian method.
The PAMAL method enjoys the sub-sequence convergence property,
and is noticeably faster than the SOC method by the experiments in \cite{chen2016augmented}.
Another variant, termed EPALAML,
was proposed by Zhu \textit{et al.} \cite{zhu2017nonconvex},
based on the proximal alternating linearized minimization (PALM) method.
Both PAMAL \cite{chen2016augmented} and EPALAML
\cite{zhu2017nonconvex}
minimize the augmented Lagrangian function approximately with different methods.

\subsection{Proximal Gradient Methods}

The development of proximal gradient on the Riemannian manifold is a rather new topic.
It started with the landmark paper from Chen \textit{et al.} \cite{chen2020proximal} which developed a proximal gradient on the Stiefel manifold with convergence guarantee to a critical point.
Experiments in \cite{chen2020proximal} show that the Riemannian proximal gradient method is more efficient than operator-splitting methods such as SOC and PAMAL.
More recently, Huang \textit{et al.}
\cite{huang2021riemannian} proposed another formulation, by using retractions on the non-smooth cost as well to define an iteration.
In \cite{huang2021riemannian}, a convergence rate analysis was also given based on this adaptation, which is $O(1/k)$ for the case without acceleration.
In addition to convergence guarantee to a critical point, Riemannian proximal gradient methods allow the possibility to design accelerated algorithms to obtain even faster convergence, see \cite{huang2021extension, huang2021riemannian}.

However, the subproblem of each iteration in \cite{chen2020proximal, huang2021riemannian} is solved by a SSNM which is expensive. In addition, the SSNM method requires the generalized Clarke differential which is difficult for advanced $h(\cdot)$~\cite{watson1992characterization} and contradicts the spirit of proximal gradient (which aims to avoid the generalized Clarke differential).

In this work, we propose the proxy step-size technique for proximal gradient on the sphere manifold by exploiting the convexity and absolute homogeneity of $h(\cdot)$.
Our method does not require the generalized Clarke differential, thus is more applicable to advanced $h(\cdot)$.
In addition, our method is much faster thanks to the closed-form evaluation.

\section{Preliminaries}
\label{section. backgrounds. Proximal Gradient on the Sphere Manifold}


\subsection{Absolute Homogeneous Function}

A function $h(\cdot)$ is said to be \textit{absolutely homogeneous} if $h \left( \alpha \vecf{x} \right) = \left\vert \alpha \right\vert h  \left( \vecf{x} \right)$ for any scalar $\alpha$ and vector $\vecf{x}$.

\begin{lemma}
\label{lemma. properties of convex and absolutely homogeneous h()}
If $h(\cdot)$ is convex and absolutely homogeneous, then:
\begin{itemize}
\item $h(\vecf{0}) = 0$;
\item $h(\cdot)$ is even, \textit{i.e.,}~$h( - \vecf{x}) = h (\vecf{x})$ for any $\vecf{x}$.
\item  $h(\cdot)$ is non-negative, \textit{i.e.,}~$h(\vecf{x}) \ge 0$ for any $\vecf{x}$.
\end{itemize}
\end{lemma}
\begin{proof}
The first two are obtained by setting $\alpha = 0$ and $\alpha = -1$ respectively in $h \left( \alpha \vecf{x} \right) = \left\vert \alpha \right\vert h  \left( \vecf{x} \right)$.
The third is true because if $h(\cdot)$ is further convex:
\begin{equation*}
\frac{ h(\vecf{x}) + h( - \vecf{x}) }{2}  \ge  h \left( \frac{\vecf{x} + (-\vecf{x}) }{2} \right) =  h(\vecf{0})
.
\end{equation*}
The proof is immediate by $h( - \vecf{x}) = h (\vecf{x})$ and $h(\vecf{0}) = 0$.
\end{proof}

\subsection{Proximal Operator}

For a convex but possibly non-smooth function $h(\cdot)$,
the proximal operator evaluates the solution of the following convex optimization problem at a given point $\vecf{w}$:
\begin{equation}
\label{eq: definition of the proximal operator}
\begin{aligned}
\mathrm{prox}_{t h}
\left(\vecf{w}\right)
& =
\argmin_{\vecf{x}}
\left\{
t h\left(\vecf{x}\right)
+
\frac{1}{2}
\left\Vert\vecf{x} - \vecf{w}\right\Vert_2^2
\right\}
\\ & =
\argmin_{\vecf{x}}
\left\{
h\left(\vecf{x}\right)
+
\frac{1}{2t}
\left\Vert\vecf{x} - \vecf{w}\right\Vert_2^2
\right\}
\defeq
\vecf{z}
,
\end{aligned}
\end{equation}
where $t \ge 0$ is a given scalar.
For $t=0$, $\vecf{z} = \vecf{w}$.

Since the proximal is a convex problem, its solution is characterized by its first-order necessary condition:
\begin{align}
\vecf{z}
=
\mathrm{prox}_{t h} \left(\vecf{w}\right) 
& \Leftrightarrow
\vecf{0} \in \partial{h}\vert_{\vecf{z}} + \frac{1}{t}(\vecf{z} - \vecf{w})
\\ 
\label{eq: relation of proximal and subdifferential relation}
& \Leftrightarrow	
\vecf{w} \in  t \partial{h}\vert_{\vecf{z}} + \vecf{z}.
\end{align}

The proximal satisfies firm non-expansiveness and non-expansiveness (see Appendix \ref{appendix. properties of proximal}).
Additionally, we show the following statements hold true.
\begin{lemma}
\label{lemma. additional properties for convex and absolute homogeneous h}
If $h(\cdot)$ is convex and absolutely homogeneous and $t \ge 0$, we have:
\begin{itemize}
\item $\mathrm{prox}_{t h} ( \vecf{0} ) = \vecf{0}$;
\item $ \left\Vert \mathrm{prox}_{t h} ( \vecf{w} )  \right\Vert_2^2  \le  \langle  \mathrm{prox}_{t h} ( \vecf{w} ),\, \vecf{w}  \rangle $;
\item $ \left\Vert \mathrm{prox}_{t h} ( \vecf{w} )  \right\Vert_2  \le  \left\Vert \vecf{w} \right\Vert_2$
\end{itemize}
\end{lemma}
\begin{proof}
See Appendix \ref{appendix. properties of proximal}.
\end{proof}

\subsection{Proximal Gradient in the Euclidean Space}
\label{sec. Proximal Gradient in the Euclidean Space}

Let $g(\cdot)$ be smooth, and $h(\cdot)$ be convex.
A proximal gradient step in the Euclidean space is defined as:
\begin{align}
\label{eq: proximal gradient formulation in the Euclidean space}
& \vecf{v}_k
= 
\argmin_{\vecf{v} \,\in\, \mathbb{R}^n} \ 
\langle \nabla g \vert_{\vecf{x}_k}, \vecf{v}  \rangle 
+
\frac{1}{2t}
\langle  \vecf{v}, \vecf{v} \rangle
+ h(\vecf{x}_k + \vecf{v})
\\
\label{eq: proximal gradient formulation in the Euclidean space. update of x}
& \vecf{x}_{k+1} = \vecf{x}_{k} + \vecf{v}_{k}
,
\end{align}
where $t \ge 0$.
Problem (\ref{eq: proximal gradient formulation in the Euclidean space}) is convex,
whose solution is characterized by its first-order necessary condition:
\begin{align*}
& \vecf{0} \in
\nabla g \vert_{\vecf{x}_k} + \frac{1}{t} \vecf{v}_k + \partial h \vert_{\vecf{x}_k + \vecf{v}_k} \Leftrightarrow
\\ & 
\vecf{x}_k - t \nabla g \vert_{\vecf{x}_k}
\in
( \vecf{x}_k + \vecf{v}_k )
+ 
t \partial h \vert_{\vecf{x}_k + \vecf{v}_k}
.
\end{align*}
Therefore, from equation (\ref{eq: relation of proximal and subdifferential relation}), we have:
\begin{equation*}
\vecf{x}_k + 
\vecf{v}_k = \mathrm{prox}_{t h} \left(\vecf{x}_k - t \nabla g \vert_{\vecf{x}_k}\right) 
=
\vecf{x}_{k+1}
.
\end{equation*}
The update step differs from the classical gradient step for minimizing the $g(\cdot)$ cost only by a proximal operation,
thus is termed as \textit{proximal gradient}.
Here $t$ works as the \textit{step-size} in the standard gradient descent methods.

\subsection{Sphere Manifold}
The unit sphere, or the sphere manifold,
is an embedded manifold:
\begin{equation}
\mathcal{S}
=
\left\{
\vecf{x} \in \mathbb{R}^n  \,:\,
\left\Vert \vecf{x} \right\Vert_2 = 1
\right\}
.
\end{equation}
The tangent space at a point $\vecf{x} \in \mathcal{S}$ is:
\begin{equation}
\label{eq: tangent space of the unit sphere}
\mathcal{T}_{\vecf{x}} \mathcal{S} =
\left\{
\vecf{v} \in \mathbb{R}^n  \,:\,
\vecf{x}^{\tran} \vecf{v} = 0
\right\}
.
\end{equation}

Let $g\,:\, \mathcal{S} \rightarrowtail \mathbb{R}$ be a function defined on the manifold.
The Riemannian gradient of $g(\cdot)$
at $\vecf{x} \in \mathcal{S}$,
denoted by
$\mathrm{grad}\,g\vert_{\vecf{x}}$,
is the unique tangent vector satisfying:
\begin{equation*}
\mathrm{D} g(\vecf{x}) [\vecf{v}]
=
\langle
\mathrm{grad}\,g\vert_{\vecf{x}},\,
\vecf{v}
\rangle_{\vecf{x}}
,\quad
\forall
\vecf{v} \in \mathcal{T}_{\vecf{x}}\mathcal{S},
\end{equation*}
where $\mathrm{D} g(\vecf{x}) [\vecf{v}]$ is the directional derivative of $g(\cdot)$ along the direction of $\vecf{v}$.
We shall use the induced Euclidean metric as the Riemannian metric, which means
$
\langle  \vecf{v}, \vecf{v} \rangle_{\vecf{x}}
= \langle  \vecf{v}, \vecf{v} \rangle
=
\vecf{v}^{\tran} \vecf{v}
$.
Practically, the Riemannian gradient $\mathrm{grad}\,g\vert_{\vecf{x}}$ can be obtained by orthogonally projecting the Euclidean gradient
$\nabla g \vert_{\vecf{x}}$ in $\mathbb{R}^n$ into the tangent space at $\vecf{x}$.
Using the so-called orthogonal projector
$\mathrm{proj}_{\mathcal{T}_{\vecf{x}} \mathcal{S}}(\cdot)$,
we can write:
\begin{align}
\mathrm{grad}\,g\vert_{\vecf{x}} 
& =
\mathrm{proj}_{\mathcal{T}_{\vecf{x}} \mathcal{S}}
\nabla g \vert_{\vecf{x}}
=
\left(
\matf{I} -  \vecf{x}\vecf{x}^{\tran}
\right)
\nabla g \vert_{\vecf{x}}
\\ & =
\nabla g \vert_{\vecf{x}}
-
\left\langle \vecf{x}, \nabla g \vert_{\vecf{x}} \right\rangle  \vecf{x}
.
\end{align}

We shall use the following \textit{retraction} to bring
an increment in the tangent space back to the manifold:
\begin{equation}
\label{eq: retraction of unit sphere}
\mathcal{R}_{\vecf{x}} \left( \vecf{v} \right)
\defeq
\frac{\vecf{x} + \vecf{v}}{ \left\Vert \vecf{x} + \vecf{v} \right\Vert_2 }
\,:\,
\mathcal{T}_{\vecf{x}} \mathcal{S}  \rightarrowtail \mathcal{S}
.
\end{equation}
For more details of these concepts, we refer to \cite{boumal2020introduction, absil2009optimization}.

\subsection{Proximal Gradient on the Sphere Manifold}

\vspace{5pt}\noindent\textbf{Formulation in the tangent space.}
Inspired by the proximal gradient in the Euclidean space,
researchers have tried to formulate the update vector $\vecf{v}_k$ in the tangent space of $\vecf{x}_k$ as an extension to the manifold setting \cite{chen2020proximal,huang2021riemannian}.
In this work, we propose the following:
\begin{align}
\label{eq: RPG formulation, with tangent domain}
& \vecf{v}_k
= 
\argmin_{\vecf{v} \,\in\, \mathcal{T}_{\vecf{x}_k} \mathcal{S}} \ 
\langle \mathrm{grad}\,g \vert_{\vecf{x}_k},  \vecf{v}  \rangle 
  +
\frac{1}{2t}
\langle  \vecf{v}, \vecf{v} \rangle
+ h(\vecf{x}_k + \vecf{v})
\\
\label{eq: RPG formulation, retraction}
& \vecf{x}_{k+1} = \mathcal{R}_{\vecf{x}_k}\left( \vecf{v}_k \right)
.
\end{align}
The subproblem (\ref{eq: RPG formulation, with tangent domain}) was first proposed by Chen \textit{et al.}~\cite{chen2020proximal} in their work on the Stiefel manifold.
However in Chen \textit{et al.}'s work,
the update equation (\ref{eq: RPG formulation, retraction}) is $\vecf{x}_{k+1} = \mathcal{R}_{\vecf{x}_k}\left( \alpha_k   \vecf{v}_k \right)$ with $\alpha_k$ acting as another step-size to control the length of $\vecf{v}_k$ similar to Riemannian subgradient methods.
We shall see that this $\alpha_k$ is not required, as in accordance with equation (\ref{eq: proximal gradient formulation in the Euclidean space. update of x}) in the Euclidean proximal gradient.

\vspace{5pt}\noindent\textbf{The solution from the KKT system.}
Problem (\ref{eq: RPG formulation, with tangent domain}) is convex, thus its solution is uniquely characterized by the KKT system.
We take the tangent space constraint explicitly, and write the Lagrange function as:
\begin{equation*}
\mathcal{L}\left(\vecf{v},\,\mu \right) =
\langle \mathrm{grad}\,g \vert_{\vecf{x}_k},\,  \vecf{v}  \rangle 
+
\frac{1}{2t}
\langle  \vecf{v}, \vecf{v} \rangle
+ h(\vecf{x}_k + \vecf{v})
+
\mu \vecf{x}_k^{\tran} \vecf{v}
.
\end{equation*}
The KKT system is
$\vecf{0} \in \partial \mathcal{L}_{\vecf{v}}$,
$\vecf{x}_k^{\tran} \vecf{v}_k = 0$.
With some trivial calculations,
the KKT system is reduced to the following (see Appendix \ref{appendix: derivation of the two optimality equations} for details):
\begin{subnumcases} { \label{eq: old KKT system in u and t} }
\vecf{v}_k
=
\mathrm{prox}_{th}
\left(
\left(1-\mu t\right) \vecf{x}_k
- t \mathrm{grad}\,g \vert_{\vecf{x}_k}
\right)
- \vecf{x}_k
\label{eq: exact RPG iteration equation 1}
\\
\vecf{x}_k^{\tran}
\mathrm{prox}_{th}
\left(
\left(1-\mu t\right) \vecf{x}_k
- t \mathrm{grad}\,g \vert_{\vecf{x}_k}
\right) = 1 .
\label{eq: exact RPG iteration equation 2}
\end{subnumcases}
This KKT system can be solved by the SSNM method \cite{xiao2018regularized}, by solving the non-smooth equation (\ref{eq: exact RPG iteration equation 2}) first to obtain the Lagrange multiplier $\mu$
and then computing the tangent update $\vecf{v}_k$ by equation (\ref{eq: exact RPG iteration equation 1}).
This approach has been discussed in \cite{xiao2018regularized} and used by Chen \textit{et al.} in \cite{chen2020proximal}.

\vspace{5pt}\noindent\textbf{Limitations of the existing solution.}
While subproblem (\ref{eq: RPG formulation, with tangent domain}) seems a trivial extension from the Euclidean case (\ref{eq: proximal gradient formulation in the Euclidean space}),
the existing solution is not as elegant as its Euclidean counterpart.
First, to apply the SSNM method, as $h(\cdot)$ is non-smooth, the generalized Clarke differential of $h(\cdot)$ is required.
This contradicts the spirit of proximal gradient, as the proximal operator is typically used to avoid the generalized differential.
Second, the generalized Clarke differential is usually difficult to obtain and only applicable for special forms,
\textit{e.g.,}~the $\ell_1$ norm used in \cite{chen2020proximal,huang2021riemannian}.
If $h(\cdot)$ is the nuclear norm or more advanced functions, the SSNM is hard to implement.
Third, solving an inner loop by an iterative method (like SSNM) can degenerate the numerical accuracy and even the overall convergence.

In this work, we propose proxy step-size, an effective technique to handle subproblem (\ref{eq: RPG formulation, with tangent domain}). Instead of solving (\ref{eq: RPG formulation, with tangent domain}) directly, we aim to generate valid solutions to problem (\ref{eq: RPG formulation, with tangent domain}) in closed-form controlled by the proxy step-size.
With this new technique, the generalized differential is never required, and we show that proximal gradient on the sphere can be elegantly formulated as a proximal gradient step with respect to the proxy step-size followed by normalizations.

\vspace{5pt}\noindent\textbf{One iteration in Chen \textit{et al.} \cite{chen2020proximal}.}
After solving $\vecf{v}_k$,
a line-search process is used to ensure the descent of the total cost $f(\cdot) = g(\cdot) + h(\cdot)$. In general, they propose:
\begin{itemize}
\item given $t$, solve $\vecf{v}_k$ using the SSNM method;
\item set $\alpha_k \gets 1$ and shrink $\alpha_k$ until the following line-search criterion is met:
\begin{equation}
\label{eq. line-search in Chen and Huang works}
f \left( \mathcal{R}_{\vecf{x}_k} \left( \alpha \vecf{v}_k \right) \right)
\le
f \left( \vecf{x}_k \right)
- \frac{\alpha_k}{2t} \langle  \vecf{v}_k, \vecf{v}_k \rangle
;
\end{equation}
\item set $\vecf{x}_{k+1} \gets \mathcal{R}_{\vecf{x}_k} \left( \alpha_k \vecf{v}_k \right)$.
\end{itemize}
The validity of the line-search (\ref{eq. line-search in Chen and Huang works}) is proved in \cite{chen2020proximal}.

\vspace{5pt}\noindent\textbf{Convex toolbox for problem (\ref{eq: RPG formulation, with tangent domain}).}
Problem (\ref{eq: RPG formulation, with tangent domain}) is convex, thus a naive idea is to simply call a convex toolbox.
However, the computation is prohibitive for practical usage for high dimensional $\vecf{x}$.
Needless to say this is only one iteration, which we want to solve as efficiently as possible.
Mathematically, this solution is not elegant.

\section{The Proxy Step-size Technique}
\label{section: proposed solution by using proxy step-size}

We now formally present our proxy step-size technique to solve one proximal gradient iteration,
by assuming $h  \left( \cdot \right)$ to be absolutely homogeneous.
We term the proposed algorithm based on the proxy step-size technique as PGS (short for Proximal Gradient on the Sphere manifold).

\subsection{Proxy Step-size}
\label{subsection: proxy step-size}

\begin{lemma}
\label{proposition: a result of proximal operator on linear operators}
If $h(\cdot)$ is convex and absolutely homogeneous and $t \ge 0$, then
$
\mathrm{prox}_{t h} \left( \alpha \vecf{w} \right) = \alpha \mathrm{prox}_{\frac{t}{\left\vert\alpha\right\vert}h}
\left( \vecf{w} \right)
$
for any scalar $\alpha$.
\end{lemma}
\begin{proof}
	See Appendix \ref{appendix: proof of lemma on proximal for absolutely scable functions}.
\end{proof}

Now we introduce
$t' = \frac{t}{1 - \mu t}$, which we term \textit{proxy step-size}.
We show that the update equations from $\vecf{x}_k$ to $\vecf{x}_{k+1}$ are completely determined by the proxy step-size $t'$.
To that end, we rewrite
equation (\ref{eq: exact RPG iteration equation 1}) as:
\begin{align*}
\vecf{x}_k + \vecf{v}_k
& =  \mathrm{prox}_{th}
\left(
\left(1-\mu t\right) \vecf{x}_k
- t \mathrm{grad}\,g \vert_{\vecf{x}_k}
\right)
\\ & = 
\mathrm{prox}_{th}
\left(
\left(1-\mu t\right) 
\left(
\vecf{x}_k - \frac{t}{1 - \mu t}
\mathrm{grad}\,g \vert_{\vecf{x}_k}
\right)
\right)
\\ & = 
\left(1-\mu t\right)
\mathrm{prox}_{\frac{t}{\left\vert 1 - \mu t \right\vert} h}
\left(
\vecf{x}_k - \frac{t}{1 - \mu t}
\mathrm{grad}\,g \vert_{\vecf{x}_k}
\right)
\\ & =
\frac{t}{t'}
\mathrm{prox}_{\left\vert t' \right\vert h}
\left(
\vecf{x}_k - t'
\mathrm{grad}\,g \vert_{\vecf{x}_k}
\right)
,
\end{align*}
where the third equality stems from Lemma \ref{proposition: a result of proximal operator on linear operators}.
By $\vecf{x}_k^{\tran} \vecf{v}_k = 0$,
we obtain
$
t
=
\frac{t'}{\vecf{x}_k^{\tran}
	\mathrm{prox}_{\left\vert t' \right\vert h}
	\left(
	\vecf{x}_k - t'
	\mathrm{grad}\,g \vert_{\vecf{x}_k}
	\right)}$.
Therefore the KKT system (\ref{eq: old KKT system in u and t}) can be rewritten as:
\begin{subnumcases} 
{ \label{eq: new KKT system in t' and t} }
\vecf{v}_k
=
\frac{\mathrm{prox}_{\left\vert t' \right\vert h}
	\left(
	\vecf{x}_k - t'
	\mathrm{grad}\,g \vert_{\vecf{x}_k}
	\right)}{\vecf{x}_k^{\tran}
	\mathrm{prox}_{\left\vert t' \right\vert h}
	\left(
	\vecf{x}_k - t'
	\mathrm{grad}\,g \vert_{\vecf{x}_k}
	\right)}
- \vecf{x}_k
\label{eq: our exact RPG iteration equation 1}
\\
t
=
\frac{t'}{\vecf{x}_k^{\tran}
	\mathrm{prox}_{\left\vert t' \right\vert h}
	\left(
	\vecf{x}_k - t'
	\mathrm{grad}\,g \vert_{\vecf{x}_k}
	\right)}
\defeq \phi (t') 
.
\label{eq: our exact RPG iteration equation 2}
\end{subnumcases}
The new KKT system (\ref{eq: new KKT system in t' and t}) can be considered as a reparameterization of the previous KKT system (\ref{eq: old KKT system in u and t}), using $t$ and $t'$.
However, in the new KKT system, both $\vecf{v}_k$ and $t$ are completely decided by the proxy step-size $t'$.

\begin{figure}
\centering
\includegraphics[width=0.47\textwidth]{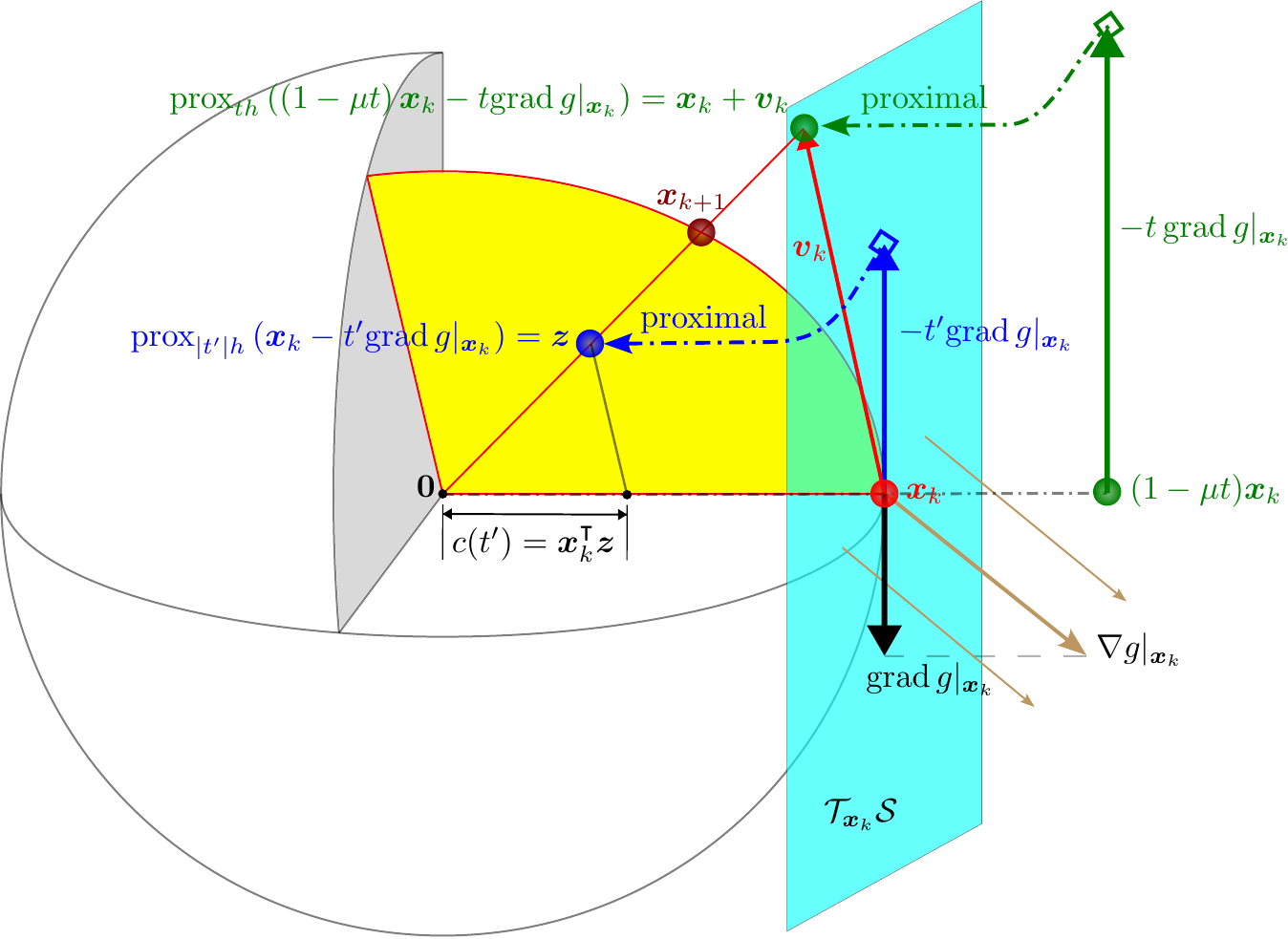}
\caption{The proxy step-size technique.
In the KKT system (\ref{eq: old KKT system in u and t}), given $t$, the essence of the non-smooth equation (\ref{eq: exact RPG iteration equation 2}) is to decide a proper Lagrange multiplier $\mu$ that lands the point $\mathrm{prox}_{th}
\left(
\left(1-\mu t\right) \vecf{x}_k
- t \mathrm{grad}\,g \vert_{\vecf{x}_k}
\right) $ into the tangent plane at $\vecf{x}_k$.
Such a process is difficult as it is hard to solve equation (\ref{eq: exact RPG iteration equation 2}) for advanced $h(\cdot)$.
Instead, we propose to use the proxy step-size $t'$ to generate valid solutions to the KKT system (\ref{eq: old KKT system in u and t}).
In specific, we first move in the tangent plane at $\vecf{x}_k$ by $ - t'\mathrm{grad}\,g\vert_{\vecf{x}_k}$, and then apply the proximal to reach the point
$\vecf{z} = \mathrm{prox}_{\left\vert t' \right\vert h}
\left(
\vecf{x}_k - t'
\mathrm{grad}\,g\vert_{\vecf{x}_k}
\right)$.
We see the point $\mathrm{prox}_{th}
\left(
\left(1-\mu t\right) \vecf{x}_k
- t \mathrm{grad}\,g \vert_{\vecf{x}_k}
\right)$ is the intersection of the line $(\vecf{0}, \vecf{z})$ and the tangent plane at $\vecf{x}_k$, given as $\frac{1}{\vecf{x}_k^{\tran} \vecf{z}} \vecf{z}$.
}
\label{fig. visualization of the proxy step-size technique}
\end{figure}

Therefore,
problem (\ref{eq: RPG formulation, with tangent domain}) can be solved in closed-form with respect to a given proxy step-size $t'$ as follow:
\begin{equation}
\label{eq: closed form solution to iteration subproblem}
\begin{cases}
\vecf{z} = \mathrm{prox}_{\left\vert t' \right\vert h}
\left(
\vecf{x}_k - t'
\mathrm{grad}\,g\vert_{\vecf{x}_k}
\right) \\
t = \frac{1}{\vecf{x}_k^{\tran}
	\vecf{z}} t' \\
\vecf{v}_k
=
\frac{1}{\vecf{x}_k^{\tran} \vecf{z}} \vecf{z}
- \vecf{x}_k 
.
\end{cases}
\end{equation}
Note that $t$ and $\vecf{v}_k$ computed from $t'$ satisfy the KKT system (\ref{eq: old KKT system in u and t}), thus they are optimal for problem (\ref{eq: RPG formulation, with tangent domain}).
An illustration of the proxy step-size technique is given in Fig.~\ref{fig. visualization of the proxy step-size technique}.

To design iterations based on the proxy step-size entirely, we need to reveal the relation between the proxy step-size $t'$ and the actual step-size $t$,
and design a line-search process to govern the convergence.

\begin{figure}[t]
\centering
\includegraphics[width = 0.42\textwidth]{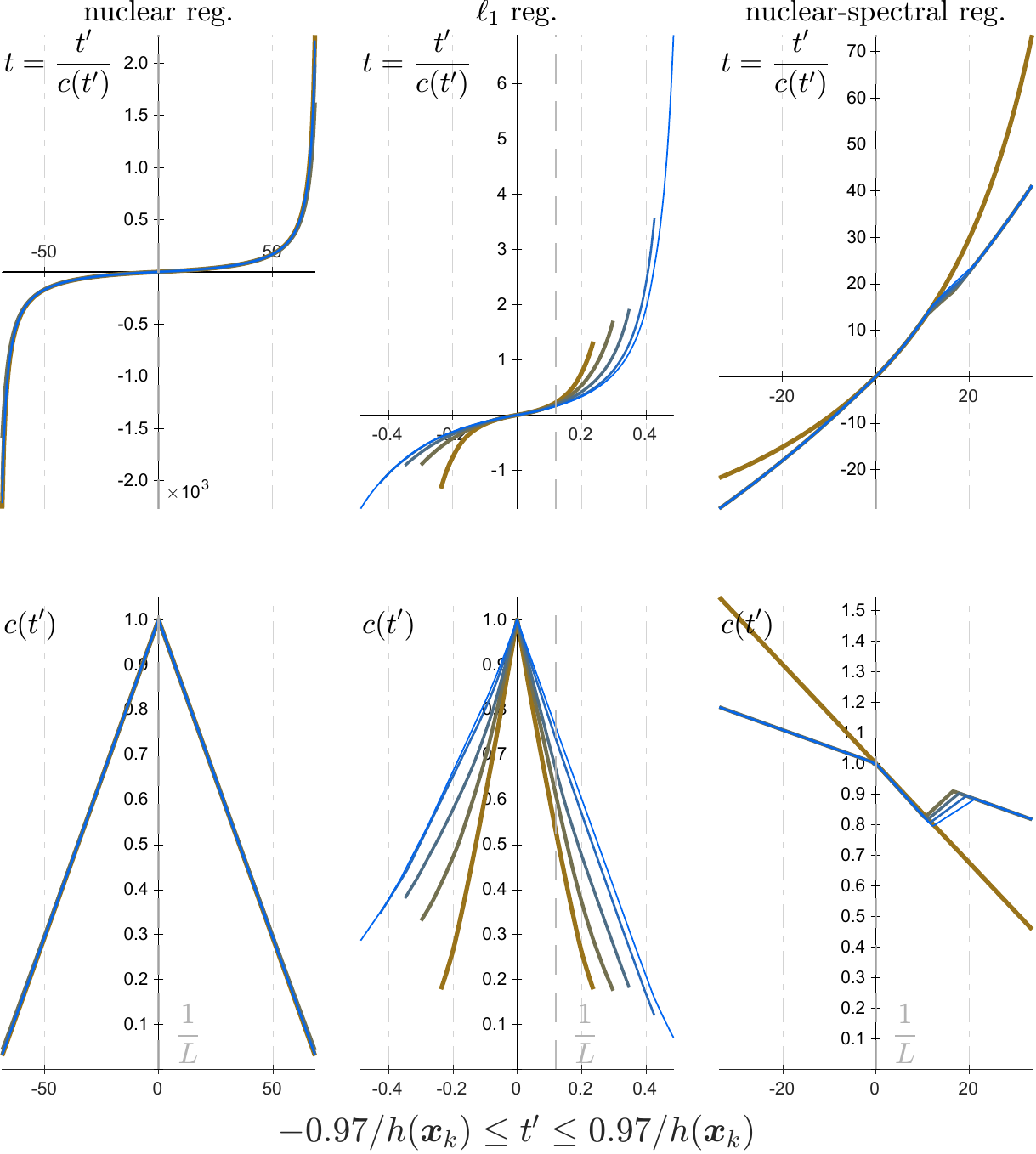}
\caption{Numerical examples for Proposition~\ref{proposition. result. c(t) satifies c(t) < 1 - th} and Theorem~\ref{theorem. the monotonicity between t and t'}.
If $h(\cdot)$ is convex and absolutely homogeneous, then for $\left\vert t' \right\vert < {1}/{h (\vecf{x}_k)}$, we have $c (t') = \vecf{x}_k^{\tran} \vecf{z} > 0$ and thus the mapping $t= \phi (t') = {t'}/c(t')$ is monotonically increasing within this interval.	
The reference proxy step-size $1/L$ to be described in Section \ref{subsection. maximum proxy step-size} is plotted as vertical dashed lines. These examples are obtained from the first $5$ iterations of:
1) fundamental matrix estimation with nuclear norm regularization,
2) correspondence association with $\ell_1$ norm regularization,
and 3) self-calibration with nuclear-spectral norm regularization to be presented in Section \ref{section. applications. theoretical description}.
}
\label{fig. numerical example for the mapping between proxy step-size and actual step-size}
\end{figure}

\subsection{The Mapping between Proxy Step-size and Actual Step-size}
\label{subsection. the mapping between proxy step-size and actual step-size}

We denote $t = \phi (t') = {t'}/{c (t')}$ the mapping between $t'$ and $t$ defined by equation (\ref{eq: our exact RPG iteration equation 2}),
where we introduce $c (t') = 
\vecf{x}_k^{\tran} \vecf{z}
=
\vecf{x}_k^{\tran} \mathrm{prox}_{\left\vert t' \right\vert h}
\left(
\vecf{x}_k - t'
\mathrm{grad}\,g \vert_{\vecf{x}_k}
\right)$.

\begin{lemma}
\label{lemma: the monotonicity of 1/t in 1/t'}
Given arbitrary $t'_1 \neq 0$, $t'_2 \neq 0$, $t'_1 \neq t'_2$, it can be shown that:
\begin{equation}
\label{eq: the monotone of reciporal t and t' proof}
\epsilon(t'_1, t'_2) = 
\left(
\frac{1}{\phi(t'_1)} - \frac{1}{\phi(t'_2)}
\right)
\left(
\frac{1}{t'_1} - \frac{1}{t'_2}
\right)
\ge 0 ,
\end{equation}
which means $1/{\phi(t')}$ is monotonically increasing with respect to ${1}/{t'}$.
\end{lemma}
\begin{proof}
See Appendix \ref{appendix: proof to theorem, relation between t and t', ineqiality elsion > 0}.
\end{proof}

Lemma \ref{lemma: the monotonicity of 1/t in 1/t'} states the monotonicity between $1/t'$ and $1/t$.
To establish the monotonicity between $t'$ and $t$ explicitly, it suffices to identify an interval where $c (t') > 0$.
\begin{proposition}
\label{proposition: the monontonicity and c>0}
If $c (t') > 0$ for $t' \in [l, u]$, then
$\phi(t')$ is monotonically increasing with respect to $t'$ for $t' \in [l, u]$.
\end{proposition}
\begin{proof}
Note that $\phi(t') = {t'}/{c (t')}$.
We thus have:
\begin{equation*}
t'_1 t'_2 \phi(t'_1) \phi(t'_2) = (t'_1 t'_2)^2/(c(t'_1) c(t'_2)) > 0 \Leftrightarrow c(t'_1) c(t'_2) > 0
\end{equation*}
Therefore if $c(t'_1) c(t'_2) > 0$, inequality (\ref{eq: the monotone of reciporal t and t' proof}) is equivalent to:
\begin{equation}
\label{eq: the monotone of t and t' proof}
\left( \phi(t'_1) - \phi(t'_2) \right)  \left( t'_1 - t'_2 \right)
\ge 0
.
\end{equation}
We see $c (t') > 0$ is sufficient for
$c(t'_1) c(t'_2) > 0$.
\end{proof}

Now we characterize an interval in which $c (t') > 0$. For this purpose, we first prove the following global inequality for convex and absolutely homogeneous $h(\cdot)$:
\begin{lemma}
	\label{eq: the general inequality about proximal of y = prox(w), for absolute homogeneous h}
	Let $h(\cdot)$ be convex and absolutely homogeneous.
	Then for any $t$, $\vecf{x}$ and $\vecf{w}$, we have:
	\begin{equation}
	\label{eq: result: key inequality <w - y, x> = t h(x), for convex and homogeneous h}
	\langle
	\vecf{w} - 
	\mathrm{prox}_{\left\vert t \right\vert h}
	\left( \vecf{w} \right)   ,\,
	\vecf{x}
	\rangle
	\le 
	\left\vert t \right\vert h (\vecf{x})
	.
	\end{equation}
	The inequality is tight for $\vecf{x} = \mathrm{prox}_{\left\vert t \right\vert h}
	\left( \vecf{w} \right)$.
\end{lemma}
\begin{proof}
See Appendix \ref{appendix. Proof of the global inequality for homogenous h}.
\end{proof}


\begin{proposition}
\label{proposition. result. c(t) satifies c(t) < 1 - th}
If $h(\cdot)$ is convex and absolutely homogeneous, then $ c (t') \ge 1 - \left\vert t' \right\vert h (\vecf{x}_k) $.
\end{proposition}
\begin{proof}
In Lemma \ref{eq: the general inequality about proximal of y = prox(w), for absolute homogeneous h},
we set $\vecf{w}= \vecf{x}_k - t' \mathrm{grad}\,g \vert_{\vecf{x}_k}$ and $\vecf{x} = \vecf{x}_k$, which yields:
\begin{multline}
\label{eq: proof: inequality realtion between x-tgrad and x_k}
\langle
\vecf{x}_k - t' \mathrm{grad}\,g \vert_{\vecf{x}_k} - 
\mathrm{prox}_{\left\vert t' \right\vert h}
\left( \vecf{x}_k - t' \mathrm{grad}\,g \vert_{\vecf{x}_k} \right)   ,\,
\vecf{x}_k
\rangle
\\ \le 
\left\vert t' \right\vert h (\vecf{x}_k)
.
\end{multline}
We reorganize inequality (\ref{eq: proof: inequality realtion between x-tgrad and x_k}) to complete the proof:
\begin{equation*}
c (t') = 
\vecf{x}_k^{\tran} \mathrm{prox}_{\left\vert t' \right\vert h}
\left(
\vecf{x}_k - t'
\mathrm{grad}\,g \vert_{\vecf{x}_k}
\right)
\ge
1 - \left\vert t' \right\vert h (\vecf{x}_k),
\end{equation*}
where we use $\langle \vecf{x}_k, \vecf{x}_k
\rangle = 1$ and $\langle \mathrm{grad}\,g \vert_{\vecf{x}_k}, \vecf{x}_k
\rangle = 0$.
\end{proof}

\begin{theorem}
\label{theorem. the monotonicity between t and t'}
If $h(\cdot)$ is convex and absolutely homogeneous, then the mapping $t= \phi (t')$ from proxy step-size $t'$ to actual step-size $t$ is monotonically increasing for $\left\vert t' \right\vert < {1}/{h (\vecf{x}_k)}$, and $\phi (0) = 0$.
\end{theorem}
\begin{proof}
By Proposition~\ref{proposition. result. c(t) satifies c(t) < 1 - th}, we see $c (t') > 0$ for $\left\vert t' \right\vert < {1}/{h (\vecf{x}_k)}$.
Then by Proposition \ref{proposition: the monontonicity and c>0},
$\phi(t')$ is monotonically increasing with respect to $t'$ for $\left\vert t' \right\vert < {1}/{h (\vecf{x}_k)}$.
By definition of the proximal operator in equation (\ref{eq: definition of the proximal operator}), we observe:
\begin{equation}
\lim_{t' \rightarrow 0} 
\mathrm{prox}_{\left\vert t' \right\vert h}
\left(
\vecf{x}_k - t'
\mathrm{grad}\,g \vert_{\vecf{x}_k}
\right) = \vecf{x}_k,
\end{equation}	
from which we conclude $c (0) = 1$,
thus $\phi(0) = 0$.
\end{proof}

Numerical examples to these results are given in Fig.~\ref{fig. numerical example for the mapping between proxy step-size and actual step-size}.
In practice, we require $t >0$, thus we use $0 <  t' < {1}/{h (\vecf{x}_k)}$.
Theorem \ref{theorem. the monotonicity between t and t'} suggests that we can control the step-size $t$ by the proxy step-size $t'$ within this interval.

\begin{proposition}
\label{proposition. result on t neq 0 for t' neq 0}
If $h(\cdot)$ is convex and absolutely homogeneous
and if $\left\Vert 
\mathrm{grad}\,g\vert_{\vecf{x}_k}
\right\Vert_2$ is bounded,
then for any $t' \neq 0$, we have $t = \phi(t') \neq 0$.
\end{proposition}
\begin{proof}
It can be shown (see Appendix~\ref{appendix. proof of results on t neq 0 for t' neq 0} for details) that:
\begin{equation*}
\vert t \vert  \ge
\frac{1}{
	\sqrt{
		\left( 1 / t' \right)^2 + 
		\left\Vert  \mathrm{grad}\,g\vert_{\vecf{x}_k}  \right\Vert_2^2
	}	
}.
\end{equation*}
Thus if $\left\Vert 
\mathrm{grad}\,g\vert_{\vecf{x}_k}
\right\Vert_2$ is bounded, then $t \neq 0$ for $t' \neq 0$.
\end{proof}

We notice that
$
\left\Vert  \mathrm{grad}\,g\vert_{\vecf{x}_k} \right\Vert_2
\le
\left\Vert  \matf{I} -  \vecf{x}_k \vecf{x}_k^{\tran} \right\Vert
\left\Vert  \nabla g \vert_{\vecf{x}_k} \right\Vert_2
=
\left\Vert \nabla g \vert_{\vecf{x}_k} \right\Vert_2
$,
thus it suffices to have a bounded $\left\Vert \nabla g \vert_{\vecf{x}_k} \right\Vert_2$.
In particular, if the Euclidean gradient $\nabla g$ is Lipschitz continuous on the sphere, $\left\Vert \nabla g \vert_{\vecf{x}_k} \right\Vert_2$ is bounded.

\subsection{Line-search and Convergence}
\label{subsection. line-search and convergence proofs}

To establish the convergence proof,
we make the following mild assumption on the cost $g(\cdot)$.
\begin{assumption}
\label{eq: Liptchtz type assumption on g cost}
We assume the pullback function $g \left( \mathcal{R}_{\vecf{x}} \left( \vecf{v} \right) \right)$
with $\vecf{v} \in \mathcal{T}_{\vecf{x}}\mathcal{S}$ satisfies:
\begin{equation}
\label{ineq. Liptchtz type g cost}
g \left( \mathcal{R}_{\vecf{x}} \left( \vecf{v} \right) \right)
\le
g (\vecf{x})
+
\langle \mathrm{grad}\,g(\vecf{x}),  \vecf{v}  \rangle 
+
\frac{L}{2}
\langle  \vecf{v}, \vecf{v} \rangle,
\end{equation}
for a (known or unknown) constant $L > 0$.	
\end{assumption}

Assumption \ref{eq: Liptchtz type assumption on g cost} holds if the ambient Euclidean space function $g (\cdot)$ on $\mathbb{R}^n$ has Lipschitz continuous gradient $\nabla g$ within the convex hull of $\mathcal{S}$ (Lemma 3 in \cite{boumal2019global}), \textit{i.e.,}~within the unit ball $\{\vecf{x} : \Vert \vecf{x} \Vert_2 \le 1 \}$.
In other words, this assumption is satisfied if $g(\cdot)$ is not changing radically within the ball (like going to infinity at some point). Thus in practice, this assumption is hardly violated.

Based on Assumption \ref{eq: Liptchtz type assumption on g cost},
we propose the following line-search.

\vspace{5pt}
\noindent\textbf{Line-search criterion.}
Search for a proxy step-size $t'$,
such that the actual step-size $t$ and the tangent update $\vecf{v}_k$ solved from equation (\ref{eq: closed form solution to iteration subproblem}) satisfy:
\begin{equation}
\label{eq: line-search equation based on t}
g \left( \mathcal{R}_{\vecf{x}_k} \left( \vecf{v}_k \right) \right)
\le
g (\vecf{x}_k)
+
\langle \mathrm{grad}\,g \vert_{\vecf{x}_k},  \vecf{v}_k \rangle 
+
\frac{1}{2t}
\langle  \vecf{v}_k,  \vecf{v}_k \rangle
.
\end{equation}
By Assumption \ref{eq: Liptchtz type assumption on g cost},
the inequality (\ref{eq: line-search equation based on t}) is satisfied for any $0 \le t \le 1/L$.
This line-search criterion is in the same spirit as the one used in the classical Euclidean proximal gradient literature, where the Euclidean gradient $\nabla g \vert_{\vecf{x}_k}$ is used instead of the Riemannian gradient here.

\vspace{5pt}
\noindent\textbf{Line-search process}.
We start with an initial proxy step-size $t' < 1/h(\vecf{x}_k)$ and reduce $t'$ until $t = \phi(t')$ satisfies  $0 \le t \le 1/L$.
This process is well-defined by the monotonicity of $t = \phi(t')$, as proved by Theorem \ref{theorem. the monotonicity between t and t'}.
One line-search iteration is given in Algorithm \ref{algorithm: our proposed algorithm for RPG on the sphere}.

\begin{algorithm}[t]
	
	\DontPrintSemicolon
	
	\SetKwFunction{lineSearchFun}{\textbf{lineSearch}}
	
	\SetKwProg{FLineSearch}{function}{}{end}
	
	\FLineSearch{$(\vecf{v}_{k},\,t,\,t') \gets$ \lineSearchFun $(\vecf{x}_k,\,t'_{\max})$}{
		
		$t' \gets \min\left\{ t'_{\max},  {1}/{h (\vecf{x}_k)} \right\}$\;

		$
		\vecf{z} \gets \mathrm{prox}_{\left\vert t' \right\vert h}
		\left(
		\vecf{x}_k - t'
		\mathrm{grad}\,g\vert_{\vecf{x}_k}
		\right)
		$
		
		$
		\vecf{v}_k
		\gets
		\frac{1}{\vecf{x}_k^{\tran}
			\vecf{z}} \vecf{z}
		- \vecf{x}_k
		$
		
		$
		t \gets \frac{1}{\vecf{x}_k^{\tran}
			\vecf{z}} t'
		$
		
		$Q_L = g (\vecf{x}_k)
		+
		\langle \mathrm{grad}\,g \vert_{\vecf{x}_k},  \vecf{v}_k \rangle 
		+
		\frac{1}{2t}
		\langle  \vecf{v}_k,  \vecf{v}_k \rangle$
		
		\lIf{$g \left( \mathcal{R}_{\vecf{x}_k} \left( \vecf{v}_k \right) \right) \le Q_L $}{											
			\Return{$(\vecf{v}_{k},\,t,\,t')$}	
		}
		\lElse{
			$t' \gets 0.8 t'$,
			goto step $3$
		}
	}
	\caption{One line-search iteration.
		\label{algorithm: our proposed algorithm for RPG on the sphere}}
\end{algorithm}

We now show the line-search criterion guarantees a descent for the total cost $f \left(\vecf{x}\right) = g \left(\vecf{x}\right)
+h \left(\vecf{x}\right)$ at each iteration.
See the PGS curve in Fig.~\ref{fig. ther convergence validation from deltaX and tangentV} for an illustration.

\begin{lemma}
\label{lemma: the h cost for given retraction and homogeinity}
For retraction (\ref{eq: retraction of unit sphere}),
if $h(\cdot)$ is absolutely homogeneous and $\vecf{v}_k \in \mathcal{T}_{\vecf{x}_k} \mathcal{S}$,
we have $h\left( \mathcal{R}_{\vecf{x}_{k}} \left( \vecf{v}_{k} \right) \right) \le h\left( \vecf{x}_{k} + \vecf{v}_{k} \right)$.
\end{lemma}
\begin{proof}
See Appendix \ref{appendix: proof of lemma h cost for given retraction and homoginity}.
\end{proof}

\begin{theorem}
	\label{theorem: the pillar reslut for convergence proofs}
	Let $f \left(\vecf{x}\right) = g \left(\vecf{x}\right)
	+h \left(\vecf{x}\right)$.
	If the line-search criterion (\ref{eq: line-search equation based on t}) holds then:
	\begin{equation}
	\label{eq: the pillar inequality for convergence proofs}
	f (\vecf{x}_{k+1})
	=
	f \left( \mathcal{R}_{\vecf{x}_k} \left( \vecf{v}_k \right) \right)
	\le
	f (\vecf{x}_k)
	-
	\frac{1}{2t}
	\langle  \vecf{v}_k, \vecf{v}_k \rangle
	.
	\end{equation}
\end{theorem}
\begin{proof}
	The first-order necessary optimality condition of problem (\ref{eq: RPG formulation, with tangent domain}) states that:
	\begin{align}
	\label{eq: the optimality condition of each iteration}
	& \vecf{0}
	\in 
	\mathrm{grad}\,g \vert_{\vecf{x}_k}
	+ 
	\frac{1}{t} \vecf{v}_k
	+
	\mathrm{proj}_{\mathcal{T}_{\vecf{x}_k} \mathcal{S}}
	\, \partial h \vert_{\vecf{x}_k +\vecf{v}_k}
	\\
	\label{eq: the optimality condition of each iteration - reorganized}
	\Leftrightarrow & - 
	\mathrm{grad}\,g \vert_{\vecf{x}_k}
	- \frac{1}{t} \vecf{v}_k
	\in
	\mathrm{proj}_{\mathcal{T}_{\vecf{x}_k} \mathcal{S}}
	\, \partial h \vert_{\vecf{x}_k +\vecf{v}_k}
	.
	\end{align}	

	Because $\vecf{v}_k \in \mathcal{T}_{\vecf{x}_k} \mathcal{S}$ then
	$
	\langle
	\vecf{z},\, \vecf{v}_k
	\rangle
	=
	\langle 
	\mathrm{proj}_{\mathcal{T}_{\vecf{x}_k} \mathcal{S}} \vecf{z},\,
	\vecf{v}_k
	\rangle
	$
	for any $\vecf{z} \in \mathbb{R}^n$.
	By the convexity of $h \left( \cdot \right)$ at $\vecf{x}_k + \vecf{v}_k$, we obtain:
	\begin{align}
	h \left( \vecf{x}_k + \vecf{v}_k \right)
	& \le
	h \left( \vecf{x}_k \right) 
	+
	\langle
	\partial h \vert_{\vecf{x}_k +\vecf{v}_k},\, \vecf{v}_k
	\rangle
	\\[4pt] &  =
	h \left( \vecf{x}_k \right) 
	+
	\langle 
	\mathrm{proj}_{\mathcal{T}_{\vecf{x}_k} \mathcal{S}} \, \partial h \vert_{\vecf{x}_k +\vecf{v}_k},\,
	\vecf{v}_k
	\rangle
	\\ & \stackrel{(\ref{eq: the optimality condition of each iteration - reorganized})}{=}
	h \left( \vecf{x}_k \right) 
	+
	\langle  - 
	\mathrm{grad}\,g \vert_{\vecf{x}_k}
	- \frac{1}{t} \vecf{v}_k,\,  \vecf{v}_k \rangle.
	\end{align}
	By considering Lemma \ref{lemma: the h cost for given retraction and homogeinity}, we have the following inequality:
	\begin{equation}
	\label{eq: intermediate result, an inequality of  iwth respect to t}
	h \left( \mathcal{R}_{\vecf{x}_k} \left( \vecf{v}_k \right) \right)
	\le
	h \left( \vecf{x}_k \right) 
	+
	\langle  - 
	\mathrm{grad}\,g \vert_{\vecf{x}_k}
	- \frac{1}{t} \vecf{v}_k,\,  \vecf{v}_k \rangle.
	\end{equation}
	Summing together inequalities (\ref{eq: line-search equation based on t}) and (\ref{eq: intermediate result, an inequality of  iwth respect to t}), we obtain inequality (\ref{eq: the pillar inequality for convergence proofs}).
\end{proof}

We denote $t'_k$ $(k=0,1,\dots,K-1)$ to be the proxy step-size satisfying the line-search criterion (\ref{eq: line-search equation based on t}) at iteration $k$.
In the line-search process, we have:
\begin{equation*}
0 < t'_k \le \min\left\{ t'_{\max},  {c}/{h (\vecf{x}_k)} \right\}
, \quad \mathrm{with\ } 0 < c < 1
.
\end{equation*}
We further denote the corresponding actual step-sizes as $t_{k}$ where $t_{k} = \phi(t'_k) = t'_k / c(t'_k)$.
From Theorem \ref{theorem. the monotonicity between t and t'}, we know $t_k \ge 0$.
By Proposition~\ref{proposition. result on t neq 0 for t' neq 0}, we know if $t'_k \neq 0$ then $t_k \neq 0$ as $\left\Vert 
\mathrm{grad}\,g\vert_{\vecf{x}_k} \right\Vert_2$ is bounded by the assumption of Lipschitz continuous gradient $\nabla g$ within the unit ball. Therefore we conclude $t_k > 0$.
To summarize, we have:
\begin{equation*}
	0 < t_{\min} \le t_k \le t_{\max},
\end{equation*}
where we denote $t_{\min} = \min \{ t_k \}$, $ t_{\max} = \max \{ t_k \}$.

At last, we show the iterations guided by the line-search process converge to a critical point of problem (\ref{eq: the general form: optimization problem studied}).

\begin{proposition}
	\label{proposition: convergence: lim v_k = 0}
	Assume $f (\vecf{x}) $ is bounded from below on $\mathcal{S}$, \textit{i.e.,} the problem is well-posed. If $g (\cdot)$ has Lipschitz continuous gradient $\nabla g$ within the unit ball $\{\vecf{x} : \Vert \vecf{x} \Vert \le 1 \}$, the line-search iterations converge to a critical point of problem (\ref{eq: the general form: optimization problem studied}).
\end{proposition}
\begin{proof}
Let $f^{\star}$ be the optimal value, \textit{i.e.,} $f^{\star} \le f(\vecf{x})$ for any $\vecf{x}$ on $\mathcal{S}$. Given $K$ iterations, inequality (\ref{eq: the pillar inequality for convergence proofs}) implies:
\begin{align*}
	f(\vecf{x}_0) - f^{\star}
	& \ge
	f(\vecf{x}_0) - f({\vecf{x}_K})
	=
	\sum_{k=0}^{K-1}
	\left(f(\vecf{x}_k) - f({\vecf{x}_{k+1}})\right)
	\\ & \ge
	\sum_{k=0}^{K-1}
	\frac{1}{2t_k}
	\langle  \vecf{v}_k,  \vecf{v}_k \rangle
	=
	\sum_{k=0}^{K-1}
	\frac{t_k}{2}
	\langle  \vecf{v}_k/{t_k},  \vecf{v}_k/{t_k} \rangle
	.
\end{align*}
Considering $0 < t_{\min} \le t_k \le t_{\max}$, we have:
\begin{align*}
& f(\vecf{x}_0) - f^{\star}  \ge \frac{1}{2 t_{\max}} \sum_{k=0}^{K-1} \left\Vert \vecf{v}_k \right\Vert_2^2
, 
\\
& f(\vecf{x}_0) - f^{\star} \ge \frac{ t_{\min}}{2} \sum_{k=0}^{K-1} \left\Vert \vecf{v}_k/{t_k} \right\Vert_2^2
.
\end{align*}
Denote $ \epsilon_0 = f(\vecf{x}_0) - f^{\star} $.
Taking $K$ to the infinity, we obtain:
\begin{align*}
\lim_{K \rightarrow \infty}\,
\sum_{k=0}^{K-1} \left\Vert \vecf{v}_k \right\Vert_2^2
\le
2 \epsilon_0 t_{\max} 
,\quad
\lim_{K \rightarrow \infty}\,
\sum_{k=0}^{K-1} \left\Vert \frac{\vecf{v}_k}{t_k} \right\Vert_2^2
\le 
\frac{2 \epsilon_0}{ t_{\min}}
.
\end{align*}
Since $f(\vecf{x})$ is bounded from below on $\mathcal{S}$,
$\epsilon_0$ is a non-negative constant.
Thus the right side of each inequality is bounded.
Noting the left side is the summation of an infinite non-negative sequence, we have
$\lim_{k \rightarrow \infty}\, \left\Vert \vecf{v}_k \right\Vert_2 = 0$, $\lim_{k \rightarrow \infty}\, \left\Vert \vecf{v}_k/{t_k} \right\Vert_2 = 0$,
which means $\lim_{k \rightarrow \infty}\,  \vecf{v}_k = \vecf{0}$, $\lim_{k \rightarrow \infty}\, \vecf{v}_k/{t_k} = \vecf{0}$.
Lastly, we notice if $\vecf{v}_k = \vecf{0}$ and $\vecf{v}_k/{t_k} = \vecf{0}$, the first-order necessary optimality condition of problem (\ref{eq: RPG formulation, with tangent domain}), \textit{i.e.,}~equation (\ref{eq: the optimality condition of each iteration}), becomes:
\begin{equation}
	\label{eq. the optimality condition of regularized optimization}
		\vecf{0}
		\in
		\mathrm{grad}\,g \vert_{\vecf{x}_k} + 
		\mathrm{proj}_{\mathcal{T}_{\vecf{x}_k} \mathcal{S}}
		\, \partial h \vert_{\vecf{x}_k}
		,
	\end{equation}
	which is exactly the first-order necessary optimality condition of problem (\ref{eq: the general form: optimization problem studied}) \cite{chen2020proximal, yang2014optimality}.
\end{proof}

Both $\left\Vert \vecf{v}_k \right\Vert$ and $\left\Vert \vecf{v}_k/{t_k} \right\Vert$ have linear convergence rates:
\begin{align*}
\min_{k=1,2,\dots,K}
\left\Vert \vecf{v}_k \right\Vert_2^2
\le
\frac{2 \epsilon_0 t_{\max}}{K} 
,
\quad
\min_{k=1,2,\dots,K}
\left\Vert  \frac{\vecf{v}_k}{t_k} \right\Vert_2^2
\le 
\frac{ \epsilon_0 }{ t_{\min} K}
,
\end{align*}
with the constant decided by line-search strategies.


\subsection{Algorithm}

The final algorithm is to repeat the line-search until convergence.
The pseudocode is provided in Algorithm \ref{algorithm: a unified implementation of PGS, A-PGS, AM-PGS methods},
which unifies the PGS method in this section and the accelerated methods (A-PGS, AM-PGS) to be discussed shortly.
Here we give additional components to complete the PGS method, before moving to its accelerated versions.

\subsubsection{Maximum proxy step-size $t'_{\max}$}
\label{subsection. maximum proxy step-size}

While the constant $L$ in Assumption \ref{eq: Liptchtz type assumption on g cost} (\textit{i.e.,}~inequality (\ref{ineq. Liptchtz type g cost})) exists ubiquitously, it is often not clear how to obtain $L$ in closed-form except for certain types of $g(\cdot)$ \textit{e.g.,}~$g(\vecf{x}) = \vecf{x}^{\tran} \matf{A} \vecf{x}$ where $L = 2 \sigma_{\max} (\matf{A})$ (see Section \ref{section. Rayleigh Quotient Optimization}).
Importantly, the existence of $L$ is used to establish proofs, but its actual value is never required to be known explicitly.
Instead, to complete Algorithm \ref{algorithm: our proposed algorithm for RPG on the sphere}, we need to decide $t'_{\max}$.

\vspace{5pt}\noindent\textbf{Setting $t'_{\max}$ from known $L$.}
In line-search criterion (\ref{eq: line-search equation based on t}),
$L$ establishes a lower-bound for the search, where if step-size $t < 1/L$ then the total cost is guaranteed to descend (Theorem \ref{theorem: the pillar reslut for convergence proofs}).
While we work with proxy step-size $t'$ and control step-size $t$ by $t = \phi(t')$ in its monotone region, the value $1/L$ provides a good reference for the maximum proxy step-size $t'_{\max}$.
Practically, if $L$ is known ahead, we recommend using $t'_{\max} = 1/L$.

\begin{figure}[t]
	\centering
	\includegraphics[width=0.46\textwidth]{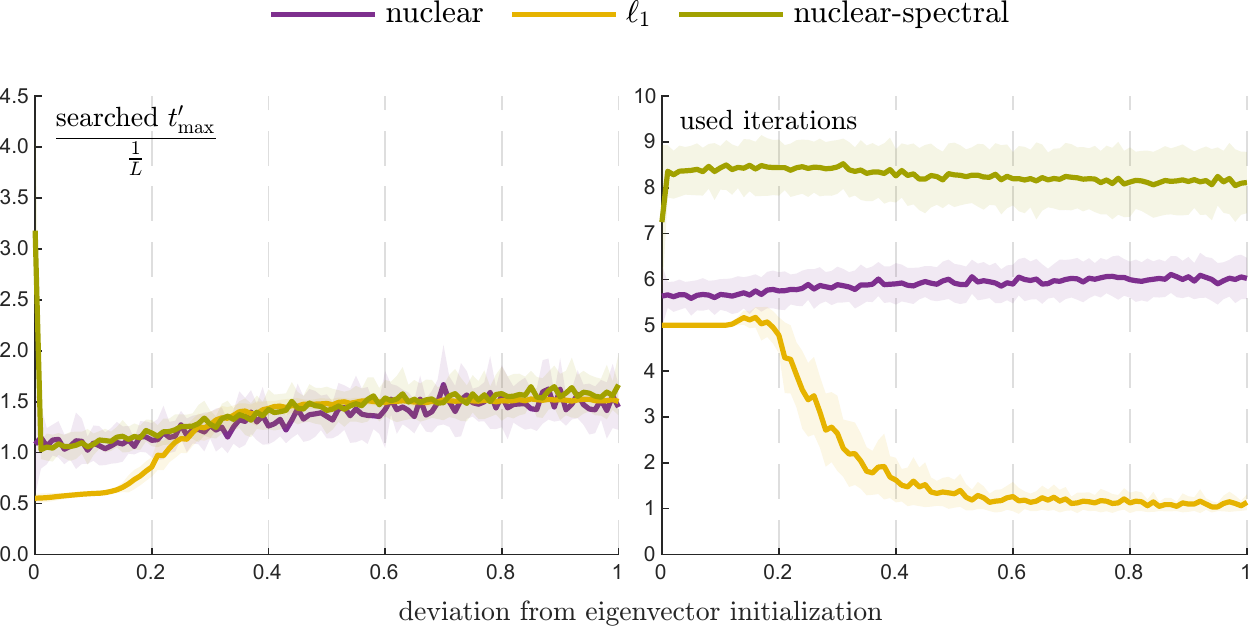}
	\caption{Examples of  line-searched maximum proxy step-size $t'_{\max}$, using different initializations $\vecf{x}_0$ (detailed in Section \ref{subsection. experiment setup of PGS methods}). On the left we report the ratio $t'_{\max}/ (1/L)$, and on the right the used iterations.}
	\label{fig. validation of adaptive search for initial proxy step-size}
\end{figure}

\vspace{5pt}\noindent\textbf{Setting $t'_{\max}$ from line-search.}
Nonetheless, if $L$ is unknown (which is the usual case), $t'_{\max}$ can be decided effectively from line-search. Such a line-search process is described in Algorithm \ref{algorithm: adaptive search for maximum proxy stepsize}, which is well-defined owing to the existence of $L$. Numerical examples are provided in Fig.~\ref{fig. validation of adaptive search for initial proxy step-size}. We compare the searched $t'_{\max}$ and the known reference proxy step-size $t' = 1/L$ using their ratio, and see that a proper $t'_{\max}$ close to $1/L$ can be found cheaply within $5$ - $10$ iterations.

\begin{algorithm}[t]
	
	\DontPrintSemicolon
	
	\SetKwFunction{searchMaxProxyStepsizeFun}{\textbf{searchMaxProxyStepsize}}
	
	\SetKwProg{FLineSearchFirst}{function}{}{end}
	
	\FLineSearchFirst{$t'_{\max} \gets$ \searchMaxProxyStepsizeFun $(\vecf{x}_0)$}{
		
		$found \gets false$, $ub \gets   {0.7}/{h (\vecf{x}_0)}$, $t' \gets   ub $\;
		
		$
		\vecf{z} \gets \mathrm{prox}_{\left\vert t' \right\vert h}
		\left(
		\vecf{x}_0 - t'
		\mathrm{grad}\,g\vert_{\vecf{x}_0}
		\right)
		$
		
		$
		\vecf{v}_0
		\gets
		\frac{1}{\vecf{x}_0^{\tran}
			\vecf{z}} \vecf{z}
		- \vecf{x}_0
		$
		
		$
		t \gets \frac{1}{\vecf{x}_0^{\tran}
			\vecf{z}} t'
		$
		
		$Q_L = g (\vecf{x}_0)
		+
		\langle \mathrm{grad}\,g \vert_{\vecf{x}_0},  \vecf{v}_0 \rangle 
		+
		\frac{1}{2t}
		\langle  \vecf{v}_0,  \vecf{v}_0 \rangle$
		
		\uIf{$g \left( \mathcal{R}_{\vecf{x}_0} \left( \vecf{v}_0 \right) \right) \le Q_L $}{
			\lIf{$t' = ub$}{ \Return{$t'_{\max} \gets ub$} }
			
			$found \gets true$,
			$t' \gets \min\left\{ 2 t',  ub \right\}$, goto step $3$
		}
		\lElseIf{$found = true$}{			
			\Return{$t'_{\max} \gets 0.5 t'$}
		}
		\lElse{
			$t' \gets 0.1 t'$, goto step $3$
		}
	}
	\caption{Search for proxy step-size $t'_{\max}$.
		\label{algorithm: adaptive search for maximum proxy stepsize}}
\end{algorithm}

\vspace{5pt}\noindent\textbf{Adaptive maximum proxy step-size $t'_{\max}$.}
We observe that typically the corresponding $t = \phi(t') = t' / c(t')$ is slightly larger than $t'$ for $t'$ around $1/L$, as $c(t')$ is slightly below $1$ as seen in Fig.~\ref{fig. numerical example for the mapping between proxy step-size and actual step-size}.
Thus setting $t' = 1/L$ does not necessarily guarantee a line-search success from criterion (\ref{eq: line-search equation based on t}) established on $t$.
On the other hand, criterion (\ref{eq: line-search equation based on t}) can also be satisfied for some $t$ greater than $1/L$.
This motivates us to use an adaptive $t'_{\max}$, where we set $t'_{\max}$ to the working $t'$ obtained at the previous iteration.
This choice is controlled by the \textit{AdaptiveMaxProxyStepsize} flag in Algorithm \ref{algorithm: a unified implementation of PGS, A-PGS, AM-PGS methods}.
We shall see such a strategy is useful to reduce subsequent total line-search iterations if $t'_{\max}$ is initially obtained from the line-search in Algorithm \ref{algorithm: adaptive search for maximum proxy stepsize} (see results in Fig.~\ref{fig. Proxy step-size with adaptive search and Lipschitz constant. Iterations}).

\subsubsection{Stop criteria for convergence}

\begin{figure}[t]
\centering
\includegraphics[width = 0.46\textwidth]{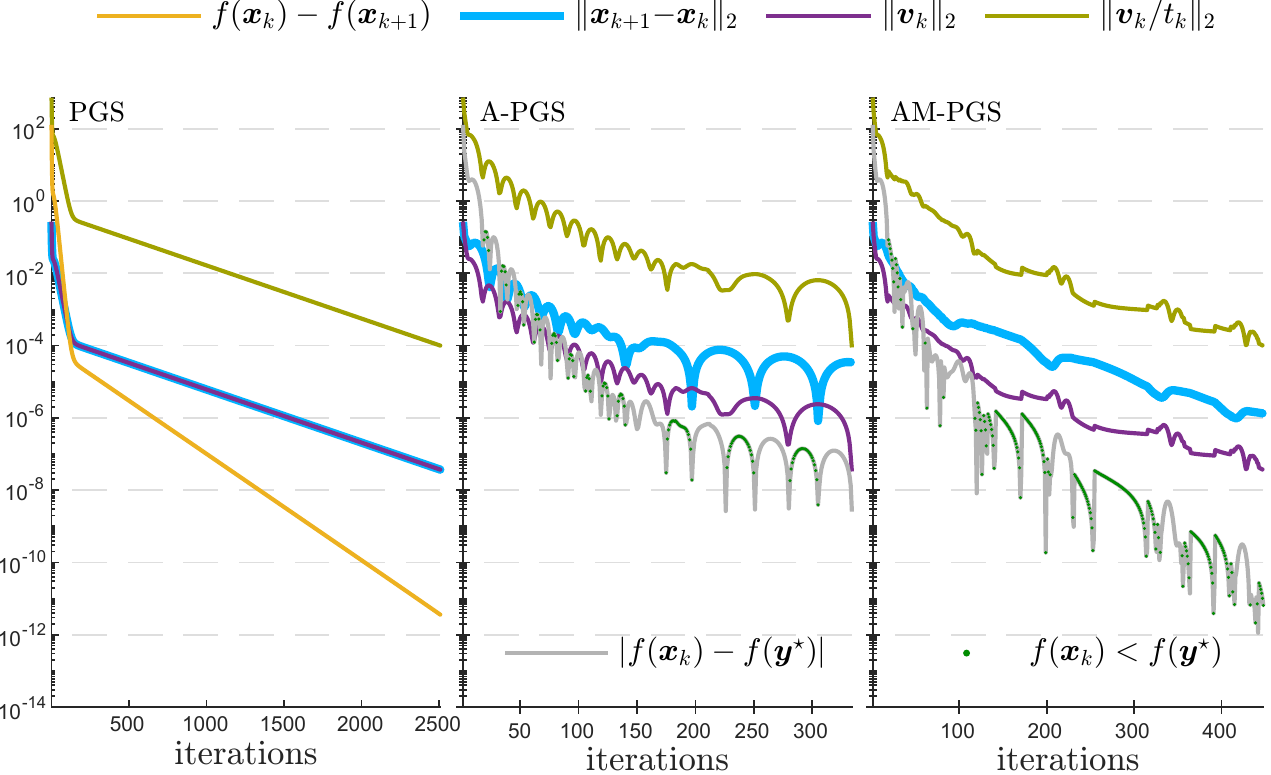}
\\[10pt]
\includegraphics[width = 0.46\textwidth]{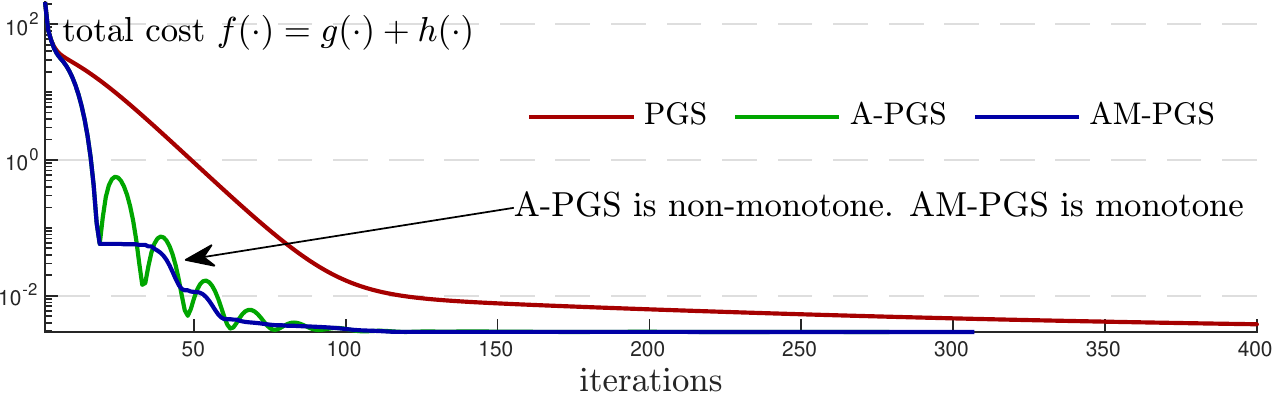}
\caption{The convergence behavior of the PGS, A-PGS and AM-PGS methods. The example is drawn from self-calibration with nuclear norm regularization.
For each method, we report the convergence of the estimate as $\left\Vert \vecf{x}_{k+1} - \vecf{x}_{k} \right\Vert_2$, of the tangent vector as $\Vert \vecf{v}_k \Vert_2$ and of the first-order optimality as $\Vert \vecf{v}_k/ t_{k} \Vert_2$.
The progression of the total cost $f(\vecf{x}_{k})$ for each method is plotted in the second figure.
In the first figure, for PGS, we report the convergence of the total cost as $f(\vecf{x}_{k}) - f(\vecf{x}_{k+1})$ which is always positive from Theorem \ref{theorem: the pillar reslut for convergence proofs}.
For A-PGS and AM-PGS, we report the absolute value $\vert f(\vecf{x}_{k}) - f(\vecf{y}^{\star}) \vert$ and mark down cases where $f(\vecf{x}_{k}) < f(\vecf{y}^{\star})$.}
\label{fig. ther convergence validation from deltaX and tangentV}
\end{figure}

In general, we propose to monitor at least $\Vert \vecf{v}_k \Vert_2$ and $\Vert \vecf{v}_k/ t_{k} \Vert_2$ as the stop criteria for convergence.
Numerical examples are provided in Fig.~\ref{fig. ther convergence validation from deltaX and tangentV}.
While other options are also possible, \textit{e.g.,}~by monitoring $f(\cdot)$, these two indicators are important for the following reasons.

\vspace{5pt}\noindent\textbf{Convergence of estimates.}
The convergence of the estimates can be determined from the distance of $\vecf{x}_{k}$ and $\vecf{x}_{k+1}$,
\textit{e.g.,}~by using the chordal distance $d = \left\Vert \vecf{x}_{k+1} - \vecf{x}_{k} \right\Vert_2$ or the angle $\theta = \mathrm{acos} \left( \vecf{x}_{k+1}^{\tran} \vecf{x}_{k} \right)$.
Here we propose to check the length of the tangent vector $\Vert \vecf{v}_k \Vert_2$.
In fact, $\vecf{v}_k  = \mathcal{R}_{\vecf{x}_{k}}^{-1} \left( \vecf{x}_{k+1} \right)$ with $\mathcal{R}_{\cdot}^{-1}(\cdot)$ to be defined in equation (\ref{eq. inverse of retraction}), thus measuring the distance of $\vecf{x}_{k}$ and $\vecf{x}_{k-1}$.

\vspace{5pt}\noindent\textbf{Optimality as critical points.}
If $\vecf{v}_k/{t_k} = \vecf{0}$ and $\vecf{v}_k = \vecf{0}$, equation (\ref{eq: the optimality condition of each iteration}) becomes equation (\ref{eq. the optimality condition of regularized optimization}), thus $\vecf{x}_k$ admits a critical point of problem (\ref{eq: the general form: optimization problem studied}) by satisfying the first-order optimality condition.
For a deeper understanding, we see from the KKT system (\ref{eq: new KKT system in t' and t}):
\begin{align*}
\frac{\vecf{v}_k}{t_k}
& = \frac{1}{t_k'}
	(\matf{I} - \vecf{x}_k \vecf{x}_k^{\tran})\,
	\mathrm{prox}_{\left\vert t_k' \right\vert h}
	\left(
	\vecf{x}_k - t_k'
	\mathrm{grad}\,g \vert_{\vecf{x}_k}
	\right)
\\ & = 
\mathrm{proj}_{\mathcal{T}_{\vecf{x}_k} \mathcal{S}} \,
\mathrm{prox}_{h}
(
\frac{1}{t_k'}
\vecf{x}_k -
\mathrm{grad}\,g \vert_{\vecf{x}_k}
)
.
\end{align*}
Thus $\vecf{v}_{k}/t_k$ is to project the proximal of $\frac{1}{t_k'} \vecf{x}_k - \mathrm{grad}\,g \vert_{\vecf{x}_k}$ into the tangent space at $\vecf{x}_k$.
A critical point is where this tangent projection goes to zero.

\subsubsection{Initialization $\vecf{x}_0$}

Typically, it is a good idea to initialize the regularized problem (\ref{eq: the general form: optimization problem studied}) from the solution of the original problem (\ref{eq: the general form: optimization on sphere without regularization}).
However this choice is problem dependent and should be discussed specifically according to the problem at hand.
For example, we initialize the regularized problems for fundamental matrix estimation and correspondence association with the solution of the Rayleigh quotient optimization (the original problem of these instances).
However, for self-calibration, the regularized problems are initialized from the canonical DAQ after quasi-calibration, as the solution space of the original problem is likely to be ambiguous due to the critical motion sequence.
Details are given in Section \ref{section. applications. theoretical description}, and numerical validations provided in Section \ref{section: experimental results: the proposed PGS algorithm}.

\section{Acceleration using the Nesterov Momentum Technique}
\label{section. Acceleration using Nesterov Momentum Technique}

The conventional PGS method evaluates the gradient and proximal at the current estimate $\vecf{x}_k$.

For the Nesterov momentum technique, the gradient and the proximal are instead evaluated at an auxiliary state $\vecf{y}_k$ defined as a linear combination of the current estimate $\vecf{x}_k$ and the previous estimate $\vecf{x}_{k-1}$ \cite{nesterov1983method, nesterov2013gradient, beck2009fast_FISTA}.
In the Euclidean case, one accelerated iteration is defined as:
\begin{equation}
\begin{cases}
\vecf{x}_{k+1} = \vecf{y}_k + \vecf{v}_k
\\[5pt]
\vecf{y}_{k+1} = \vecf{x}_{k+1} + \frac{1-\alpha_{k}}{\alpha_{k+1}} \left( \vecf{x}_{k} - \vecf{x}_{k+1} \right)
,\quad 
\vecf{y}_{0} = \vecf{x}_{0}
,
\end{cases}
\end{equation}
where $\vecf{v}_k$ is obtained by evaluating iteration (\ref{eq: proximal gradient formulation in the Euclidean space}) at $\vecf{y}_k$ as:
\begin{equation*}
\vecf{v}_k
= 
\argmin_{\vecf{v} \,\in\, \mathbb{R}^n} \ 
\langle \nabla g \vert_{\vecf{y}_k},  \vecf{v}  \rangle 
+
\frac{1}{2t}
\langle  \vecf{v}, \vecf{v} \rangle
+ h(\vecf{y}_k + \vecf{v}),
\end{equation*}
and the scalar sequence $\alpha_{k}$ $(k=1,2,\dots)$ is defined as:
\begin{equation}
\label{eq: Nestereov sequence, the alpha update}
\alpha_{k+1} = \frac{1+\sqrt{1+4 \alpha_k^2}}{2}
,\quad
\alpha_{0} = 1.
\end{equation}
The iteration defined above was first proposed by Nesterov for smooth optimization \cite{nesterov1983method}, and later extended to non-smooth composite optimization in
\cite{nesterov2013gradient, beck2009fast_FISTA}.
The original proof shows that
the accelerated iteration attains quadratic convergence for convex cost functions.
These results are also valid if the cost function is locally convex around a local minimum.
Although Riemannian manifolds are typically non-convex,
it has been shown that the Nesterov sequence can attain quadratic convergence rate for geodesically convex optimization problems on the manifold \cite{liu2017accelerated, zhang2018estimate}.

To extend the above result to the sphere manifold, we need to evaluate the difference between $\vecf{x}_{k+1}$ and $\vecf{x}_{k}$ on the sphere.
Inspired by \cite{huang2021extension},
we define this difference as a vector $\Delta \vecf{v}$ in the tangent space of $\vecf{x}_{k+1}$,
thus the first summation can be extended by the retraction at $\vecf{x}_{k+1}$.
Such $\Delta \vecf{v}$ must satisfy
$\mathcal{R}_{\vecf{x}_{k+1}}
(\Delta \vecf{v})
= \vecf{x}_{k}
$.
Abusing notations, we define the inverse of the retraction,
$
\Delta \vecf{v}
=
\mathcal{R}_{\vecf{x}_{k+1}}^{-1} \left( \vecf{x}_{k} \right)
$,
which can be calculated in closed-form as:
\begin{equation}
\label{eq. inverse of retraction}
\mathcal{R}_{\vecf{x}_{k+1}}^{-1} \left( \vecf{x}_{k} \right)
\defeq
\frac{1}{\vecf{x}_{k}^{\tran} \vecf{x}_{k+1}} \vecf{x}_{k}
- \vecf{x}_{k+1}
.
\end{equation}
It can be easily verified that
$ \mathcal{R}_{\vecf{x}_{k+1}}
(
\mathcal{R}_{\vecf{x}_{k+1}}^{-1} (\vecf{x}_{k})
)
=
\vecf{x}_k $
and
$
\mathcal{R}_{\vecf{x}_{k+1}}^{-1} (\vecf{x}_{k})
\in
\mathcal{T}_{\vecf{x}_{x+1}} \mathcal{S}
$
for any $\vecf{x}_{k},\, \vecf{x}_{k+1} \in \mathcal{S}$.

Overall we extend the Nesterov sequence to the sphere manifold as:
\begin{equation}
\label{eq: APGS iterations}
\begin{cases}
\vecf{x}_{k+1} = \mathcal{R}_{\vecf{y}_k}\left( \vecf{v}_k \right) \\[5pt]
\vecf{y}_{k+1} = \mathcal{R}_{\vecf{x}_{k+1}} \left(
\frac{1-\alpha_{k}}{\alpha_{k+1}}
\mathcal{R}_{\vecf{x}_{k+1}}^{-1} \left( \vecf{x}_{k} \right)
\right) 
,\quad 
\vecf{y}_{0} = \vecf{x}_{0}
.
\end{cases}
\end{equation}
The extension is to replace Euclidean addition and subtraction with Riemannian retraction $\mathcal{R}_{\cdot}(\cdot)$ and its inverse $\mathcal{R}^{-1}_{\cdot}(\cdot)$.
Here $\vecf{x}_k$ is obtained by evaluating iteration (\ref{eq: RPG formulation, with tangent domain}) at $\vecf{y}_k$ as:
\begin{equation*}
\vecf{v}_k
= 
\argmin_{\vecf{v} \,\in\, \mathcal{T}_{\vecf{y}_k} \mathcal{S}} \ 
\langle \mathrm{grad}\,g \vert_{\vecf{y}_k},  \vecf{v}  \rangle 
+
\frac{1}{2t}
\langle  \vecf{v}, \vecf{v} \rangle
+ h(\vecf{y}_k + \vecf{v})
,
\end{equation*}
and the scalar $\alpha_k$ is defined as in equation (\ref{eq: Nestereov sequence, the alpha update}).
The tangent update $\vecf{v}_k$ can be solved in closed-form using the proxy step-size technique proposed in Section \ref{section: proposed solution by using proxy step-size}.

The estimates $\vecf{x}_{k}$ $(k=1,2,\dots)$ generated from equation (\ref{eq: APGS iterations})
do not guarantee the monotonicity of the total cost $f(\vecf{x}_k) = g(\vecf{x}_k) + h(\vecf{x}_k)$ \cite{beck2009fast_MFISTA}.
On the manifold setting,
this can sometimes lead to divergence if $\vecf{v}_k$ is computed based on the SSNM method \cite{huang2021extension}.
To detect and recover from potential failures,
the authors in \cite{huang2021extension} introduced a safeguard by monitoring the progression of the cost function within several iterations.
One potential reason for the divergence is the inexact computation of each iteration \cite{beck2009fast_MFISTA}:
\begin{quote}
``In our case, where the denoising subproblems are not solved exactly, monotonicity becomes an important issue. It might happen that due to the inexact computations of the denoising subproblems, the algorithm might become extremely non-monotone and in fact can even diverge!''
\end{quote}
Practically in our experiments to be presented in Section \ref{section. experimental results}, we did not observe the divergence of iterations by using the proxy step-size to obtain the tangent update $\vecf{v}_k$.
This may be due to the fact that we solve each iteration exactly (in closed-form) while the SSNM method used in \cite{chen2020proximal, huang2021extension,huang2021riemannian} is iterative thus incurring inexact solutions.

Nonetheless, we propose a monotone algorithm for the sphere manifold based on Beck \textit{et al.}'s \cite{beck2009fast_MFISTA} Euclidean version
which has been proved to retain the quadratic convergence rate.
Beck \textit{et al.}'s \cite{beck2009fast_MFISTA} monotone algorithm is defined as:
\begin{equation*}
\begin{cases}
\vecf{y}^{\star} = \vecf{y}_k + \vecf{v}_k
\\[5pt]
\vecf{x}_{k+1} = 
\begin{cases}
\vecf{y}^{\star} & \quad\mathrm{if\ }  f \left( \vecf{y}^{\star} \right) < f (\vecf{x}_k) \\
\vecf{x}_{k}  & \quad\mathrm{otherwise}
\end{cases}
\\[15pt]
\vecf{y}_{k+1} = 
\begin{cases}
\vecf{y}^{\star} + \frac{1-\alpha_{k}}{\alpha_{k+1}} \left( \vecf{x}_{k} - \vecf{y}^{\star} \right) & \quad\mathrm{if\ }  f \left( \vecf{y}^{\star} \right) < f (\vecf{x}_k) \\
\vecf{x}_{k} + \frac{\alpha_{k}}{\alpha_{k+1}} \left( \vecf{y}^{\star} - \vecf{x}_{k} \right) &
\quad\mathrm{otherwise} .\\
\end{cases}
\end{cases}
\end{equation*}
The above iteration ensures the monotonicity of the total cost $f(\cdot)$ by leveraging between the new estimate $\vecf{y}^{\star}$ and the previous estimate $\vecf{x}_{k}$.
We extend Beck \textit{et al.}'s monotone algorithm to the sphere manifold as follow:
\begin{equation*}
\begin{cases}
 \vecf{y}^{\star} = \mathcal{R}_{\vecf{y}_k}\left( \vecf{v}_k \right) 
\\[5pt]
 \vecf{x}_{k+1} = 
\begin{cases}
\vecf{y}^{\star} & \quad\mathrm{if\ }  f \left( \vecf{y}^{\star} \right) < f (\vecf{x}_k) \\
\vecf{x}_{k}  & \quad\mathrm{otherwise} 
\end{cases}
\\[15pt]
 \vecf{y}_{k+1} = 
\begin{cases}
\mathcal{R}_{\vecf{y}^{\star}} \left(
\frac{1-\alpha_{k}}{\alpha_{k+1}}
\mathcal{R}^{-1}_{\vecf{y}^{\star}} \left( \vecf{x}_k \right)
\right)
& \quad\mathrm{if\ }  f \left( \vecf{y}^{\star} \right) < f (\vecf{x}_k) \\[5pt]
\mathcal{R}_{\vecf{x}_{k}} \left(
\frac{\alpha_{k}}{\alpha_{k+1}}
\mathcal{R}^{-1}_{\vecf{x}_k} \left(
\vecf{y}^{\star}
\right)
\right)
& \quad\mathrm{otherwise} .
\end{cases}
\end{cases}
\end{equation*}

The accelerated PGS (A-PGS), \textit{i.e.,} Nesterov sequence,
and the accelerated monotone PGS (AM-PGS), \textit{i.e.,} Beck's sequence, are implemented as pseudocode in Algorithm \ref{algorithm: a unified implementation of PGS, A-PGS, AM-PGS methods}.
An illustration of the convergence is given in Fig.~\ref{fig. ther convergence validation from deltaX and tangentV}.

\begin{algorithm}[t]
	
	\DontPrintSemicolon
	
	\SetKwData{Method}{method}
	\SetKwData{PGS}{PGS}
	\SetKwData{APGS}{A-PGS}
	\SetKwData{AMPGS}{AM-PGS}
			
	\SetKwInOut{Input}{input}
	\SetKwInOut{Output}{output}
	\SetKwInOut{Methods}{method}
	
	\SetKwFunction{lineSearchFun}{\textbf{lineSearch}}
	\SetKwFunction{searchMaxProxyStepsizeFun}{\textbf{searchMaxProxyStepsize}}
	
	\Input{\Method = \PGS or \APGS or \AMPGS \;}
	
	\Input{$\vecf{x}_0$\;}
	
	\lIf{known Lipschitz constant $L$}{
	$t'_{\max} \gets 1/L$
	}
	\lElse{
		$t'_{\max} \gets$ \searchMaxProxyStepsizeFun $(\vecf{x}_0)$
	}

$\vecf{y}_0 \gets \vecf{x}_0$, $\alpha_0 \gets 1$, $k \gets 0$\;
	
\While{ $k < maxIterations$ }{
	
	{
			
		$(\vecf{v}_{k},\,t,\,t') \gets $
		\lineSearchFun$(\vecf{y}_k,\,t'_{\max})$\;

		$\vecf{y}^{\star} \gets \mathcal{R}_{\vecf{y}_k}\left( \vecf{v}_k \right)$\;

		\lIf{AdaptiveMaxProxyStepsize}{ 
			$t'_{\max} \gets t'$
		}	
	}
					
\uIf{\Method = \PGS}{
	$\vecf{x}_{k+1} \gets \vecf{y}^{\star}$\;
	$\vecf{y}_{k+1} \gets \vecf{y}^{\star}$
}
\ElseIf{\Method = \APGS or \AMPGS}{
	$\alpha_{k+1} \gets \frac{1+\sqrt{1+4 \alpha_k^2}}{2}$\;
	$\vecf{x}_{k+1} \gets \vecf{y}^{\star} $\;
	$\vecf{y}_{k+1} \gets \mathcal{R}_{\vecf{y}^{\star}} \left( \frac{1-\alpha_k}{\alpha_{k+1}} \mathcal{R}_{\vecf{y}^{\star}}^{-1} \left( \vecf{x}_{k} \right) \right)$\;
	
	\If{\Method = \AMPGS }{
		\If{$g(\vecf{y}^{\star}) + h(\vecf{y}^{\star}) > g(\matf{x}_{k}) + h(\matf{x}_{k}) $}{
		$\vecf{x}_{k+1} \gets \vecf{x}_{k} $\;
		$\vecf{y}_{k+1} \gets \mathcal{R}_{\vecf{x}_{k}} \left( \frac{\alpha_k}{\alpha_{k+1}} \mathcal{R}_{\vecf{x}_{k}}^{-1} \left( \vecf{y}^{\star} \right) \right)$\;
		}
	}
}

\uIf{ $\Vert \vecf{v}_{k} \Vert_2 < 1e-5$ and $\Vert \vecf{v}_{k} / t \Vert_2 < 1e-3$ }{	
	\Return{$\vecf{x}_{k+1}$ }
}
\lElse{
	$k \gets k + 1$
}

}
		
\caption{A unified implementation of the PGS, A-PGS and AM-PGS methods.
\label{algorithm: a unified implementation of PGS, A-PGS, AM-PGS methods}}
\end{algorithm}


\section{Applications}
\label{section. applications. theoretical description}

\subsection{Rayleigh Quotient Optimization}
\label{section. Rayleigh Quotient Optimization}

We consider $g (\vecf{x}) = \vecf{x}^{\tran} \matf{A}  \vecf{x}$ with $\matf{A}$ being a symmetric matrix.
The Euclidean gradient at $\vecf{x}\in\mathbb{R}^n$ is $\nabla g \vert_{\vecf{x}} = 2 \matf{A} \vecf{x}$.
The Riemannian gradient at $x \in \mathcal{S}$ is:
\begin{multline}
\label{eq: Riemannian gradient of ralyleigh quotient optimization}
\mathrm{grad}\,g (\vecf{x})
=
\nabla g \vert_{\vecf{x}}
-
\left\langle \vecf{x}, \nabla g \vert_{\vecf{x}} \right\rangle  \vecf{x}
\\ =
2 \matf{A} \vecf{x}
-
2 \left\langle \vecf{x}, \matf{A} \vecf{x} \right\rangle  \vecf{x}
=
2 \matf{A} \vecf{x} - 2 (\vecf{x}^{\tran}  \matf{A} \vecf{x}) \vecf{x}.
\end{multline}
Using retraction (\ref{eq: retraction of unit sphere}), the Lipschitz-type constant of $g \left( \mathcal{R}_{\vecf{x}} \left( \vecf{v} \right) \right)$ is $L = 2 \sigma_{\max} (\matf{A})$ (see Appendix \ref{appendix: Lipschitz-type Constant of rayleigh quotient})
where $\sigma_{\max} (\matf{A})$ denotes the largest singular value of $\matf{A}$.

Minimizing $g (\vecf{x}) = \vecf{x}^{\tran} \matf{A}  \vecf{x}$ on the sphere manifold is
called Rayleigh quotient optimization:
\begin{equation}
\label{eq: ralyleight quotient optimization}
\vecf{x}_0 = 
\argmin_{\matf{x}\, \in\, \mathbb{R}^n}\ 
\vecf{x}^{\tran} \matf{A}  \vecf{x}
\qquad\mathrm{s.t.}\quad
\left\Vert
\vecf{x}
\right\Vert_2
=1.
\end{equation}
The solution to this problem is given in closed-form, which is the bottom eigenvector of $\matf{A}$ (\textit{i.e.,} the eigenvector associated with the smallest eigenvalue), which often gives as a good initialization for its regularized versions.

\subsection{Fundamental Matrix Estimation}

\subsubsection{Problem Statement}

The fundamental matrix is a key algebraic model of the two-view geometry \cite{faugeras1995geometry, Hartley2004, zheng2013practical}.
We denote $\vecf{p}_i \leftrightarrow {\vecf{p}_i'}$ $\left(i = \left[1:m\right]\right)$ the homogeneous coordinates of corresponding points in two images. In the noise-free case, the epipolar constraint holds as:
	\begin{equation}
	\label{eq: the fundamental matrix estimation of two view geometry}
	{\vecf{p}_i'}^{\tran} \matf{F} \vecf{p}_i
	=
	0,
	\quad
	i = \left[1:m\right].
	\end{equation}
The matrix $\matf{F} \in \mathbb{R}^{3 \times 3}$ is called the \textit{fundamental matrix}, defined up to scale, thus we seek for $\matf{F}$ on the unit sphere such that $\left\Vert \matf{F} \right\Vert_F = 1$.
Importantly, $\mathrm{rank}(\matf{F}) = 2$ is required.
A related concept is the essential matrix for normalized calibration, for which we refer to a recent work~\cite{zhao2020efficient}.

To estimate $\matf{F}$, we can formulate a cost function based on the algebraic error as:
	\begin{align*}
	\label{eq: The fundamental matrix problem in matrix form}
	\varphi(\matf{F}) & =   
	\frac{1}{m}
	\sum_{i=1}^{m}
	\left(
	{\vecf{p}_i'}^{\tran} \vecf{F} \vecf{p}_i \right)^2
	 =
	\frac{1}{m}
	\sum_{i=1}^{m}
	\left(
	\left(
	\vecf{p}_i^{\tran}
	\otimes
	{\vecf{p}_i'}^{\tran}
	\right)
	\mathrm{vec}\left( \matf{F} \right)
	\right)^2
	\\ & =
	\frac{1}{m}
	\left\Vert
	\matf{H} \mathrm{vec}\left(\matf{{F}}\right)
	\right\Vert_{2}^2
	,
	\end{align*}
where the $i$-th row of $\matf{H}$ is
$
\vecf{p}_i^{\tran}
\otimes
{\vecf{p}_i'}^{\tran}
$.
We denote $\vecf{x} = 
\mathrm{vec}\left(\matf{{F}}\right)$
and $\matf{F} = \mathrm{mat}\left(\vecf{{x}}\right)$,
where $\mathrm{vec}(\cdot)$ is the standard column-wise matrix vectorization and $\mathrm{mat}(\cdot)$ is its inverse operation. Defining $\matf{A} = \frac{1}{m} \matf{H}^{\tran}\matf{H}$, we see the fundamental matrix problem is an instance of problem (\ref{eq: ralyleight quotient optimization}).

The solution $\matf{F}_0 = \mathrm{mat}(\vecf{x}_{0})$ solved from problem (\ref{eq: ralyleight quotient optimization}) is usually of rank $3$. A remedy is to subsequently round the solution $\matf{F}_0$ using the rank-$2$ approximation via the Singular Value Decomposition (SVD). This two-stage solution is the eight-point algorithm. Instead, we give a low-rank solution using the nuclear norm regularization.

\subsubsection{Nuclear Norm Regularization}
\label{subsubsection: nuclear norm regularization for low rank, applied in fundamental matrix estimation}

The nuclear norm of a matrix $\left\Vert  \matf{X} \right\Vert_{\ast}$,
defined as the summation of its singular values, is the tightest convex envelop of the rank function
within the unit ball
$
\left\{
\matf{X} \in \mathbb{R}^{m \times n} \,:\, \left\Vert  \matf{X} \right\Vert_2 \le 1
\right\}
$~\cite{fazel2002matrix}.
Nuclear norm regularization has been widely used as a technique to promote low-rank
in Euclidean optimization problems \cite{fazel2013hankel},
	while a direct deployment to the manifold setting seems to be obscure with the results in \cite{chen2020proximal, xiao2018regularized},
	mainly due to the challenge incurred in evaluating the generalized Jacobian matrix of a non-smooth function.
	In contrast, our technique can handle nuclear norm regularization with no effort.
	We apply nuclear norm regularization to problem (\ref{eq: ralyleight quotient optimization}) as:
\begin{equation}
\label{eq: formulation for fundamental matrix problem with nuclear norm regularization}
\argmin_{\matf{x}\, \in\, \mathbb{R}^n}\ 
\vecf{x}^{\tran} \matf{A}  \vecf{x}
+
\lambda
\Vert \mathrm{mat}\left(
\vecf{x}
\right) \Vert_{\ast}
\quad \mathrm{s.t.}\quad \Vert \vecf{x} \Vert_2 = 1
,
\end{equation}
where $\lambda \ge 0$ is a given constant controlling the strength of regularization.

	In problem (\ref{eq: formulation for fundamental matrix problem with nuclear norm regularization}),
	the regularization term is
	$
	h \left( \vecf{x} \right) = \lambda \left\Vert \mathrm{mat} \left( \vecf{x} \right) \right\Vert_{\ast}
	$,
	which is convex but non-smooth.
	Besides, $h \left( \cdot \right)$
	is absolutely homogeneous,
	\textit{i.e.,}
	$
	h \left( \alpha \vecf{x} \right) = \left\vert \alpha \right\vert h  \left( \vecf{x} \right)
	$,
	because:
\begin{equation}
\lambda
\left\Vert \mathrm{mat} \left( \alpha \vecf{x} \right) \right\Vert_{\ast}
=
\lambda
\left\Vert \alpha \mathrm{mat} \left( \vecf{x} \right) \right\Vert_{\ast}
=
\lambda
\left\vert \alpha \right\vert \left\Vert \mathrm{mat} \left( \vecf{x} \right) \right\Vert_{\ast}
.
\end{equation}
The first equality holds because $\mathrm{mat} \left( \vecf{x} \right)$ is a linear operator,
and
the second because norms are absolutely homogeneous.

	Let $\matf{X} = \mathrm{mat}\left(
	\vecf{x}
	\right)$
	and denote its SVD as
	$\matf{X} = \matf{U} \boldsymbol{\Sigma} \matf{V}^{\tran}$.
For the nuclear norm function $\lambda \left\Vert \matf{X} \right\Vert_{\ast}$,
the proximal is given in closed-form \cite{parikh2014proximal},
where
$\mathrm{prox}_{t \lambda \left\Vert \cdot \right\Vert_{\ast}}
(\matf{X})
=
\matf{U}
(\boldsymbol{\Sigma}
-
t \lambda \matf{I})_{+}
\matf{V}^{\tran}
$.
Here $\matf{A}_{+} \defeq \max\{\matf{A},\,\matf{0}\}$ proceeds element-wise.
Thus we obtain the proximal of $h \left(\cdot \right)$ as:
\begin{equation}
\mathrm{prox}_{t h}
\left(\vecf{x} \right)
= 
\mathrm{vec}
\left(
\matf{U}
\left( \boldsymbol{\Sigma}
-
t \lambda \matf{I} \right)_{+}
\matf{V}^{\tran}
\right)
.
\end{equation}

\subsection{Correspondence Association}

\subsubsection{Problem Statement}

The correspondence association problem using pairwise constraints can be formulated as a Rayleigh quotient optimization as well \cite{leordeanu2005spectral}.
We denote the association hypothesis that a point $i$ in the point-cloud $\mathcal{Q}$ is matched with a point $i'$ in the point-cloud $\mathcal{Q}'$ as $\mathfrak{h}_{ii'}$.
The correspondence problem is to estimate the likelihood $p(\mathfrak{h}_{ii'})$ $(i \in \mathcal{Q},\, i' \in \mathcal{Q}' )$ of all possible association hypotheses,
collected as components of the state vector $\vecf{x}$.

To that end, we can design an adjacency matrix $\matf{M}$ from the pairwise consistency of hypotheses \cite{leordeanu2005spectral}.
A common practice is based on the change of distance:
\begin{equation}
\label{eq. the adjacency matrix of grpah association}
\matf{M}(\mathfrak{h}_{ii'}, \mathfrak{h}_{jj'}) =
\begin{cases}
4.5 - \frac{(d_{ij} - d_{i'j'})^2}{2 \delta_d^2}  & \mathrm{if}\ \left\vert  d_{ij} - d_{i'j'} \right\vert < 3 \delta_d \\
0  & \mathrm{otherwise},
\end{cases}
\end{equation}
where $d_{ij}$ is the Euclidean distance between the points $i$ and $j$ in $\mathcal{Q}$, $d_{i'j'}$ the distance between the points $i'$ and $j'$ in $\mathcal{Q}'$, and $\delta_d$ a tuning parameter.

The resulting problem is formalized as maximizing the overall consistency
$\vecf{x}^{\tran} \matf{M}  \vecf{x}$ on the unit sphere, as an instance of problem (\ref{eq: ralyleight quotient optimization})
by letting $\matf{A} = - \matf{M}$.
The estimate of $\vecf{x}$ is further used to decide the final correspondences, based on various assumptions, \textit{e.g.,} one point in $\mathcal{Q}$ can only be matched with one point in $\mathcal{Q}'$ \cite{leordeanu2005spectral}.

The match hypotheses solved this way are dense, while many of them present with contradictions or low probabilities.
It is thus favorable to have a sparse $\vecf{x}$, where some unlikely hypotheses and contradictions are pruned away.
We give such a sparse solution by $\ell_1$ norm regularization.

\subsubsection{$\ell_1$ Norm Regularization}

It has been known that $\ell_1$ norm regularization can favor sparsity in Euclidean \cite{beck2009fast_FISTA} and manifold optimization \cite{chen2020proximal, xiao2018regularized}.
We apply $\ell_1$ norm regularization to problem (\ref{eq: ralyleight quotient optimization}) as:
\begin{equation}
\label{eq: formulation for graph correspondence association problem with L1 regularization}
\argmin_{\matf{x}\, \in\, \mathbb{R}^n}\ 
\vecf{x}^{\tran} \matf{A}  \vecf{x}
+
\lambda
\Vert \vecf{x} \Vert_{\ell_1}
\quad \mathrm{s.t.}\quad \Vert \vecf{x} \Vert_2 = 1
,
\end{equation}
where $\lambda \ge 0$ is a given constant controlling the strength of regularization.
In this case, $h (\vecf{x}) = \lambda
\Vert \vecf{x} \Vert_{\ell_1}$,
which is convex and absolutely homogeneous.
The proximal of $h (\vecf{x})$, often termed soft shrinkage operator, is given element-wise as:
\begin{equation}
\mathrm{prox}_{t h}
\left(\vecf{x} \right)_{i}
= 
\mathrm{sgn}(\vecf{x}_i)
(
\left\vert \vecf{x}_i \right\vert
- t \lambda
)_{+}
.
\end{equation}

\subsection{Camera Self-calibration}
\label{subsection: camera self-calibration}

\subsubsection{Problem Statement}

In projective reconstruction, we obtain a set of projective cameras $\tilde{\matf{P}}_i \in \mathbb{R}^{3 \times 4}$:
\begin{equation}
\tilde{\matf{P}}_i \propto \matf{P}_i \matf{H}^{-1}
\quad
(i \in \left[1:n\right]),
\end{equation}
which differ from the Euclidean cameras $\matf{P}_i \in \mathbb{R}^{3 \times 4}$ by a common projective transformation $\matf{H} \in \mathbb{R}^{4 \times 4}$.
The camera self-calibration problem is to infer $\matf{H}$ from $\tilde{\matf{P}}_i$ \cite{Hartley2004}.

The key algebraic model to this task is the \textit{Dual Absolute Quadric (DAQ)},
a rank-$3$ symmetric matrix in $\mathbb{R}^{4 \times 4}$ defined up to scale \cite{triggs1997autocalibration}.
In specific,
the DAQ in the Euclidean space, termed the canonical DAQ, takes the form $\Omega_{\infty}^{\ast} = \mathrm{diag}(1,1,1,0)$.
We denote the DAQ in the projective space (where $\tilde{\matf{P}}_i$ is defined) by $\matf{Q}_{\infty}^{\ast} = \matf{H} \Omega_{\infty}^{\ast} \matf{H}^{\tran}$.
The image of the DAQ, denoted by $\omega_i^{\ast}$, is invariant under $\matf{H}$:
\begin{equation}
\label{eq: key property used for self-calibration}
\omega_i^{\ast}
\propto
\tilde{\matf{P}}_i \matf{Q}_{\infty}^{\ast} \tilde{\matf{P}}_i^{\tran}
\propto
\matf{P}_i
\Omega_{\infty}^{\ast}
\matf{P}_i^{\tran}
\propto
\matf{K}_i \matf{K}_i^{\tran}
,
\end{equation}
where $\matf{K}_i \in \mathbb{R}^{3 \times 3}$ is the intrinsic matrix of $\matf{P}_i$.
At its core, the self-calibration problem is to estimate $\matf{Q}_{\infty}^{\ast}$ from equation (\ref{eq: key property used for self-calibration}) using various constraints on $\matf{K}_i$.

Here we consider a linear approach developed for cameras with varying focal lengths \cite{pollefeys1999self}.
In this case, we have
$\matf{K}_i \matf{K}_i^{\tran} \propto \mathrm{diag}(f_i^2, f_i^2, 1)$. Based on equation (\ref{eq: key property used for self-calibration}), we have:
\begin{subnumcases}{ \label{eq: set of equations fo r linear self-calibration}}
\tilde{\matf{a}}_i^{\tran} \matf{Q}_{\infty}^{\ast} \tilde{\matf{a}}_i
=
\tilde{\matf{b}}_i^{\tran} \matf{Q}_{\infty}^{\ast} \tilde{\matf{b}}_i  \\
\tilde{\matf{a}}_i^{\tran} \matf{Q}_{\infty}^{\ast} \tilde{\matf{b}}_i = 0,\ 
\tilde{\matf{a}}_i^{\tran} \matf{Q}_{\infty}^{\ast} \tilde{\matf{c}}_i = 0 ,\ 
\tilde{\matf{b}}_i^{\tran} \matf{Q}_{\infty}^{\ast} \tilde{\matf{c}}_i = 0,
\end{subnumcases}
where
$
\tilde{\matf{P}}_i^{\tran} = 
\begin{bmatrix}
\tilde{\matf{a}}_i &
\tilde{\matf{b}}_i &
\tilde{\matf{c}}_i
\end{bmatrix}
$.
Equation (\ref{eq: set of equations fo r linear self-calibration}) is linear in $\matf{Q}_{\infty}^{\ast}$ thus can be rewritten as:
\begin{equation}
\matf{M}_i \,
\mathrm{vec_t}\left( \matf{Q}_{\infty}^{\ast} \right)
=
\matf{M}_i \, \vecf{x} = \vecf{0}
,
\end{equation}
where we have defined
$\mathrm{vec_t}\left( \matf{Q}_{\infty}^{\ast} \right)
= \vecf{x} \in \mathbb{R}^{10}
$
comprising of the upper triangular elements of $\matf{Q}_{\infty}^{\ast}$,
and
$\mathrm{mat_t}\left( \vecf{x} \right)
=
\matf{Q}_{\infty}^{\ast}$ its inverse operation.
Since $\matf{Q}_{\infty}^{\ast}$ is defined up to scale so is $\vecf{x}$ and we minimize the cost
$
\varphi(\vecf{x})
=
\frac{1}{n}
\sum_{i=1}^{n}
\left\Vert \matf{M}_i \, \vecf{x} \right\Vert_2^2
$
on the unit sphere.
Upon defining
$\matf{A}
=
\frac{1}{n}
\sum_{i=1}^{n}
\matf{M}_i^{\tran}
\matf{M}_i
$,
we see it is another instance of problem (\ref{eq: ralyleight quotient optimization}).

Just as the case of the fundamental matrix estimation, the DAQ estimated from solving problem (\ref{eq: ralyleight quotient optimization}) is typically of rank $4$ instead of rank $3$, thus an SVD based rounding process is used subsequently.

Once obtaining a rank-$3$ estimate of $\matf{Q}_{\infty}^{\ast}$,
we can recover $\matf{H}$ up to a similarity transformation.
This is usually done by the eigen decomposition of
$\matf{Q}_{\infty}^{\ast}$.
Let
$
\matf{Q}_{\infty}^{\ast} = \matf{U} \matf{\Lambda} \matf{U}^{\tran}
$,
with $\matf{\Lambda} = \mathrm{diag}(\lambda_1, \lambda_2, \lambda_3, 0)$.
We then set $\matf{H} \propto \matf{U}
\, \mathrm{diag}(\sqrt{\lambda_1}, \sqrt{\lambda_2}, \sqrt{\lambda_3}, 1)
$.
The estimate of $\matf{H}$ determines camera $\matf{P}_i$, and $\matf{K}_i$ afterwards by decomposing $\matf{P}_i$.

\subsubsection{Nuclear Norm Regularization}

One issue regarding self-calibration is the critical motion sequences (CMS) \cite{gurdjos2009dual}.
The CMSs are camera configurations where self-calibration is ambiguous due to the lack of sufficient constraints.
Among which, we consider the artificial CMS that can be resolved by enforcing the rank deficiency of the DAQ during the estimation rather than a posteriori.
For the linear self-calibration described in equation (\ref{eq: set of equations fo r linear self-calibration}),
one such CMS is that - all the cameras' principal axes intersect at a fixed point,
\textit{i.e.,}~all the cameras look towards a common point.
In this case, there only exist two rank deficient solutions \cite{gurdjos2009dual}: one rank-$3$ solution (desired) and one rank-$1$ solution (undesired).

In analogy to the fundamental matrix estimation, we use nuclear norm regularization to promote low rank.
This case is similar to what we have discussed in Section \ref{subsubsection: nuclear norm regularization for low rank, applied in fundamental matrix estimation}.
We omit the details as they can be readily derived
by replacing $\mathrm{mat}(\cdot)$ with $\mathrm{mat}_t(\cdot)$ and $\mathrm{vec}(\cdot)$ with $\mathrm{vec}_t(\cdot)$
in problem (\ref{eq: formulation for fundamental matrix problem with nuclear norm regularization}).

\subsubsection{Nuclear-Spectral Norm Regularization}

The aforementioned nuclear norm regularization resolves the CMS only partly due to the existence of the rank-$1$ solution.
We propose to avoid the rank-$1$ solution by additionally including a spectral norm to penalize the largest singular value.
The nuclear-spectral norm regularizer is:
\begin{equation}
\label{eq. nuclear-spectral norm regularizer}
h(\vecf{x}) = \lambda_1
\Vert \mathrm{mat_t}\left(
\vecf{x}
\right) \Vert_{\ast}
+
\lambda_2
\Vert \mathrm{mat_t}\left(
\vecf{x}
\right) \Vert_2,
\end{equation}
where $\lambda_1$ and $\lambda_2$ control the regularization strength of each part. 
We apply this regularizer to problem (\ref{eq: ralyleight quotient optimization}) as:
\begin{equation}
\label{eq: formulation for auto-calibration with nuclear plues spectral norm regularization}
\begin{aligned}
\argmin_{\matf{x}\, \in\, \mathbb{R}^n}\ 
& \vecf{x}^{\tran} \matf{A}  \vecf{x}
+
\lambda_1
\Vert \mathrm{mat_t}\left(
\vecf{x}
\right) \Vert_{\ast}
+
\lambda_2
\Vert \mathrm{mat_t}\left(
\vecf{x}
\right) \Vert_2
\\
\mathrm{s.t.}\quad & \Vert \vecf{x} \Vert_2 = 1 .
\end{aligned}
\end{equation}

Denote the SVD of $\mathrm{mat_t}\left(
\vecf{x} \right)$ as
$\mathrm{mat_t}\left(\vecf{x} \right) = \matf{U} \boldsymbol{\Sigma} \matf{V}^{\tran}$.
Since both the nuclear and spectral norm are orthogonal invariant, we can derive the proximal of $h(\cdot)$ in a similar manner as the derivation used for the nuclear norm.
The proximal of $h(\cdot)$ in equation (\ref{eq. nuclear-spectral norm regularizer}) is given as follow:
\begin{equation}
\mathrm{prox}_{t h}
\left(\vecf{x} \right)
= 
\mathrm{vec_t}
\left(
\matf{U}
\left( \boldsymbol{\Sigma}
-
t \lambda_1 \matf{I} - t \lambda_2 \matf{E}_1  \right)_{+}
\matf{V}^{\tran}
\right),
\end{equation}
where $\matf{E}_1=\mathrm{diag}(1,0,\cdots,0)$ is a diagonal matrix where the top-left element is $1$ and the rests are all-zeros.

\section{Experimental Results}
\label{section. experimental results}

\subsection{The Proposed PGS Algorithm}
\label{section: experimental results: the proposed PGS algorithm}

We first provide an evaluation of the PGS methods, \textit{i.e.,}~PGS, A-PGS and AM-PGS methods in Algorithm \ref{algorithm: a unified implementation of PGS, A-PGS, AM-PGS methods}.

\subsubsection{Experiment Setup}
\label{subsection. experiment setup of PGS methods}

\noindent\textbf{Numerical instances.}
We experiment with the Rayleigh quotient optimization,
and draw numerical examples from different applications which essentially form different $\matf{A}$ matrices in problem (\ref{eq: ralyleight quotient optimization}):
\begin{itemize}
\item nuclear norm reg. --- fundamental matrix estimation with nuclear norm regularization,
\item $\ell_1$ norm reg. --- correspondence association with $\ell_1$ norm regularization,
\item nuclear-spectral norm reg. --- self-calibration with nuclear-spectral norm regularization.
\end{itemize}
In this section, we distinguish these instances by the type of regularization used, \textit{i.e.,}~nuclear, $\ell_1$ and nuclear-spectral.

\vspace{5pt}\noindent\textbf{Initialization.}
The numerical examples used in this section are special forms ot the Rayleigh quotient optimization problem (\ref{eq: ralyleight quotient optimization}),
and its optimal solution is known to be the bottom eigenvector of matrix $\matf{A}$ which we denote by $\boldsymbol{\xi}_{\matf{A}}$.
To examine the convergence behavior with respect to different initializations, we create a range of initial values $\vecf{x}_0$ by adding independent zero-mean Gaussian noise element-wisely to $\boldsymbol{\xi}_{\matf{A}}$. For the $k$-th vector element, we let:
\begin{equation}
\label{eq. initialization by sigma}
\boldsymbol{\xi}[k] \gets
\boldsymbol{\xi}_{\matf{A}} [k] + \epsilon
,\quad
\epsilon \sim \mathcal{N}(0, \delta_{\mathrm{init}})
,
\end{equation}
and normalize $\boldsymbol{\xi}$ to get the initial value $\vecf{x}_0 = \boldsymbol{\xi} / \Vert \boldsymbol{\xi} \Vert_2$.
Intuitively, $\delta_{\mathrm{init}}$ controls the deviation from the eigenvector initialization $\boldsymbol{\xi}_{\matf{A}}$, and if $\delta_{\mathrm{init}}$ is large, the above process simulates random initialization.

\vspace{5pt}\noindent\textbf{Proxy step-size strategies.}
We examine the following proxy step-size strategies.
\begin{itemize}
\item LipschitzFixed --- $t'_{\max} = 1/L$ initially, and $t'_{\max}$ is kept fixed in the following iterations;
\item LipschitzAdaptive --- $t'_{\max} = 1/L$ initially, and at each iteration $t'_{\max}$ is updated to the previous working proxy step-size $t'$;
\item SearchedFixed --- $t'_{\max}$ is obtained from Algorithm \ref{algorithm: adaptive search for maximum proxy stepsize} initially, and $t'_{\max}$ is kept fixed in the following iterations;
\item SearchedAdaptive --- $t'_{\max}$ is obtained from Algorithm \ref{algorithm: adaptive search for maximum proxy stepsize}, and at each iteration $t'_{\max}$ is updated to the previous working proxy step-size $t'$.
\end{itemize}

\subsubsection{Results}

We report the results with a $20$ run Monte-Carlo simulation for each PGS method and each proxy step-size strategy with respect to different initializations.

\begin{figure}[t]
\centering
\includegraphics[width=0.48\textwidth]{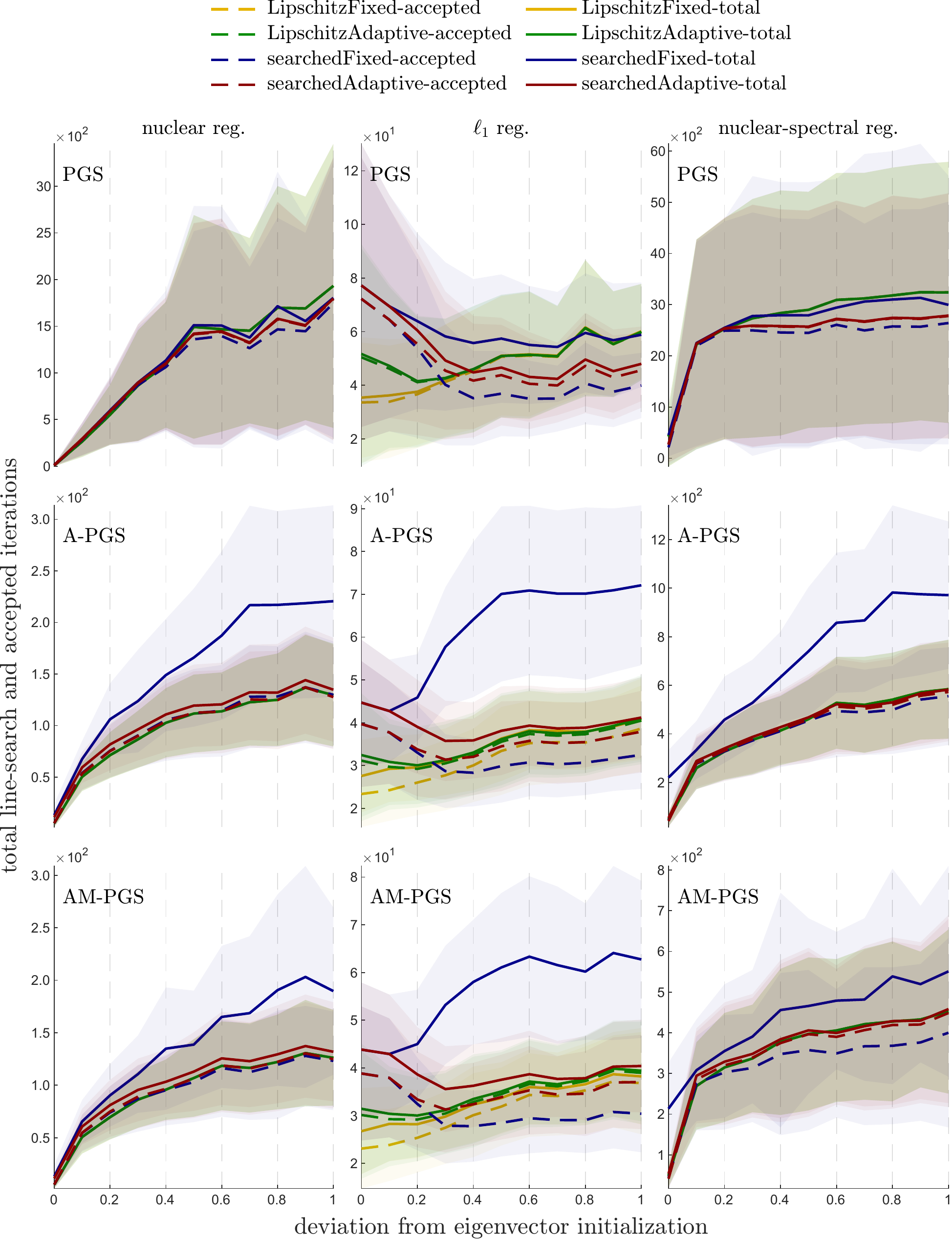}
\caption{The convergence of the PGS, A-PGS and AM-PGS methods by using different proxy step-size strategies. Results are reported as the total line-search and the accepted iterations with respect to different initialization.}
\label{fig. Proxy step-size with adaptive search and Lipschitz constant. Iterations}	
\end{figure}

\vspace{5pt}\noindent\textbf{Convergence by used iterations.}
We examine the convergence by the used iterations in Fig.~\ref{fig. Proxy step-size with adaptive search and Lipschitz constant. Iterations}.
In particular, we distinguish the total line-search iterations (the overall iterations processed in Algorithm \ref{algorithm: our proposed algorithm for RPG on the sphere}) and the accepted iterations (the iterations used in Algorithm \ref{algorithm: a unified implementation of PGS, A-PGS, AM-PGS methods}).
The final verdicts is as follows:
a) In general, we observe no significant differences for the accepted iterations across different proxy step-size strategies.
b) However, the proxy step-size strategy SearchedFixed is not recommended, as it often leads to line-search failures especially when the initialization is bad, as reflected by the number of total line-search iterations.
Therefore, if $t'_{\max}$ is obtained from line-search in Algorithm \ref{algorithm: adaptive search for maximum proxy stepsize}, we recommend at each iteration updating $t'_{\max}$ to the previous working proxy step-size $t'$.
c) If $t'_{\max}$ is initialized as $1/L$ from the Lipschitz constant $L$, both strategies LipschitzFixed and LipschitzAdaptive give similar results.
d) We observe accelerated methods A-PGS and AM-PGS converge much faster than the unaccelerated PGS method, thus are generally recommended.

\begin{figure}[t]
\centering
\includegraphics[width=0.47\textwidth]{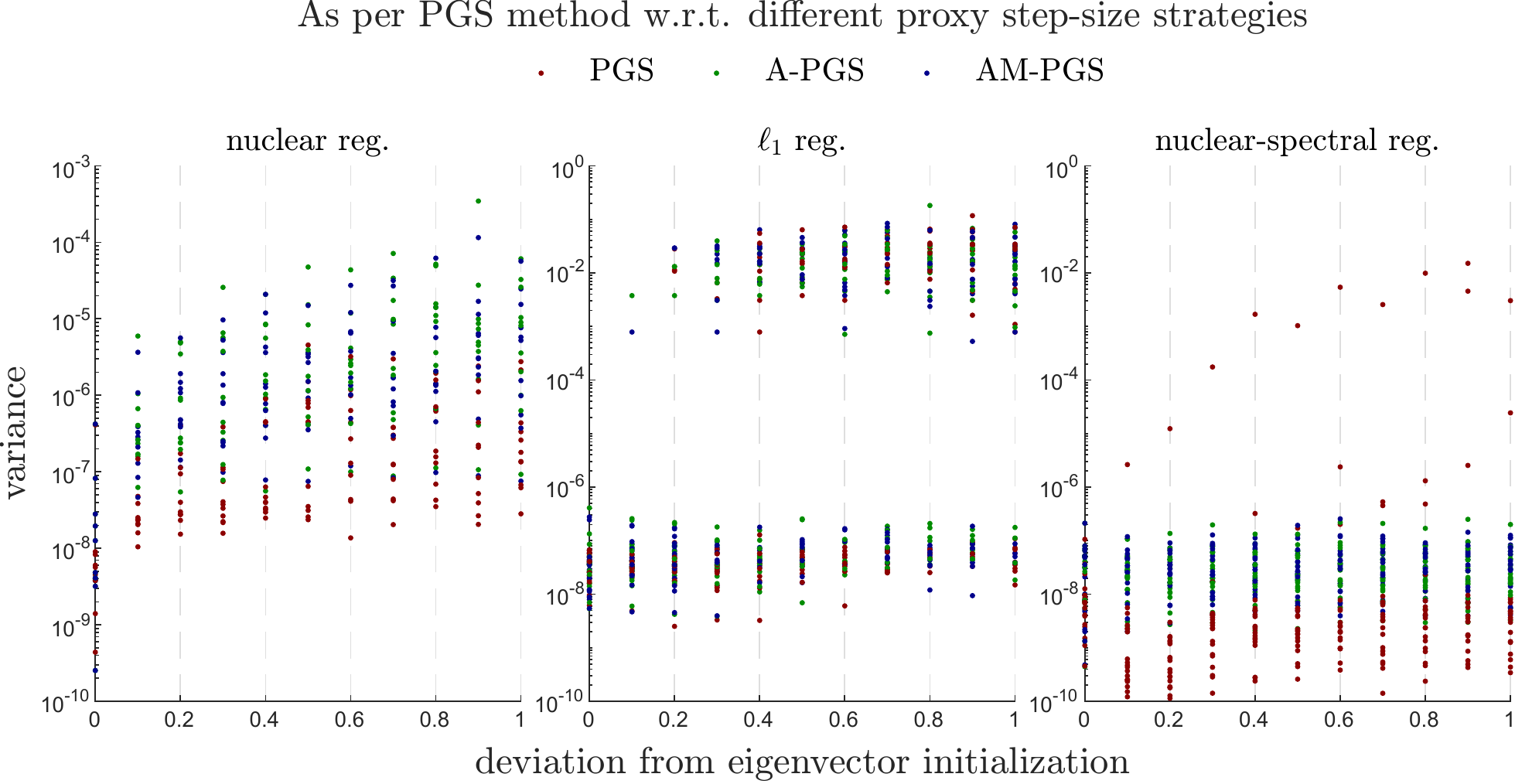}
\\[10pt]
\includegraphics[width=0.47\textwidth]{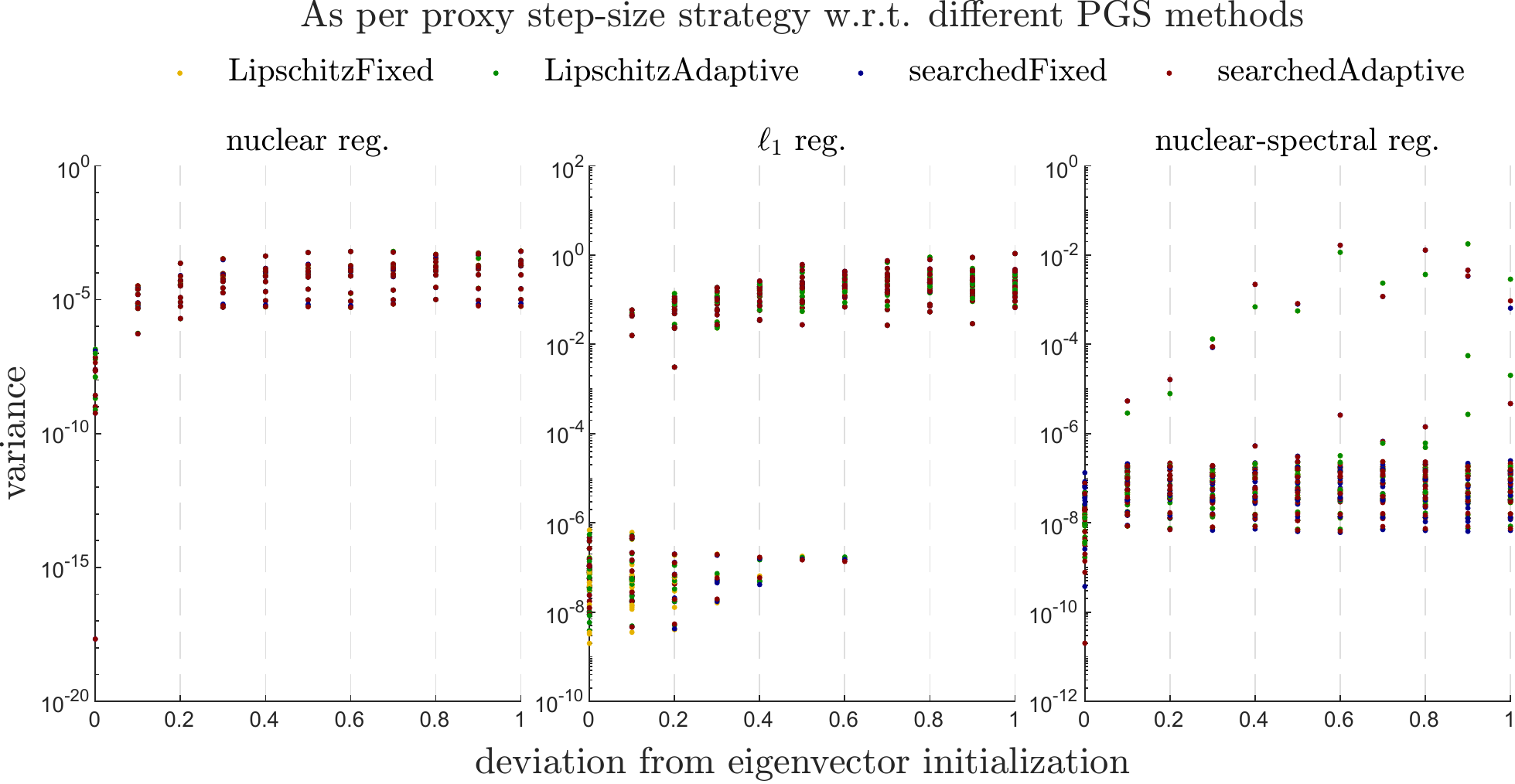}
\caption{The variation of the optimal costs for each case.}
\label{fig. Proxy step-size with adaptive search and Lipschitz constant. Variation of the final costs}	
\end{figure}

\vspace{5pt}\noindent\textbf{Variation of the optimal costs.}
We compare the optimal costs of each Monte-Carlo run in Fig.~\ref{fig. Proxy step-size with adaptive search and Lipschitz constant. Variation of the final costs}.
In specific, in the top figure, for each PGS method (\textit{i.e.,} PGS/A-PGS/AM-PGS), we compare the costs obtained from different proxy step-size strategies.
Likewise, in the bottom figure, for each proxy step-size strategy, we compare the costs obtained from different PGS methods.
The difference is evaluated as the variation of the optimal costs.
From Fig.~\ref{fig. Proxy step-size with adaptive search and Lipschitz constant. Variation of the final costs}, we see that the optimal costs are mostly the same if $\delta_{\mathrm{init}}$ in equation (\ref{eq. initialization by sigma}) is small, or otherwise stated if the initialization $\vecf{x}_0$ is close to the eigenvector initialization $\boldsymbol{\xi}_{\matf{A}}$ where we simulate good initializations.
As $\delta_{\mathrm{init}}$ grows where we simulate bad initializations, the differences grow as different PGS methods and proxy step-size strategies can lead to the convergence to different local minima.
This phenomenon is extremely clear for the $\ell_1$ norm regularized instances, where $\vecf{x}$ is valued mostly below $0.5$ (see the example in Fig.~\ref{fig. graph assoication hypothesis visualization and mc simulation}) and $\delta_{\mathrm{init}} > 0.5$ sets $\vecf{x}_0$ almost to random.

\begin{figure}[t]
\centering
\includegraphics[width=0.47\textwidth]{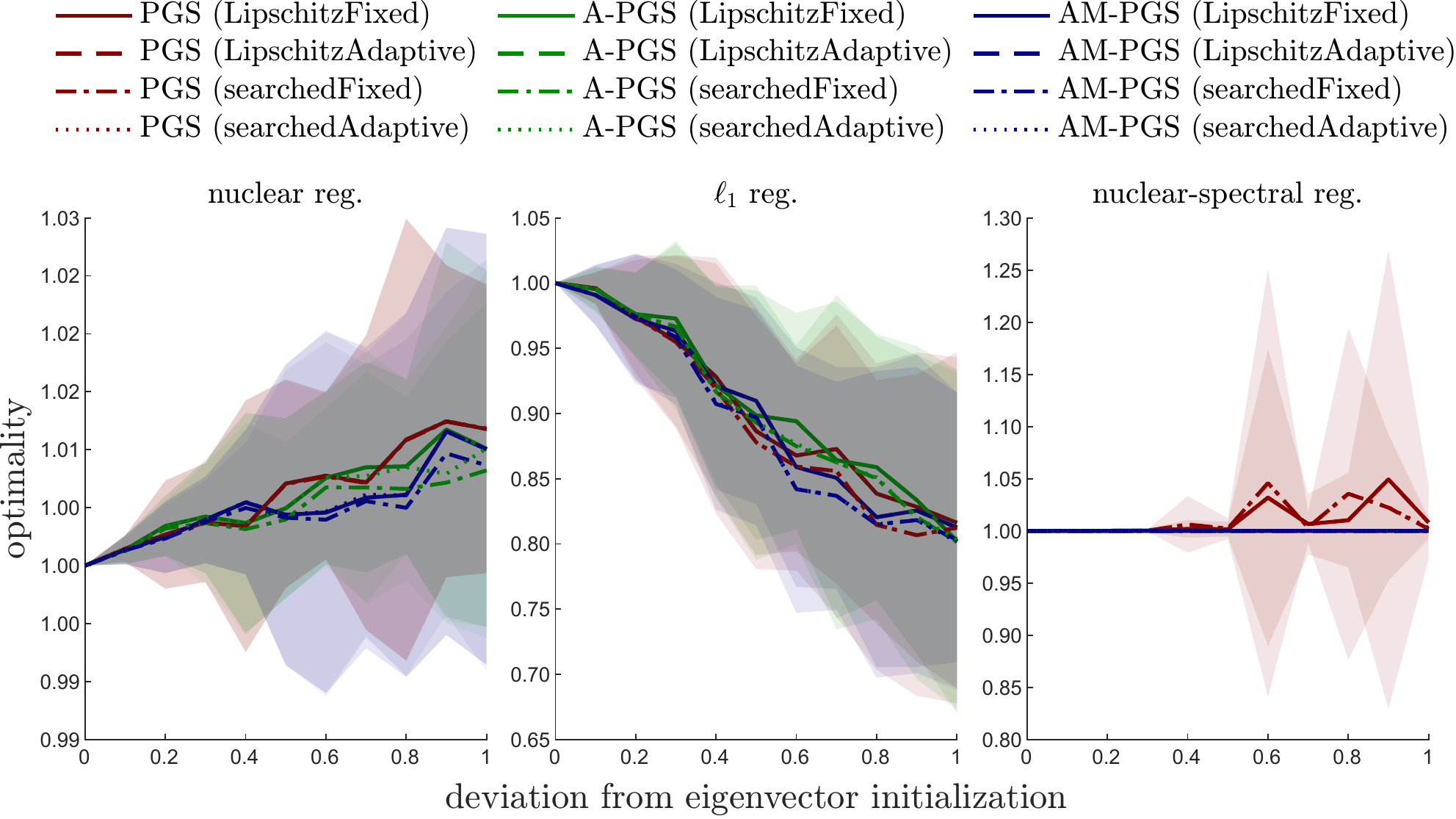}
\caption{The impact of different initializations.}
\label{fig. Proxy step-size with adaptive search and Lipschitz constant. Optimality of the final costs}	
\end{figure}

\vspace{5pt}\noindent\textbf{Optimality of the optimal costs.}
We evaluate the optimal cost obtained from the initialization $\vecf{x}_0$ by comparing with the optimal cost obtained from the eigenvector initialization $\boldsymbol{\xi}_{\matf{A}}$.
For each Monte-Carlo run, we define the optimality:
\begin{equation*}
\mathrm{optimality} = \frac{ \mathrm{optimal\ cost\ initialized\ from}\ \vecf{x}_0 }{ \mathrm{optimal\ cost\ initialized\ from}\ \boldsymbol{\xi}_{\matf{A}} }
,
\end{equation*}
to benchmark the influence of different initializations.
If the optimality metric is close to $1$, then the optimal cost is close to the one obtained from the eigenvector initialization,
and otherwise if this metric deviates from $1$ then the computed solution is considered to be suboptimal.
Intuitively, the optimality curve defines the robustness against bad initializations.
The statistics for each tested case are plotted in Fig.~\ref{fig. Proxy step-size with adaptive search and Lipschitz constant. Optimality of the final costs}.
It is worth noting that the costs of the $\ell_1$ norm regularized instances are negative (as $\matf{A} = - \matf{M}$), thus the optimality metric is below $1$.
For the problem instances used in this paper, the PGS, A-PGS and AM-PGS methods are in general robust to a large range of initializations,
while the $\ell_1$ norm regularized instances are more sensitive to initializations.

\subsection{Comparison with ManPG \cite{chen2020proximal} and AManPG \cite{huang2021extension}}

\begin{figure}[t]
	\centering
	\includegraphics[width=0.46\textwidth]{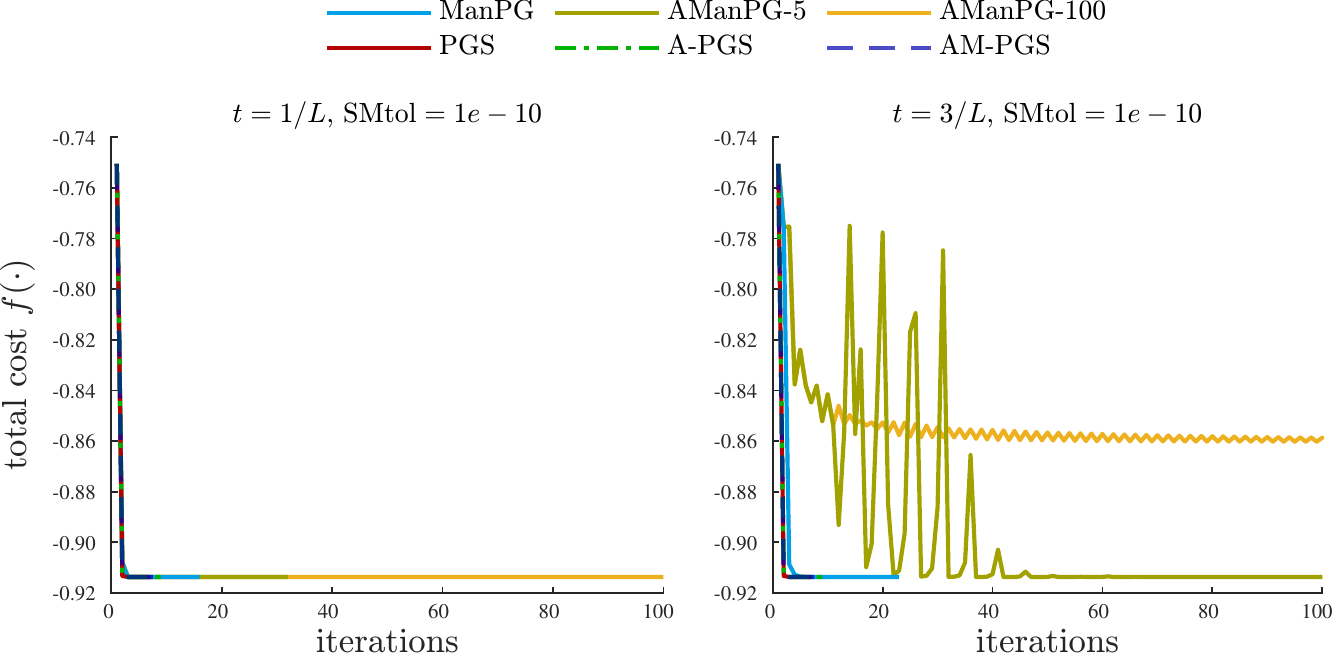}
	\caption{The convergence of ManPG and AManPG using the line-search criterion (\ref{eq. line-search in Chen and Huang works}) in comparison to the proposed one.}
	\label{fig. numerical example of comparing PGS and ManPG convergence}
\end{figure}
\begin{figure}[t]
	\centering
	\includegraphics[width=0.46\textwidth]{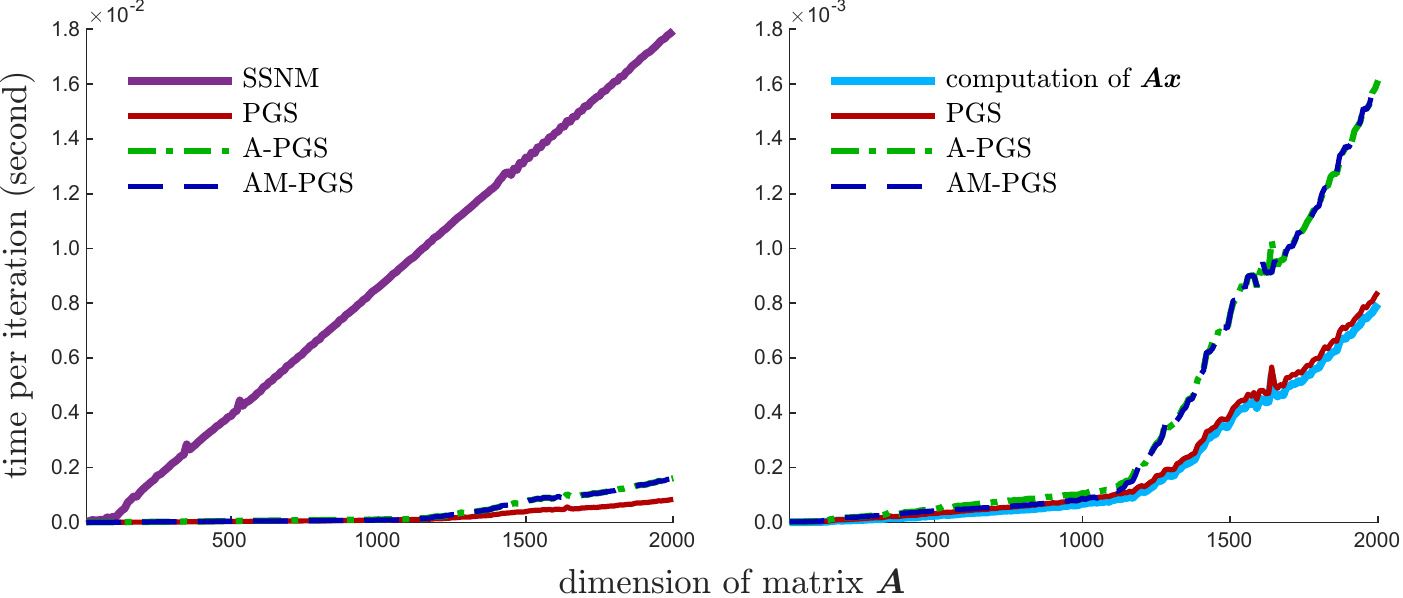}
	\caption{The computational time per iteration by using the SSNM method \cite{xiao2018regularized} (used in \cite{chen2020proximal,huang2021extension,huang2021riemannian}) and the proposed proxy step-size technique (used in the PGS, A-PGS and AM-PGS methods). The matrix-vector multiplication $\boldsymbol{A} \boldsymbol{x}$ is computed by the level-2 BLAS routine ``dsymv".}
	\label{fig. numerical example of comparing the timing of ssnm and proxy-stepsize techniques}
\end{figure}

We compare with the ManPG method \cite{chen2020proximal} and its accelerated version AManPG \cite{huang2021extension}.
Due to the difficulties of implementing the generalized Clarke differential for nuclear norm regularization, the comparison is only performed by applying $\ell_1$ regularization to problem (\ref{eq: ralyleight quotient optimization}).
We use the C++ implementation of ManPG and AManPG released in \cite{huang2021extension}.
AManPG uses a safeguard mechanism every $N$ iterations, where we set $N=5$ and $N=100$ and term the resulting methods as AManPG-5 and AManPG-100.
The error-tolerance of the SSNM method~\cite{xiao2018regularized} is set to $1e-10$.

We implemented our methods PGS, A-PGS, and AM-PGS in C++ as well for a fair comparison.


\begin{table*}[t]
	\centering
	\begin{tabular}{@{}  c | ccccc | ccccc @{} }
		\cline{2-11} 
		& \multicolumn{5}{c|}{$e_{dist}$ (pixels)$\downarrow$} & \multicolumn{5}{c}{$e_{rep}$ (pixels)$\downarrow$} \\ 
		& \textsf{8pt} & \textsf{PGS}$_5$ & \textsf{PGS}$_{10}$ &  \textsf{PGS} & \textsf{Gp} & \textsf{8pt} & \textsf{PGS}$_5$ & \textsf{PGS}$_{10}$ & \textsf{PGS} & \textsf{Gp} \\ 
		\hline
		Chapel(0,1) & 0.386 & 0.380 & \textit{0.380} & \underline{0.379} & \textbf{0.376} & 0.260 & 0.256 & \textit{0.256} & \underline{0.255} & \textbf{0.254} \\ 
		Keble(0,3) & 0.248 & 0.247 & \textit{0.247} & \underline{0.247} & \textbf{0.247} & 0.175 & 0.175 & \textit{0.175} & \underline{0.175} & \textbf{0.175} \\ 
		Desktop(C,D) & 0.574 & 0.336 & \textit{0.303} & \underline{0.288} & \textbf{0.268} & 0.406 & 0.238 & \textit{0.214} & \underline{0.204} & \textbf{0.190} \\ 
		Library(1,3) & 0.428 & 0.409 & \textit{0.405} & \underline{0.405} & \textbf{0.400} & 0.302 & 0.289 & \textit{0.286} & \underline{0.286} & \textbf{0.282} \\ 
		Merton1(1,3) & 0.308 & 0.295 & \textit{0.291} & \underline{0.288} & \textbf{0.277} & 0.217 & 0.209 & \textit{0.205} & \underline{0.203} & \textbf{0.196} \\ 
		Merton2(1,3) & 0.596 & 0.528 & \textit{0.498} & \underline{0.472} & \textbf{0.404} & 0.421 & 0.373 & \textit{0.352} & \underline{0.334} & \textbf{0.286} \\ 
		Arch & 0.304 & 0.299 & \textit{0.299} & \underline{0.299} & \textbf{0.298} & 0.215 & 0.211 & \textit{0.211} & \underline{0.211} & \textbf{0.211} \\ 
		Yard & 0.433 & 0.429 & \textit{0.429} & \underline{0.428} & \textbf{0.426} & 0.306 & 0.303 & \textit{0.303} & \underline{0.302} & \textbf{0.301} \\ 
		Slate & 0.246 & 0.184 & \textit{0.177} & \underline{0.170} & \textbf{0.163} & 0.129 & 0.097 & \textit{0.093} & \underline{0.089} & \textbf{0.085} \\ 
		Ben1 & 0.203 & 0.144 & \textit{0.128} & \underline{0.102} & \textbf{0.101} & 0.142 & 0.101 & \textit{0.089} & \underline{0.071} & \textbf{0.070} \\ 
		Ben2 & 0.139 & 0.086 & \textit{0.063} & \underline{0.050} & \textbf{0.048} & 0.097 & 0.060 & \textit{0.044} & \underline{0.035} & \textbf{0.033} \\ 
		\hline
	\end{tabular}
	\caption{Fundamental matrix estimation with nuclear norm regularization (with data in \cite{bugarin2015rank}).}
	\label{table: the error statistics of different methods on different dataset}
\end{table*}

\subsubsection{Line-search Criteria and Convergence}

\noindent\textbf{Theoretical justification.}
In ManPG's line-search criterion (\ref{eq. line-search in Chen and Huang works}), we first need to assign $t$ and then perform line-search for $\alpha_{k}$ to ensure the descent of the total cost $f(\cdot)$.
From Theorem \ref{theorem: the pillar reslut for convergence proofs} of our work, we see that if $t \le 1/L$ with $L$ being the Lipschitz constant, it suffices to set $\alpha_k = 1$ in the ManPG's line-search criterion (\ref{eq. line-search in Chen and Huang works}).
The authors in \cite{chen2020proximal, huang2021extension} assume known Lipschitz constant $L$ and suggest to use $t = 1/L$ as a reference.
Although for arbitrary $t > 0$, the existence of $\alpha_k$ for criterion (\ref{eq. line-search in Chen and Huang works}) is proved in \cite{chen2020proximal}, we observe a proper $t$ in ManPG/AManPG is required.

\vspace{5pt}\noindent\textbf{Numerical validation.}
On the left of Fig.~\ref{fig. numerical example of comparing PGS and ManPG convergence},
we see by setting $t = 1/L$, ManPG's line-search (\ref{eq. line-search in Chen and Huang works}) works almost the same way as the proposed PGS line-search.
On the right of Fig.~\ref{fig. numerical example of comparing PGS and ManPG convergence},
we run the same numerical instance again by setting $t = 3/L$ in ManPG and AManPG.
In this case, ManPG converges slower as seen from its curve being slightly shifted right.
With some fluctuations, the accelerated method AManPG-5 manages to converge while the convergence is even slower than ManPG. AManPG-10 simply does not converge.

In practice, if the Lipschitz constant $L$ is unknown, it is expected that an approximated $L$ may cause a lot of trouble in ManPG's line-search criterion (\ref{eq. line-search in Chen and Huang works}) as used in~\cite{chen2020proximal,huang2021extension}.
In the proposed line-search, this is never a problem. Intuitively, we fix $\alpha_{k} = 1$ in criterion (\ref{eq. line-search in Chen and Huang works}) and use proxy step-size $t'$ to find a working $t$ from the line-search criterion (\ref{eq: line-search equation based on t}).
This process is well-defined by the Lipschitz type assumption (\textit{i.e.,}~Assumption \ref{eq: Liptchtz type assumption on g cost}), and ensures the descent of the total cost $f(\cdot)$ by Theorem \ref{theorem: the pillar reslut for convergence proofs}.

\subsubsection{Computational Complexity}

The proposed PGS methods are much faster than ManPG and AManPG. Both ManPG and AManPG rely on the SSNM method \cite{xiao2018regularized} to solve the non-smooth KKT system, thus we compare the computation time per iteration of the proposed proxy step-size technique with that of the SSNM method.
We report the timing statistics per iteration in
Fig.~\ref{fig. numerical example of comparing the timing of ssnm and proxy-stepsize techniques}, with respect to the dimension of the $\matf{A}$ matrix in problem (\ref{eq: ralyleight quotient optimization}).
It is clearly seen that the proposed proxy step-size technique is substantially faster than the SSNM method.
An ablation study shows that the computation time per iteration of the PGS, A-PGS and AM-PGS methods is mostly decided by the matrix-vector multiplication $\matf{A} \vecf{x}_k$ used in evaluating the Euclidean gradient $\nabla g \vert_{\vecf{x}} = 2 \matf{A} \vecf{x}_k$ and the cost function $g (\vecf{x}_k) = \vecf{x}_k^{\tran} \matf{A}  \vecf{x}_k$.
The A-PGS and AM-PGS methods have the same per iteration complexity, while being twice more expensive than the PGS method due to an extra evaluation at the auxiliary state $\vecf{y}_k$.

\subsection{Fundamental Matrix Estimation}

We find $\lambda = 0.01$ works well in general,
after normalizing the image points \cite{Hartley2004}.
For this problem, we find that the PGS algorithm converges mostly within $20$ iterations.
For most of the cases, $10$ iterations or even $5$ are sufficient to reduce the last singular value $\sigma_3$ close enough to zero.
Therefore aside from the PGS method (with full convergence), aiming for an efficient engineering design, we propose the following two truncated PGS algorithms:
\begin{itemize}
	\item[-]PGS$_{5}$. $5$ PGS iterations, followed by a rank-$2$ rounding by setting $\sigma_3 = 0$.
	\item[-]PGS$_{10}$. $10$ PGS iterations, followed by a rank-$2$ rounding by setting $\sigma_3 = 0$.
\end{itemize}

We use two benchmark algorithms: a) the normalized eight point algorithm (denoted as 8pt), b) the global polynomial optimization \cite{bugarin2015rank} (denoted as Gp) with a formulation based on $\det(\matf{F}) = 0$.

We run the 8pt, PGS$_{5}$, PGS$_{10}$, PGS and Gp methods on a list of standard benchmarks,
and report in Table \ref{table: the error statistics of different methods on different dataset}
a) $e_{dist}$ the distance between the epipolar line and the corresponding image feature point \cite{zhang1998determining},
and b) $e_{rep}$ the reprojection error of the triangulated 3D points.
Overall, the $e_{dist}$ and $e_{rep}$ statistics decrease consistently over the 8pt, PGS$_{5}$, PGS$_{10}$, PGS, and Gp methods.
Table \ref{table: the error statistics of different methods on different dataset} shows that
the PGS$_{5}$ and PGS$_{10}$ methods, with a close performance towards the PGS method, consistently outperform the 8pt method, 
and they can give almost the same accuracy as the global method Gp.

\begin{figure}[t]
	\centering
	\includegraphics[width=0.42\textwidth]{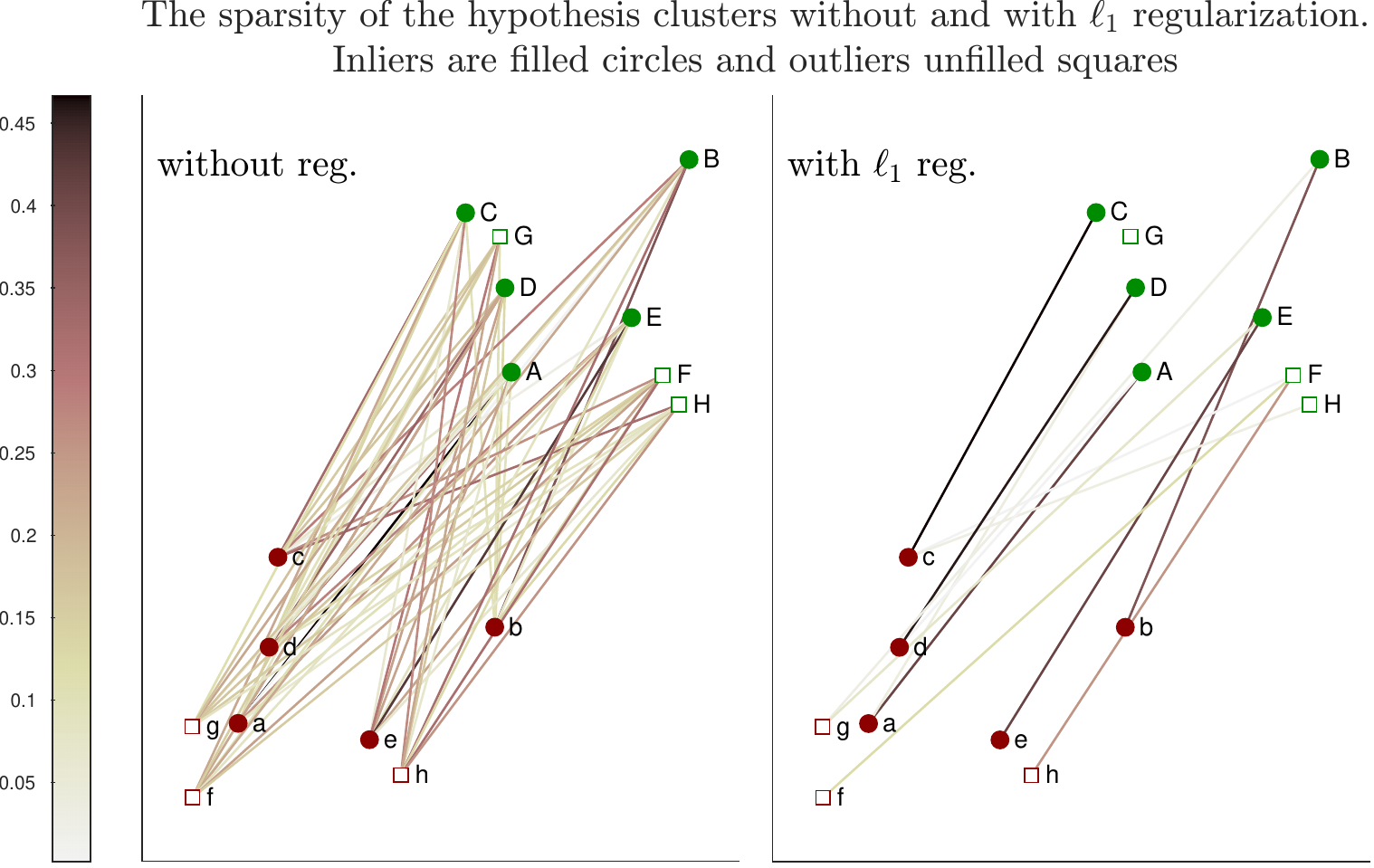}
	\\[10pt]
	\includegraphics[width=0.45\textwidth]{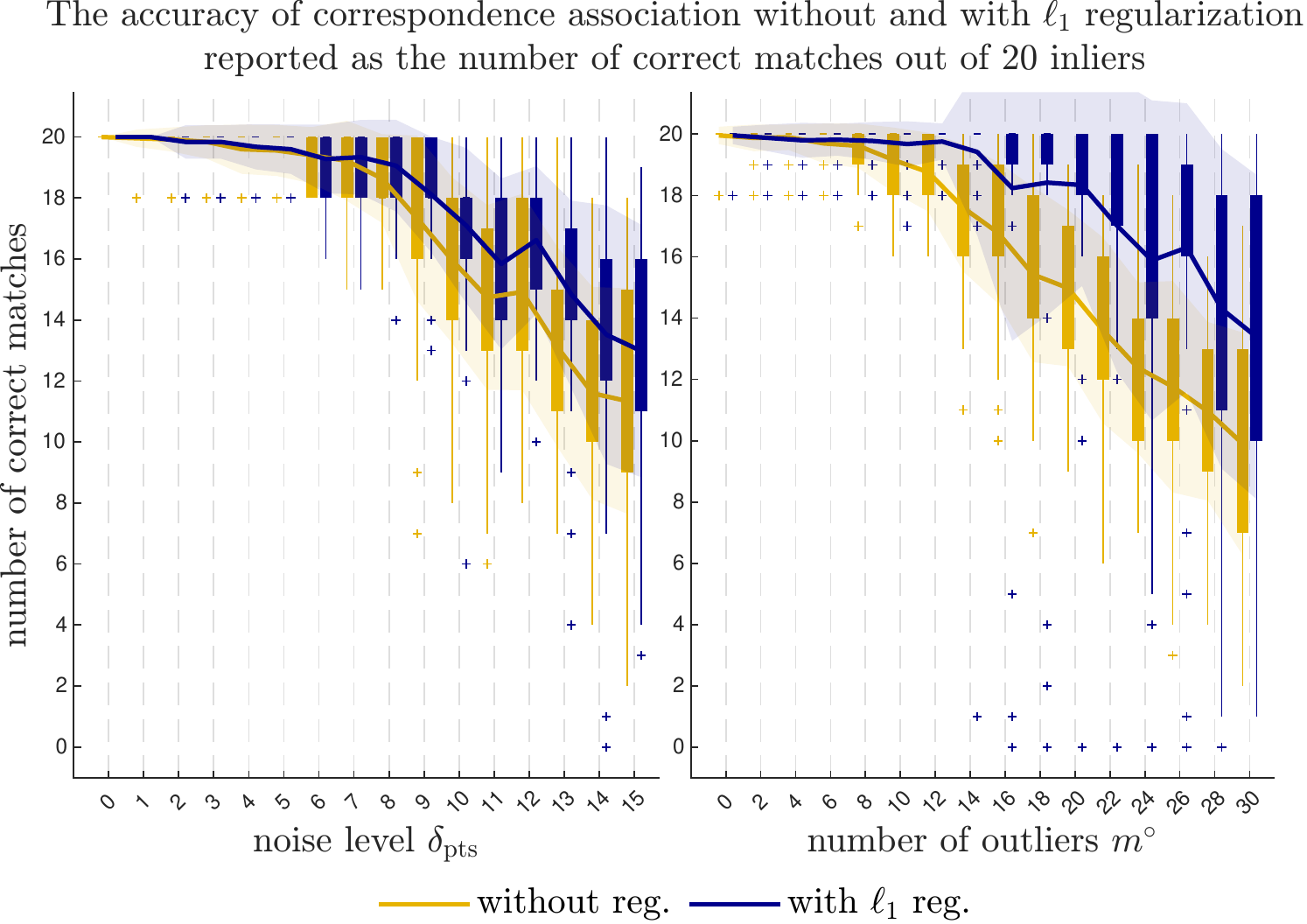}
	\caption{Correspondence association with $\ell_1$ norm regularization. With $\ell_1$ regularization, high probability hypotheses are enhanced while low probability hypotheses are trimmed off, resulting in a sparser and more consistent cluster of  association hypotheses.
		Therefore the correct correspondences (\textit{i.e.,}~inliers) can be identified more robustly.}
	\label{fig. graph assoication hypothesis visualization and mc simulation}
\end{figure}

\subsection{Correspondence Association}

The regularization strength is set as
$
\lambda = - \sigma_{\min} (\matf{A}) / (\sqrt{n} - 1)
=
\sigma_{\max} (\matf{M}) / (\sqrt{n} - 1)
$ with $n$ the dimension of $\vecf{A}$.
This choice is motivated by the canonical basis vector $\vecf{e}_k^{\tran}$.
Since the diagonal elements of $\matf{A}$ are zero by construction, the total cost at $\vecf{e}_k$ is
$
f (\vecf{e}_k)
= 
\vecf{e}_k \matf{A} \vecf{e}_k^{\tran}
+
\lambda
\Vert \vecf{e}_k \Vert_{\ell_1}
=
\lambda
$.
Noting that
$
1 \le \Vert \vecf{x} \Vert_{\ell_1} \le \sqrt{n}
$,
we thus have the following relation:
\begin{align*}
f(\vecf{x}) < f (\vecf{e}_k)
& \Leftrightarrow
\vecf{x}^{\tran} \matf{A}  \vecf{x}
+
\lambda
\Vert \vecf{x} \Vert_{\ell_1}
\le
\vecf{x}^{\tran} \matf{A}  \vecf{x}
+
\lambda \sqrt{n}
< \lambda
\\ & \Leftrightarrow
\lambda
\le
-
\frac{\vecf{x}^{\tran} \matf{A}  \vecf{x}}{\sqrt{n}-1}
\le
-
\frac{\sigma_{\min}(\matf{A})}{\sqrt{n}-1}
,
\end{align*}
where the used $\lambda$ is chosen as the largest possible value.

Our experimental setup is similar to the one used in Section 5.1 of \cite{leordeanu2005spectral}.
We simulate a 2-dimensional point-cloud $\mathcal{Q}$ of $m$ points, with $x-y$ coordinates uniformly distributed in $[0,\, 256 \sqrt{m/10}]$.
Then we add zero-mean Gaussian noise with standard deviation $\delta_{\mathrm{pts}}$ to each point in $\mathcal{Q}$, and rotate and translate the whole $\mathcal{Q}$ to obtain $\mathcal{Q}'$.
We generate $m^{\circ}$ outliers in both $\mathcal{Q}$ and $\mathcal{Q}'$ uniformly in the same region.
We use
$\delta_d = 5$ in equation (\ref{eq. the adjacency matrix of grpah association}) as in \cite{leordeanu2005spectral}.

The results are reported in Fig.~\ref{fig. graph assoication hypothesis visualization and mc simulation}.
In the first figure, we give an illustration (using $m=5$) that with $\ell_1$ norm regularization, many unlikely hypotheses are trimmed off, thus yielding a sparse hypothesis cluster.
In the second figure, we set $m = 20$ and report the number of correct matches with respect to different noise levels
and with respect to different outlier ratios.
For each case, we use a $50$ run Monte-Carlo simulation.
It is clear that the association accuracy is consistently improved by using $\ell_1$ norm regularization.

\subsection{Linear Self-calibration}

\begin{figure}[t]	
	\centering
	\includegraphics[width = 0.45\textwidth]{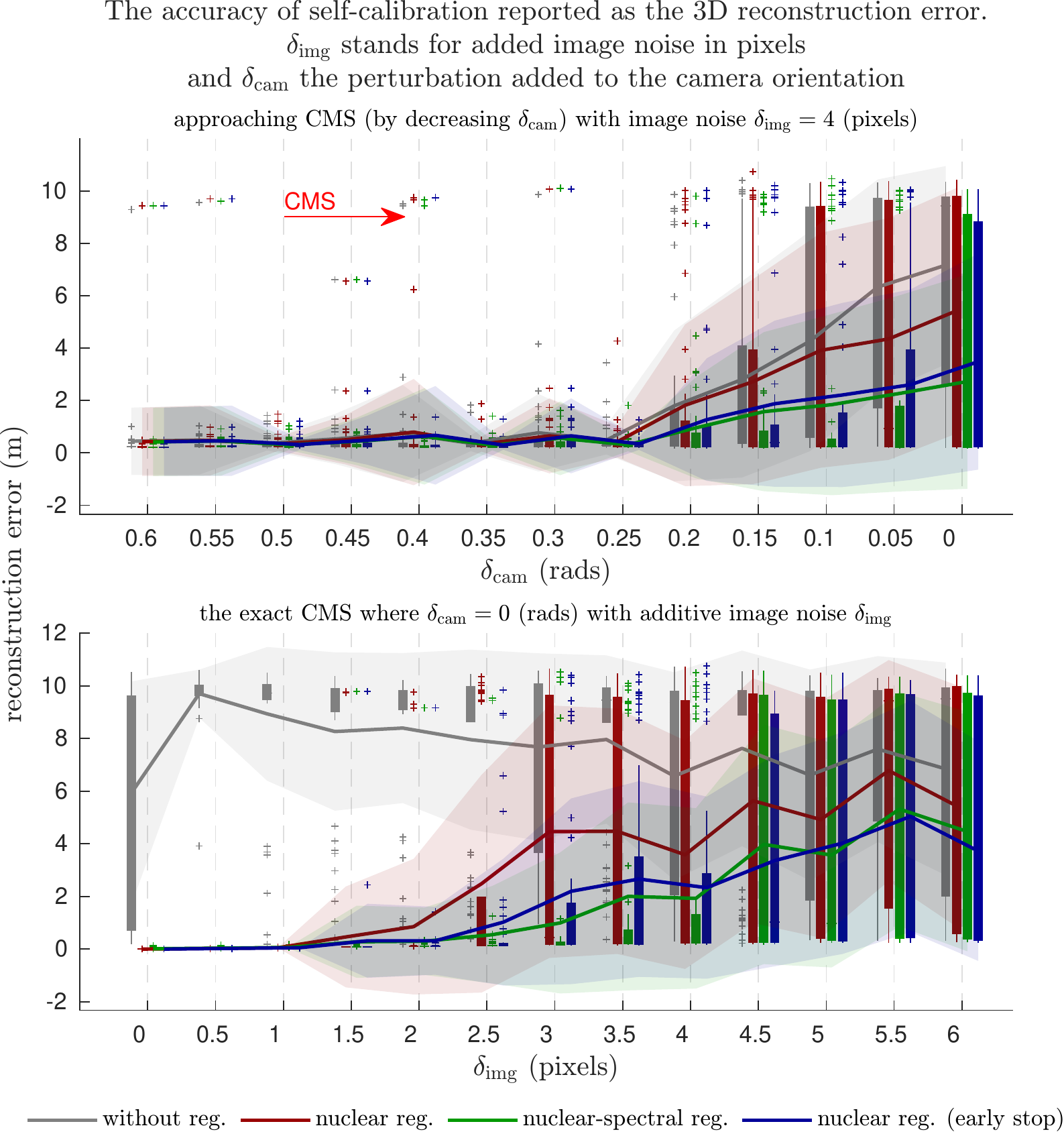}
	\\[10pt]
	\includegraphics[width=0.42\textwidth]{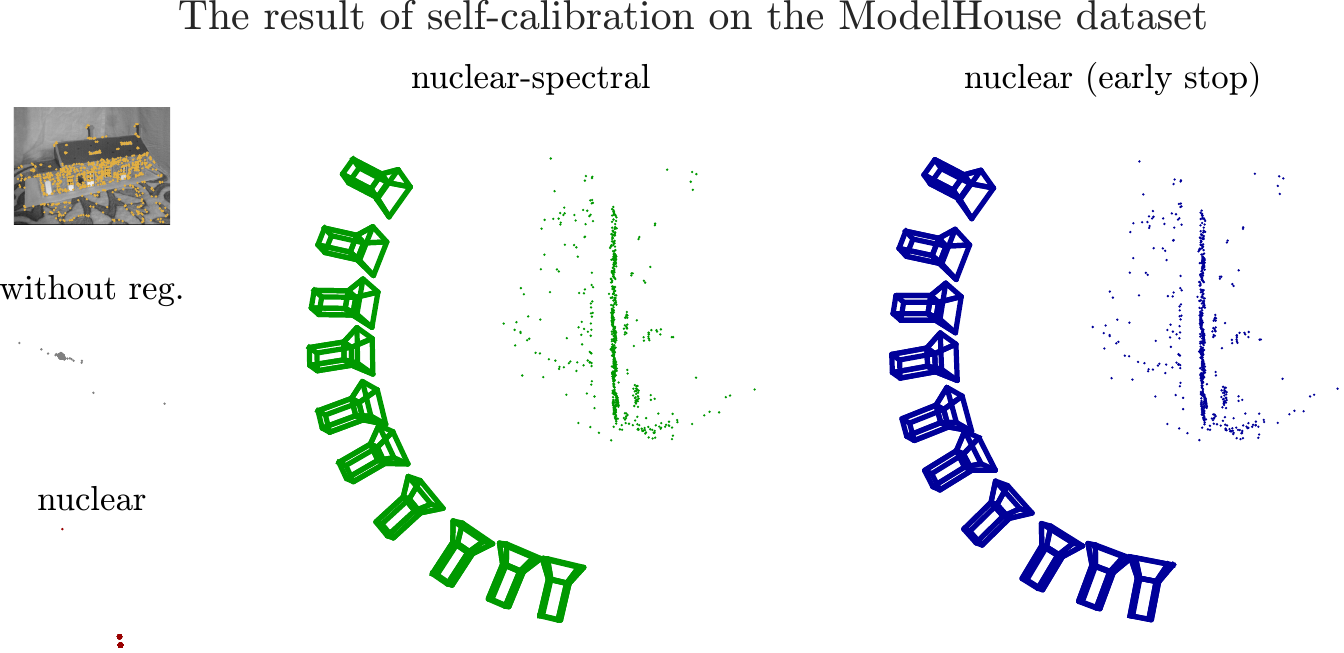}
	\caption{Self-calibration with the nuclear norm and the nuclear-spectral norm regularization.
		We show the results with the critical motion sequence (CMS) where the camera principal axes pass approximately through the geometric center of the observed object (case R4 of Table~$1$ in \cite{gurdjos2009dual}).}
	\label{fig. ther result of auto-calibration on simulation data}
\end{figure}

\subsubsection{Implementation}

We find normalization is essential to obtain stable self-calibration results.
The key points of our implementation are sketched as follows:
\begin{enumerate}
\item Image point normalization.
Let $(p_x, p_y)$ be the principal point of the camera.
If this is unknown, we approximate $(p_x, p_y) = (I_x/2, I_y/2)$, where $I_x$ and $I_y$ are the width and height of the image.
The average distance of all image points to $(p_x, p_y)$ is denoted by $s$.
We normalize all image points by a common transformation $\matf{T}^{-1}$:
\begin{equation*}
\matf{T}
= 
\begin{pmatrix}
s & 0 &  p_x \\ 
0  &  s & p_y \\ 
0 & 0 & 1
\end{pmatrix}
, \ 
\matf{T}^{-1} 
= 
\begin{pmatrix}
1/s & 0 & -{p_x}/{s} \\ 
0  &  1/s & - {p_y}/{s} \\ 
0 & 0 & 1
\end{pmatrix}
.
\end{equation*}
\item Projective reconstruction using projective bundle-adjustment from the normalized image points.
\item Quasi-Euclidean rectification
\cite{beardsley1997sequential}.
We approximate the intrinsic matrix of each camera computed in step $2)$ as
$\matf{K}_i = \mathrm{diag}(f_i, f_i, 1)$,
with $f_i = 2\sqrt{m_x^2 + m_y^2}$ where $m_x$ and $m_y$ are the maximum range in the $x$- and $y$- coordinates.
Using this approximate $\matf{K}_i$,
we compute an approximate estimate of $\matf{Q}_{\infty}^{\ast}$ from the DAQ constraint (\ref{eq: key property used for self-calibration}) \textit{i.e.,}~$
\tilde{\matf{P}}_i \matf{Q}_{\infty}^{\ast} \tilde{\matf{P}}_i^{\tran}
\propto
\matf{K}_i \matf{K}_i^{\tran}
$,
and rectify the projective reconstruction approximately.
\item Linear self-calibration from the quasi-Euclidean rectification in step $3)$.
We initialize the PGS methods from the canonical DAQ $\Omega_{\infty}^{\ast} = \mathrm{diag}(1,1,1,0)$.
\item Transforming cameras $\matf{P}_i$ obtained from step $4)$ by $\matf{T}$. Lastly, $\matf{T} \matf{P}_i$ are the final estimate of Euclidean cameras.
\end{enumerate}
We use $\lambda = 0.01$ for the nuclear norm regularization, and $\lambda_1 = 0.01$, $\lambda_2 = 2 \lambda_1$ for the nuclear-spectral norm regularization.

\subsubsection{Simulated Data}

We simulate a scene comprising: a) $50$ points spreading randomly in a diameter of $3$ meters;
b) $7$ cameras circularly distributed $30$ meters away from the point-cloud.
All cameras are oriented towards the centroid of the point-cloud,
thus the simulated scenario is an artificial CMS whose ambiguity can be removed using the rank deficiency of the DAQ.
We evaluate the performance of each method with respect to the perturbation of camera orientations $\delta_{\mathrm{cam}}$ and the noise of image points $\delta_{\mathrm{img}}$.
We use the 3D reconstruction error as the evaluation metric,
which is computed by the similarity Procrustes analysis between the ground-truth point-cloud and the estimated point-cloud.

We conduct two sets of experiments.
First, we use the fixed image noise $\delta_{\mathrm{img}}=4$ pixels
and decrease $\delta_{\mathrm{cam}}$ to gradually bring the camera configuration to the CMS.
Second, we set the camera configuration to the exact CMS where $\delta_{\mathrm{cam}} = 0$ and then test the performance with respect to different image noise $\delta_{\mathrm{img}}$.
We report the reconstruction error with a $50$ run Monte-Carlo simulation in Fig.~\ref{fig. ther result of auto-calibration on simulation data}.

As shown in Fig.~\ref{fig. ther result of auto-calibration on simulation data}, when facing the CMS, the classical method without regularization fails,
and the method with nuclear-spectral norm regularization is more robust than the one with nuclear norm regularization \textit{e.g.,}~in case of the exact CMS where $\delta_{\mathrm{cam}} = 0$ and $\delta_{\mathrm{img}} > 3$.
It is interesting to see that the nuclear-norm regularization works well for less noisy scenarios of the CMS, \textit{e.g.,}~when $\delta_{\mathrm{img}} < 2$. For these cases, it seems that the solution of the nuclear-norm regularized problem is well-trapped at the local minimum (the rank-$3$ DAQ), while the gradient is not large enough to go to the global minimum (the rank-$1$ DAQ).
To examine this hypothesis, we use an early-stop trick (by setting the maximum PGS iterations to $1000$) in the nuclear norm regularization, and observe that with the early-stop trick the nuclear-norm regularized method mostly performs well.

\subsubsection{Real Data}

A qualitative example of the CMS is given in
Fig.~\ref{fig. ther result of auto-calibration on simulation data} using a real dataset called ModelHouse
where all cameras look towards the geometric center of a model house\footnote{\url{https://www.robots.ox.ac.uk/~vgg/data/mview/}}.
In this scenario, matrix $\matf{A}$ has two eigenvalues close to zero.
The classic method fails because the solution space is ambiguous.
The nuclear norm regularized method fails by converging to the rank-$1$ solution.
The nuclear-spectral norm regularized method can find the correct rank-$3$ solution thus recover the correct Euclidean geometry.
The nuclear norm regularized method with the early stop trick also works.
Intuitively, the spectral norm in the nuclear-spectral norm regularization prevents the algorithm from gliding to the rank-$1$ solution, and the early stop trick has the similar functionality.

\section{Conclusion}
\label{section. conclusion}

We have proposed the proxy step-size technique, and presented an effective solution to problem (\ref{eq: the general form: optimization problem studied}) for convex and absolutely homogeneous $h(\cdot)$.
The proposed solution is: exact (satisfying the first-order necessary optimality condition), elegant (simple and in closed-form), and easily applicable (to nuclear norm regularization etc.).
Future work includes extending the proxy-step size technique to the oblique and the Stiefel manifolds, and analyzing the convergence rate of the accelerated methods in Algorithm~\ref{algorithm: a unified implementation of PGS, A-PGS, AM-PGS methods}.

\begin{IEEEbiography}[{\includegraphics[width=1in,height=1.25in,clip,keepaspectratio]{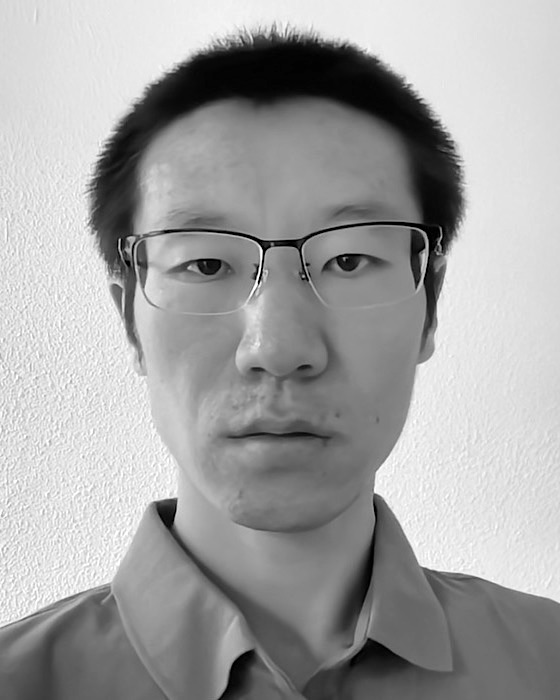}}]{Fang Bai}
Fang Bai was born in Ningxia Province, China, in 1988. He received the B.Sc. degree in computer science and technology from Nankai University, China, in 2010, and the Ph.D. degree in robotics from University of Technology Sydney, Australia, in 2020.
His research has been focused on mathematical abstractions in robotics and computer vision. He has conducted several fundamental breakthroughs on related topics as the first-author, e.g., the cycle based pose graph optimization, the equation to predict the change of optimal values, and the closed-form solution for template-free deformable Procrustes analysis.
\end{IEEEbiography}
\begin{IEEEbiography}[{\includegraphics[width=1in,height=1.25in,clip,keepaspectratio]{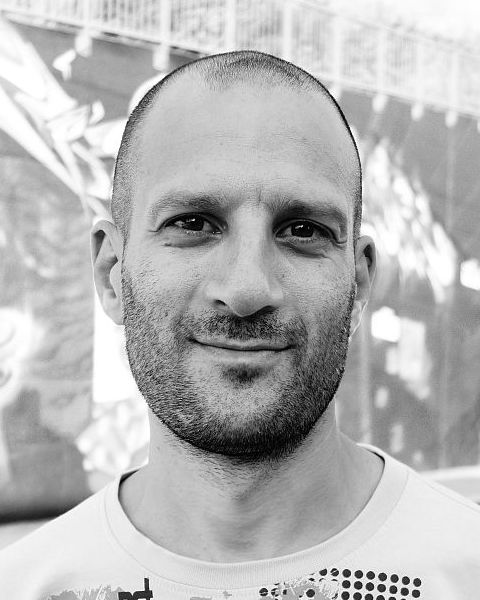}}]{Adrien Bartoli}
Adrien Bartoli has held the position of Professor of Computer Science at Université Clermont Auvergne since fall 2009 and has been a member of Institut Universitaire de France since 2016. He is currently on leave as research scientist at the University Hospital of Clermont-Ferrand and as Chief Scientific Officer at SurgAR. He leads the Endoscopy and Computer Vision (EnCoV) research group at the University and Hospital of Clermont-Ferrand. His main research interests are in computer vision, including image registration and Shape-from-X for deformable environments, and their application to computer-aided medical interventions.
\end{IEEEbiography}

%
%
%

\clearpage

\appendices

\section{Properties of proximal}
\label{appendix. properties of proximal}

\subsection{Convexity}

If $h(\cdot)$ is convex, then for any $\vecf{x}$ and $\vecf{x}_0$, we have:
\begin{equation*}
h (\vecf{x}) \ge  h(\vecf{x}_0) + 
\langle \partial{h}\vert_{\vecf{x}_0},\, \vecf{x} - \vecf{x}_0 \rangle
.
\end{equation*}

\subsection{Firm non-expansiveness and non-expansiveness}

For a convex $h(\cdot)$,
let $\vecf{z}_1 = \mathrm{prox}_{th} \left( \vecf{w}_1 \right)$ and $\vecf{z}_2 = \mathrm{prox}_{th} \left( \vecf{w}_2 \right)$ with $t\ge 0$. Then the following holds:
\begin{itemize}
\item Firm non-expansiveness:
\begin{equation}
\label{eq. standard Firm non-expansiveness}
\left\Vert \vecf{z}_1 - \vecf{z}_2 \right\Vert_2^2
\le
\langle
\vecf{z}_1 - \vecf{z}_2 ,\, \vecf{w}_1 - \vecf{w}_2
\rangle
;
\end{equation}
\item Non-expansiveness:
\begin{equation}
\label{eq. standard non-expansiveness}
\left\Vert \vecf{z}_1 - \vecf{z}_2 \right\Vert_2^2  \le  \left\Vert \vecf{w}_1 - \vecf{w}_2 \right\Vert_2^2
.
\end{equation}
\end{itemize}
\begin{proof}
The solution of the proximal is characterized by its first-order necessary condition:
\begin{equation}
\label{eq. first-order-optimality 1}
\vecf{0} \in \partial{h}\vert_{\vecf{z}_1} + \frac{1}{t}(\vecf{z}_1 - \vecf{w}_1)
\Leftrightarrow
\frac{1}{t}(\vecf{w}_1 - \vecf{z}_1) \in  \partial{h}\vert_{\vecf{z}_1}
,
\end{equation}
\begin{equation}
\label{eq. first-order-optimality 2}
\vecf{0} \in \partial{h}\vert_{\vecf{z}_2} + \frac{1}{t}(\vecf{z}_2 - \vecf{w}_2)
\Leftrightarrow
\frac{1}{t}(\vecf{w}_2 - \vecf{z}_2) \in  \partial{h}\vert_{\vecf{z}_2}
.
\end{equation}
By the convexity of $h(\cdot)$, we have:
\begin{subnumcases} {\label{eq. convexity of h}}
h (\vecf{z}_1) \ge  h(\vecf{z}_2) + 
\langle \partial{h}\vert_{\vecf{z}_2},\, \vecf{z}_1 - \vecf{z}_2 \rangle
\\[5pt]
h (\vecf{z}_2) \ge  h(\vecf{z}_1) + 
\langle \partial{h}\vert_{\vecf{z}_1},\, \vecf{z}_2 - \vecf{z}_1 \rangle
.
\end{subnumcases}
Substituting equations (\ref{eq. first-order-optimality 1}) and (\ref{eq. first-order-optimality 2}) into (\ref{eq. convexity of h}), we have:
\begin{subnumcases} {\label{eq. convexity and optimality}}
h (\vecf{z}_1) \ge  h(\vecf{z}_2) + 
\langle \frac{1}{t}(\vecf{w}_2 - \vecf{z}_2),\, \vecf{z}_1 - \vecf{z}_2 \rangle
\label{eq. convexity and optimality 1}
\\[2pt]
h (\vecf{z}_2) \ge  h(\vecf{z}_1) + 
\langle \frac{1}{t}(\vecf{w}_1 - \vecf{z}_1),\, \vecf{z}_2 - \vecf{z}_1 \rangle
\label{eq. convexity and optimality 2}
.
\end{subnumcases}
Summing together inequalities (\ref{eq. convexity and optimality 1}) and (\ref{eq. convexity and optimality 2}), we have:
\begin{equation*}
\langle \vecf{w}_2 - \vecf{z}_2 - \vecf{w}_1 + \vecf{z}_1,\, \vecf{z}_1 - \vecf{z}_2 \rangle
\le 0.
\end{equation*}
Expanding the above, we obtain the firm non-expansiveness (\ref{eq. standard Firm non-expansiveness}).
By the Cauchy–Schwarz inequality, we obtain:
\begin{equation*}
\left\vert
\langle
\vecf{z}_1 - \vecf{z}_2 ,\, \vecf{w}_1 - \vecf{w}_2
\rangle
\right\vert
\le
\left\Vert \vecf{z}_1 - \vecf{z}_2 \right\Vert_2
\left\Vert \vecf{w}_1 - \vecf{w}_2 \right\Vert_2
.
\end{equation*}
After canceling $\left\Vert \vecf{z}_1 - \vecf{z}_2 \right\Vert_2$ in the firm non-expansiveness (\ref{eq. standard Firm non-expansiveness}), we obtain non-expansiveness (\ref{eq. standard non-expansiveness}).

\end{proof}

\subsection{Proof of Lemma \ref{lemma. additional properties for convex and absolute homogeneous h}}

By definition of the proximal, we write:
\begin{equation}
\label{eq. problem. proximal at w = 0}
\mathrm{prox}_{ t h} \left(  \vecf{0} \right)
 = 
\argmin_{\vecf{x}} \left\{ t h \left( \vecf{x} \right) +
\frac{1}{2}\left\Vert \vecf{x}  \right\Vert^2  \right\}
,\quad
t \ge 0
\end{equation}
Since $h(\cdot)$ is convex and absolutely homogeneous, from Lemma~\ref{lemma. properties of convex and absolutely homogeneous h()}, we have $h \left( \vecf{x} \right) \ge 0$ and $h \left( \vecf{0} \right) = 0$.
Thus, we observe that the optimal cost of problem (\ref{eq. problem. proximal at w = 0}) is $0$ which is attained at $\vecf{x} = \vecf{0}$ as $h \left( \vecf{0} \right) = 0$ and $\left\Vert \vecf{x}  \right\Vert = 0$.
This concludes $\mathrm{prox}_{ t h} \left(  \vecf{0} \right) = \vecf{0}$ if $h(\cdot)$ is convex and absolutely homogeneous.
The inequalities
$\left\Vert \mathrm{prox}_{th} ( \vecf{w} )  \right\Vert_2^2  \le  \langle \mathrm{prox}_{th} ( \vecf{w} ),\, \vecf{w}  \rangle $ and $ \left\Vert \mathrm{prox}_{th} ( \vecf{w} ) \right\Vert_2  \le  \left\Vert \vecf{w} \right\Vert_2$ follow from firm non-expansiveness (\ref{eq. standard Firm non-expansiveness}) and non-expansiveness (\ref{eq. standard non-expansiveness}), by setting $\vecf{w}_1 = \vecf{w}$ and $\vecf{w}_2 = \vecf{0}$.

\section{Derivation of KKT System (\ref{eq: old KKT system in u and t})}
\label{appendix: derivation of the two optimality equations}

We write $\vecf{0} \in \partial \mathcal{L}_{\vecf{v}_k}$ as follow:
\begin{align*}
& \vecf{0} \in 
\mathrm{grad}\,g \vert_{\vecf{x}_k}
+ \frac{1}{t} \vecf{v}_k
+ \partial h \vert_{\vecf{x}_k + \vecf{v}_k}
+ \mu \vecf{x}_k
\\[3pt] & \Leftrightarrow  \left(1-\mu t\right) \vecf{x}_k
- t \mathrm{grad}\,g \vert_{\vecf{x}_k}
\in 
t \partial h \vert_{\vecf{x}_k + \vecf{v}_k}
+ \vecf{x}_k + \vecf{v}_k
\\[3pt] & \Leftrightarrow
\vecf{x}_k + \vecf{v}_k
=
\mathrm{prox}_{th}
\left(
\left(1-\mu t\right) \vecf{x}_k
- t \mathrm{grad}\,g \vert_{\vecf{x}_k}
\right)
.
\end{align*}
Therefore:
\begin{equation*}
\vecf{v}_k
=
\mathrm{prox}_{th}
\left(
\left(1-\mu t\right) \vecf{x}_k
- t \mathrm{grad}\,g \vert_{\vecf{x}_k}
\right)
- \vecf{x}_k
,
\end{equation*}
which is equation (\ref{eq: exact RPG iteration equation 1}).
Considering $\vecf{x}_k^{\tran} \vecf{v}_k = 0$ with equation (\ref{eq: exact RPG iteration equation 1}), we obtain equation (\ref{eq: exact RPG iteration equation 2}).

\section{Proof of Lemma \ref{proposition: a result of proximal operator on linear operators}}
\label{appendix: proof of lemma on proximal for absolutely scable functions}

By definition of the proximal, we write:
\begin{equation*}
\mathrm{prox}_{t h} \left( \alpha \vecf{w} \right)
 = 
\argmin_{\vecf{x}} \left\{ h \left( \vecf{x} \right) +
\frac{1}{2t}\left\Vert \vecf{x} - \alpha \vecf{w} \right\Vert^2  \right\}
\defeq
\vecf{z}
.
\end{equation*}
Since $h(\cdot)$ is absolutely homogeneous,
$
 h \left( \vecf{x} \right)
 =
 \left\vert\alpha\right\vert h \left( \frac{1}{\alpha}\vecf{x} \right)
$.
The above equation can thus be written as:
\begin{align*}
\vecf{z}
&  =
\argmin_{\vecf{x}} \left\{
\left\vert\alpha\right\vert h \left( \frac{1}{\alpha}\vecf{x} \right) +
\frac{\alpha^2}{2t}\left\Vert \frac{1}{\alpha} \vecf{x} - \vecf{w} \right\Vert^2 \right\}
\\[4pt] & =
\argmin_{ \vecf{x}} \left\{
h \left( \frac{1}{\alpha}\vecf{x} \right) +
\frac{1}{2 \frac{t}{\left\vert\alpha\right\vert}} \left\Vert \frac{1}{\alpha} \vecf{x} - \vecf{w} \right\Vert^2 
\right\}
\\[4pt]
& \Leftrightarrow
\frac{1}{\alpha} 
\vecf{z}
 = 	\mathrm{prox}_{\frac{t}{\left\vert\alpha\right\vert}h}
\left( \vecf{w} \right)
.
\end{align*}
Therefore
$
\vecf{z}
=
\mathrm{prox}_{t h} \left( \alpha \vecf{w} \right)
= 
\alpha \mathrm{prox}_{\frac{t}{\left\vert\alpha\right\vert}h}
\left( \vecf{w} \right)
$.

\section{Proof of Lemma \ref{lemma: the monotonicity of 1/t in 1/t'}}
\label{appendix: proof to theorem, relation between t and t', ineqiality elsion > 0}

From equation (\ref{eq: our exact RPG iteration equation 2}), we obtain:
\begin{align*}
\frac{1}{\phi(t')}
& =
\frac{1}{t}
=
\frac{1}{t'}
\vecf{x}_k^{\tran}
\mathrm{prox}_{\left\vert t' \right\vert h}
\left(
t'
\left(
\frac{1}{t'} \vecf{x}_k -
\mathrm{grad}\,g \vert_{\vecf{x}_k}
\right)
\right)
\\ & =
\vecf{x}_k^{\tran}
\mathrm{prox}_{h}
\left(
\frac{1}{t'} \vecf{x}_k -
\mathrm{grad}\,g \vert_{\vecf{x}_k}
\right)
.
\end{align*}
The last equality is due to Lemma \ref{proposition: a result of proximal operator on linear operators}.
We denote
$\vecf{w}_1 = \frac{1}{t'_1} \vecf{x}_k -
\mathrm{grad}\,g \vert_{\vecf{x}_k}
$
and
$
\vecf{w}_2 =
\frac{1}{t'_2} \vecf{x}_k - \mathrm{grad}\,g \vert_{\vecf{x}_k}
$.
Then
$
{1}/{t'_1} = \vecf{x}_k^{\tran} \vecf{w}_1
$
and
$
{1}/{t'_2} = \vecf{x}_k^{\tran} \vecf{w}_2
$.
Moreover:
\begin{equation}
\vecf{x}_k \vecf{x}_k^{\tran} (\vecf{w}_1 - \vecf{w}_2)=
\frac{1}{t'_1} \vecf{x}_k - \frac{1}{t'_2} \vecf{x}_k
= \vecf{w}_1 - \vecf{w}_2
.
\end{equation}
The proof is given by summarizing the above facts as:
\begin{align*}
\epsilon(t'_1, t'_2)  & =
\langle
\vecf{x}_k^{\tran}
\mathrm{prox}_{h}
\left(
\vecf{w}_1
\right)
-
\vecf{x}_k^{\tran}
\mathrm{prox}_{h}
\left(
\vecf{w}_2
\right)
,\,
\vecf{x}_k^{\tran}
\left( \vecf{w}_1 - \vecf{w}_2 \right)
\rangle
\\ & =
\langle
\mathrm{prox}_{h}
\left(
\vecf{w}_1
\right)
-
\mathrm{prox}_{h}
\left(
\vecf{w}_2
\right)
,\,
\vecf{x}_k \vecf{x}_k^{\tran}
\left( \vecf{w}_1 - \vecf{w}_2 \right)
\rangle
\\ & =
\langle
\mathrm{prox}_{h}
\left(
\vecf{w}_1
\right)
-
\mathrm{prox}_{h}
\left(
\vecf{w}_2
\right)
,\,
\vecf{w}_1 - \vecf{w}_2
\rangle
\\ & \ge
\left\Vert
\mathrm{prox}_{h}
\left(
\vecf{w}_1
\right)
-
\mathrm{prox}_{h}
\left(
\vecf{w}_2
\right)
\right\Vert_2^2,
\end{align*}
where the inequality is due to the firm non-expansiveness of $\mathrm{prox}_{h}\left(\cdot\right)$.

\section{Proof of Lemma \ref{eq: the general inequality about proximal of y = prox(w), for absolute homogeneous h}}
\label{appendix. Proof of the global inequality for homogenous h}

\subsection{A General Inequality by Convexity}

\begin{lemma}
	\label{lemma: intermediate result: the general inequality about proximal of y = prox(w)}
	Let $h(\cdot)$ be convex.
	Then for any $t$, $\vecf{x}$ and $\vecf{w}$, we have:
	\begin{multline}
	\label{eq: the general inequality about proximal of y = prox(w)}
	\langle
	\vecf{w} - 
	\mathrm{prox}_{\left\vert t \right\vert h}
	\left( \vecf{w} \right)  ,\,
	\mathrm{prox}_{\left\vert t \right\vert h}
	\left( \vecf{w} \right) - \vecf{x}
	\rangle
	\\ \ge 
	\left\vert t \right\vert
	\left( h (\mathrm{prox}_{\left\vert t \right\vert h}
	\left( \vecf{w} \right)) - h (\vecf{x}) \right)
	.
	\end{multline}
\end{lemma}
\begin{proof}
	Let $\vecf{z} = \mathrm{prox}_{\left\vert t \right\vert h}
	\left( \vecf{w} \right)$.
	By definition of the proximal, we write:
	\begin{equation}
	\label{eq: proximal of y = prox(w)}
	\vecf{z} = \argmin_{\vecf{y}} \left\{
	h(\vecf{y}) + \frac{1}{2 \left\vert t \right\vert}
	\left\Vert
	\vecf{y} - \vecf{w}
	\right\Vert^2
	\right\}
	.
	\end{equation}
	The first-order necessary condition of problem (\ref{eq: proximal of y = prox(w)}) states:
	\begin{equation}
	\label{eq: proof: the first order necessary condition of proximal of y = prox(w)}
	- \frac{1}{\left\vert t \right\vert}
	\left(
	\vecf{z} - \vecf{w}
	\right)
	\in
	\partial h\vert_{\vecf{z}}
	.
	\end{equation}
	By the convexity of $h (\vecf{\cdot})$ at $\vecf{z}$, the following inequality holds for any $\vecf{x}$:
	\begin{equation*}
	h (\vecf{x}) \ge 
	h (\vecf{z}) + \langle
	\partial h\vert_{\vecf{z}},\,
	\vecf{x} - \vecf{z}
	\rangle
	.
	\end{equation*}
	Therefore for any $\vecf{x}$, we have the following inequality:
	\begin{equation*}
	h (\vecf{x}) \ge 
	h (\vecf{z}) 
	- \frac{1}{\left\vert t \right\vert}
	\langle
	\vecf{z} - \vecf{w} ,\,
	\vecf{x} - \vecf{z}
	\rangle
	.
	\end{equation*}
	Reorganizing this inequality, we obtain inequality (\ref{eq: the general inequality about proximal of y = prox(w)}).
\end{proof}

\subsection{Proof of Lemma \ref{eq: the general inequality about proximal of y = prox(w), for absolute homogeneous h}}

	Following Lemma \ref{lemma: intermediate result: the general inequality about proximal of y = prox(w)},
	in inequality (\ref{eq: the general inequality about proximal of y = prox(w)}), we let $\vecf{x} = \alpha \mathrm{prox}_{\left\vert t \right\vert h}
	\left( \vecf{w} \right)$:
	\begin{multline}
	\langle
	\vecf{w} - 
	\mathrm{prox}_{\left\vert t \right\vert h}
	\left( \vecf{w} \right)  ,\,
	(1 - \alpha)
	\mathrm{prox}_{\left\vert t \right\vert h}
	\left( \vecf{w} \right)
	\rangle
	\\ \ge 
	\left\vert t \right\vert (1 - \left\vert \alpha \right\vert)
	h (\mathrm{prox}_{\left\vert t \right\vert h}
	\left( \vecf{w} \right))
	.
	\end{multline}
	If $0 \le \alpha < 1$, then $1 - \left\vert \alpha \right\vert = 1 - \alpha > 0$:
	\begin{equation}
	\label{eq: proof: intermidate reulst on proximal of absolute homogenous function, ineq 1}
	\langle
	\vecf{w} - 
	\mathrm{prox}_{\left\vert t \right\vert h}
	\left( \vecf{w} \right)  ,\,
	\mathrm{prox}_{\left\vert t \right\vert h}
	\left( \vecf{w} \right)
	\rangle
	\ge 
	\left\vert t \right\vert
	h (\mathrm{prox}_{\left\vert t \right\vert h}
	\left( \vecf{w} \right))
	.
	\end{equation}
	If $\alpha > 1$, then
	$1 - \left\vert \alpha \right\vert = 1 - \alpha < 0$:
	\begin{equation}
	\label{eq: proof: intermidate reulst on proximal of absolute homogenous function, ineq 2}
	\langle
	\vecf{w} - 
	\mathrm{prox}_{\left\vert t \right\vert h}
	\left( \vecf{w} \right)  ,\,
	\mathrm{prox}_{\left\vert t \right\vert h}
	\left( \vecf{w} \right)
	\rangle
	\le 
	\left\vert t \right\vert
	h (\mathrm{prox}_{\left\vert t \right\vert h}
	\left( \vecf{w} \right))
	.
	\end{equation}
	By considering inequalities (\ref{eq: proof: intermidate reulst on proximal of absolute homogenous function, ineq 1}) and (\ref{eq: proof: intermidate reulst on proximal of absolute homogenous function, ineq 2}) together,
	we obtain:
	\begin{equation}
	\label{eq: result: tight relation of proximal, <w - y, y> = t h(y), for convex and homogeneous h}
	\langle
	\vecf{w} - 
	\mathrm{prox}_{\left\vert t \right\vert h}
	\left( \vecf{w} \right)   ,\,
	\mathrm{prox}_{\left\vert t \right\vert h}
	\left( \vecf{w} \right)
	\rangle
	=
	\left\vert t \right\vert h ( \mathrm{prox}_{\left\vert t \right\vert h}
	\left( \vecf{w} \right) )
	.
	\end{equation}	
	Subtracting equation (\ref{eq: result: tight relation of proximal, <w - y, y> = t h(y), for convex and homogeneous h}) in inequality (\ref{eq: the general inequality about proximal of y = prox(w)}),
	we obtain inequality (\ref{eq: result: key inequality <w - y, x> = t h(x), for convex and homogeneous h}).

\section{Proof of Proposition \ref{proposition. result on t neq 0 for t' neq 0}}
\label{appendix. proof of results on t neq 0 for t' neq 0}

It can be shown that $t = t' / c(t')$ satisfies:
\begin{align*}
\vert t \vert 
& = \frac{\left\vert t' \right\vert}{\left\vert \vecf{x}_k^{\tran} \mathrm{prox}_{\left\vert t' \right\vert h}
	\left(
	\vecf{x}_k - t' \mathrm{grad}\,g\vert_{\vecf{x}_k}
	\right) \right\vert}
\\ & \ge 
\frac{\left\vert t' \right\vert}{
	\left\Vert \vecf{x}_k \right\Vert_2 
	\left\Vert \mathrm{prox}_{\left\vert t' \right\vert h}
	\left(
	\vecf{x}_k - t' \mathrm{grad}\,g\vert_{\vecf{x}_k}
	\right) \right\Vert_2}
\quad(\mathrm{Cauchy–Schwarz})
\\ & \ge 
\frac{\left\vert t' \right\vert}{\left\Vert \vecf{x}_k \right\Vert_2 \left\Vert
	\vecf{x}_k - t' \mathrm{grad}\,g\vert_{\vecf{x}_k}
	\right\Vert_2}
\quad(\mathbf{Lemma}~\ref{lemma. additional properties for convex and absolute homogeneous h})
\\ & = 
\frac{\left\vert t' \right\vert}{
	\left\Vert \vecf{x}_k  \right\Vert_2
	\sqrt{
		\left\Vert \vecf{x}_k  \right\Vert_2^2 + 
		\left\Vert t' \mathrm{grad}\,g\vert_{\vecf{x}_k}  \right\Vert_2^2
	}	
}
\quad(\mathrm{since\ } \vecf{x}_k^{\tran} \mathrm{grad}\,g\vert_{\vecf{x}_k} = 0)
\\ & =  
\frac{1}{
	\sqrt{
		\left( 1 / t' \right)^2 + 
		\left\Vert  \mathrm{grad}\,g\vert_{\vecf{x}_k}  \right\Vert_2^2
	}	
}
\quad(\mathrm{since\ }\left\Vert \vecf{x}_k  \right\Vert_2 = 1)
.
\end{align*}

\section{Proof of Lemma \ref{lemma: the h cost for given retraction and homogeinity}}
\label{appendix: proof of lemma h cost for given retraction and homoginity}

Since $\langle \vecf{x}_k, \vecf{v}_k \rangle = 0$ as  $\vecf{v}_k \in \mathcal{T}_{\vecf{x}_k} \mathcal{S}$, the norm satisfies:
\begin{align*}
\left\Vert \vecf{x}_{k} + \vecf{v}_{k} \right\Vert_2
& =
\sqrt{
	\left\Vert \vecf{x}_{k} \right\Vert_2^2 + \left\Vert \vecf{v}_{k} \right\Vert_2^2
	+ 2 \langle \vecf{x}_k, \vecf{v}_k \rangle
} 
\\ & =
\sqrt{
	\left\Vert \vecf{x}_{k} \right\Vert_2^2 + \left\Vert \vecf{v}_{k} \right\Vert_2^2}
\ge
\sqrt{ \left\Vert \vecf{x}_{k} \right\Vert_2^2 }
= 1	.
\end{align*}
By retraction (\ref{eq: retraction of unit sphere}),
since $h(\cdot)$ is absolutely homogeneous,
we write:
\begin{align*}
h\left( \mathcal{R}_{\vecf{x}_{k}} \left( \vecf{v}_{k} \right) \right)
& =
h \left ( \frac{\vecf{x}_{k} + \vecf{v}_{k}}{ \left\Vert \vecf{x}_{k} + \vecf{v}_{k} \right\Vert_2 }
\right)
\\ & =
\frac{1}{ \left\Vert \vecf{x}_{k} + \vecf{v}_{k} \right\Vert_2 }
h\left( \vecf{x}_{k} + \vecf{v}_{k} \right)
\le
h\left( \vecf{x}_{k} + \vecf{v}_{k} \right)
.
\end{align*}

\section{Lipschitz-type Constant $L$}
\label{appendix: Lipschitz-type Constant of rayleigh quotient}

It is well-known in the Rayleigh quotient literature that:
\begin{equation*}
\max_{\vecf{v}} \frac{\vecf{v}^{\tran} \matf{A}  \vecf{v}}{ \vecf{v}^{\tran} \vecf{v} }
=
\sigma_{\max} (\matf{A})
,
\end{equation*}
where $\sigma_{\max} (\matf{A})$ is the largest singular value of $\matf{A}$.
Therefore for any $\vecf{v}$, we have:
\begin{equation}
\frac{\vecf{v}^{\tran} \matf{A}  \vecf{v}}{ \vecf{v}^{\tran} \vecf{v} }
\le
\sigma_{\max} (\matf{A})
\Leftrightarrow
\vecf{v}^{\tran} \matf{A}  \vecf{v}
\le 
\sigma_{\max} (\matf{A}) \left\Vert \vecf{v} \right\Vert_2^2
.
\end{equation}
From equation (\ref{eq: Riemannian gradient of ralyleigh quotient optimization}) and the fact that $\vecf{x}^{\tran} \vecf{v} = 0$, we have:
\begin{equation}
\langle
\mathrm{grad}\,g (\vecf{x}) ,\, \vecf{v}
\rangle
=
\langle
2 \matf{A} \vecf{x} - 2 (\vecf{x}^{\tran}  \matf{A} \vecf{x}) \vecf{x} ,\, \vecf{v}
\rangle
=
2 \vecf{x}^{\tran} \matf{A}  \vecf{v}
.
\end{equation}
Recall that
$\frac{1}{ \left\Vert \vecf{x} + \vecf{v} \right\Vert_2^2 }
=
\frac{1}{ 1 + \left\Vert \vecf{v} \right\Vert_2^2 }
\le 1
$.
For $g (\vecf{x}) = \vecf{x}^{\tran} \matf{A}  \vecf{x}$ and retraction (\ref{eq: retraction of unit sphere}),
we write
\begin{align*}
g \left( \mathcal{R}_{\vecf{x}} \left( \vecf{v} \right) \right)
& =
\left( \frac{\vecf{x} + \vecf{v}}{ \left\Vert \vecf{x} + \vecf{v} \right\Vert_2 } \right)^{\tran}
\matf{A}
\left( \frac{\vecf{x} + \vecf{v}}{ \left\Vert \vecf{x} + \vecf{v} \right\Vert_2 } \right)
\\ & =
\frac{1}{ \left\Vert \vecf{x} + \vecf{v} \right\Vert_2^2 }
\left(
\vecf{x}^{\tran} \matf{A}  \vecf{x}
+
2 \vecf{x}^{\tran} \matf{A}  \vecf{v}
+ 
\vecf{v}^{\tran} \matf{A}  \vecf{v}
\right)
\\ & \le
\vecf{x}^{\tran} \matf{A}  \vecf{x}
+
2 \vecf{x}^{\tran} \matf{A}  \vecf{v}
+ 
\vecf{v}^{\tran} \matf{A}  \vecf{v}
\\ & \le 
g (\vecf{x})
+
\langle
\mathrm{grad}\,g (\vecf{x}) ,\, \vecf{v}
\rangle
+
\frac{2 \sigma_{\max} (\matf{A})}{2}
\left\Vert \vecf{v} \right\Vert_2^2
.
\end{align*}
Therefore we obtain $L = 2 \sigma_{\max} (\matf{A})$.

\end{document}